\documentclass{amsart}

\usepackage{amsfonts,amscd,amsmath,latexsym}
\numberwithin{equation}{section}
\theoremstyle{plain}
\newtheorem{thm}{Theorem}[section]

\newtheorem{lemma}[thm]{Lemma}
\newtheorem{prop}[thm]{Proposition}
\newtheorem{cor}[thm]{Corollary}

\theoremstyle{remark}
\newtheorem{remark}[thm]{Remark}
\theoremstyle{definition}
\newtheorem{defn}[thm]{Definition}
\newtheorem{example}[thm]{Example}

\newtheorem{notation}[thm]{Notation}

\newcommand{\N}{\ensuremath{\mathbb{N}}}

\newcommand{\HOM}{\operatorname{\mathcal{H}\textit{om}}}
\newcommand{\EXT}{\operatorname{\underline{\mathcal{E}\textit{xt}}}}

\newcommand{\Lex}{\operatorname{Lex}}
\newcommand{\ann}{\operatorname{ann}}
\newcommand{\sh}[1]{\ensuremath{{\mathcal #1}}}

\newcommand{\Proj}{\operatorname{Proj}}
\newcommand{\Supp}{\operatorname{Supp}}

\newcommand{\Hom}{\operatorname{Hom}}
\renewcommand{\mod}{\operatorname{mod}}
\newcommand{\Mod}{\operatorname{Mod}}

\newcommand{\End}{\operatorname{End}}

\newcommand{\Id}{\operatorname{Id}}

\newcommand{\QCoh}{\operatorname{QCoh}}

\renewcommand{\O}{{\mathcal O}}

\newcommand{\Inj}{\operatorname{Inj}}
\newcommand{\ass}{\operatorname{ass}}
\newcommand{\Kdim}{\operatorname{Kdim}}
\newcommand{\M}{\mathcal{M}}
\newcommand{\lrk}{\operatorname{lrk}}
\newcommand{\rrk}{\operatorname{rrk}}

\newcommand{\vphi}{\varphi}
\newcommand{\Frob}{\operatorname{Frob}}
\newcommand{\Bimod}{\operatorname{Bimod}}
\newcommand{\BIMOD}{\operatorname{BIMOD}}
\newcommand{\A}{\mathcal{A}}
\newcommand{\B}{\mathcal{B}}

\newcommand{\Vect}{\operatorname{Vect}}

\newcommand{\E}{\mathcal{E}}
\newcommand{\I}{{\normalfont\bf I}}

\newcommand{\DimPres}{\operatorname{DimPres}}
\newcommand{\RLoc}{\operatorname{RLoc}}
\newcommand{\LLoc}{\operatorname{LLoc}}
\newcommand{\Loc}{\operatorname{Loc}}
\newcommand{\LCent}{\operatorname{LCent}}
\renewcommand{\a}{\alpha}
\renewcommand{\b}{\beta}
\newcommand{\qgr}{\operatorname{qgr}}
\renewcommand{\d}{\delta}
\newcommand{\Cent}{\operatorname{Cent}}
\newcommand{\opp}{\rm opp}
\newcommand{\END}{\operatorname{{\mathcal{E}\textit{nd}}}}

\begin{document}
\title[Frobenius bimodules between spaces]{Frobenius bimodules between noncommutative spaces}
\author{Christopher J. Pappacena}
\address{Department of Mathematics, Baylor University, Waco, TX 76798}
\email{Chris\_$\,$Pappacena@baylor.edu}

\thanks{This research was partially supported by the National Security Agency under grant NSA 032-5529.}

\begin{abstract} In this paper we study Frobenius bimodules between noncommutative spaces (quasi-schemes), developing some of their basic properties.  If $X$ and $Y$ are spaces, we study those Frobenius $X,Y$-bimodules $_X\M_Y$ satisfying properties that are natural in the context of noncommutative algebraic geometry, focusing in particular on cartain ``local" conditions on $\M$.  As applications, we prove decomposition and gluing theorems for those Frobenius bimodules which have good local properties. Additionally, when $X$ and $Y$ are schemes we relate Frobenius $X,Y$-bimodules  to the sheaf $X,Y$-bimodules introduced by Van den Bergh in \cite{Van Sklyanin}.
\end{abstract}

\maketitle

\tableofcontents

\section{Introduction}  
As is sometimes the case in mathematics, the main motivation for beginning work on this paper plays only a relatively small role in the finished product.  Our original (and somewhat immodest) goal was to develop a general definition of \emph{vector bundle} for the spaces studied in noncommutative algebraic geometry. Here we mean noncommutative algebraic geometry in the spirit of Rosenberg \cite{Rosenberg local, Rosenberg book}, Van den Bergh \cite{Van blowup}, and Smith \cite{Smith integral, Smith subspaces}.  Thus the basic geometric object is a Grothendieck category $\Mod X$, to be viewed as the category of sheaves on some (nonexistent) space $X$.  

Before discussing our approach to the issue of defining vector bundles on noncommutative spaces, it is worthwhile to recall some approaches that have been taken by other authors.  If $X$ and $Y$ are schemes, then Van den Bergh studies sheaf bimodules $_X\E_Y$ which are locally free of finite rank on each side \cite{Van Sklyanin, Van P1} and defines noncommutative symmetric algebras associated to certain locally free sheaf bimodules $_X\E_X$ in \cite{Van P1}.  In another direction, \cite{BGK} and \cite{KKO} consider certain noncommutative projective schemes $\Proj R$ (in the sense of \cite{Artin Zhang}), and $\E\in\qgr(R)$ is defined to be locally free if $\EXT^i(\E,\O)=0$ for all $i>0$. See \cite[Definition 1.1.4]{BGK} or \cite[Definition 5.4]{KKO} for more details.

While the approaches taken by the above authors are undoubtedly the correct ones for the problems that they are interested in, neither formulation of ``locally free sheaf" can be readily generalized to more general noncommutative spaces.  For instance, the sheaf bimodules of \cite{Van Sklyanin, Van P1} are actual sheaves over the fiber product of schemes $X\times Y$.  Similarly, the definition of ``locally free sheaf" given in \cite{BGK,KKO} applies only to noncommutative projective schemes $\Proj R$, where $R$ is a ``strongly regular algebra" \cite[Definition 1.1.1]{BGK}.

Our approach to the problem of defining vector bundles on noncommutative spaces is to seek a definition that is well-suited for algebraic $K$-theory. (Indeed, one of the long-term goals of the author is to study intersection theory for noncommutative spaces, and it is expected that the $K$-theory of vector bundles, once defined, will play an important role.)  In particular, vector bundles should be ``two-sided," in the sense that ``tensoring with a vector bundle $\E$" should define an exact functor on $\Mod X$, making the $K$-theory of $\Mod X$ into a module over the $K$-theory of $\Vect(X)$ (the category of vector bundles on $X$).  This two-sidedness would also make it possible to define the tensor product of two vector bundles, an important ingredient in such basic constructions as tensor algebras, symmetric algebras, and so on.

Fortunately, Ven den Bergh has developed a theory of bimodules between noncommutative spaces $X$ and $Y$ \cite{Van blowup} as a way of providing these ``two-sided" objects. Briefly, an $X,Y$-bimodule $_X\M_Y$ is a left exact functor $\HOM_Y(\M,-):\Mod Y\rightarrow\Mod X$ having a left adjoint.  (See section 2 below for more details.) In particular, ``tensoring with $\M$" means applying the left adjoint functor, and ``tensoring two bimodules" means composition of functors. It therefore seems natural to define a vector bundle on a noncommutative space $X$ to be an $X,X$-bimodule $_X\E_X$ satisfying certain properties.  In particular, if $X$ is a scheme and $\E$ is a locally free sheaf of finite rank on $X$, then the functor $\HOM_{\O_X}(\E,-)$ should enjoy all of these properties.   

The following functorial property is particularly intriguing:  If $\E$ is a locally free sheaf of finite rank on a scheme $X$, then the dual sheaf is defined as $\E^*=\HOM_{\O_X}(\E,\O_X)$, and there are functorial isomorphisms $\HOM_{\O_X}(\E,\sh{F})\cong\sh{F}\otimes_{\O_X}\E^*$ and $\sh{F}\otimes_{\O_X}\E\cong\HOM_{\O_X}(\E^*,\sh{F})$ for all quasicoherent $\O_X$-modules $\sh{F}$.  Thus the functor $\HOM_{\O_X}(\E,-)$ is both a left and a right adjoint to $-\otimes_X\E$.  Functors $F:\sh{A}\rightarrow \sh{B}$ and $G:\sh{B}\rightarrow\sh{A}$ betweeen categories $\A$ and $\B$ such that $(F,G)$ and $(G,F)$ are both adjoint pairs are called \emph{Frobenius pairs}, and $F$ and $G$ are individually called \emph{Frobenius functors} \cite{CGN}.  Thus if $\E$ is a locally free sheaf of finite rank on a scheme $X$, then $-\otimes_{\O_X}\E$ is a Frobenius functor.  
  
Thus, we were initially led to consider vector bundles on noncommutative spaces to be $X,X$-bimodules $_X\E_X$ such that $\HOM_X(\E,-)$ is a Frobenius functor, perhaps with additional properties as well.  During the course of investigating this idea, it seemed clear that a more natural setting was to study the property of bimodules $_X\M_Y$ between two noncommutative spaces $X$ and $Y$ with $\HOM_Y(\M,-)$ a Frobenius functor, and this paper is devoted to the study of such bimodules.  Many of the results on Frobenius $X,Y$-bimdodules proved in this paper were inspired by a desire to establish desirable properties of vector bundles (e.g. gluing theorems, local analysis, behavior of rank functions, etc.).  

Frobenius functors have been studied recently by a number of authors, often in the context of Frobenius functors between module and comodule categories. (See for instance \cite{CIM, CK, CMZ, CGN, K1, K2, Kadison, MM} and the references contained therein.)  What is new in this paper is the study of Frobenius functors in a geometric context, most notably in the consideration of local behavior of Frobenius functors and Frobenius functors between schemes.  

We now describe the contents of the paper in more detail.  In section 2, we summarize some of the salient definitions and results of noncommutative algebraic geometry that we use in subsequent sections, including a discussion of weakly closed and weakly open subspaces in the sense of Smith \cite{Smith subspaces}, local spaces as defined by Rosenberg \cite{Rosenberg local} and their generalizations, Van den Bergh's notion of a bimodule between noncommutative spaces \cite{Van blowup}, and the injective spectrum \cite{Pappacena injective}.  In section 3, we formally introduce the notion of a Frobenius bimodule between noncommutative spaces and prove some of the basic properties they possess.  

Section 4 begins our study of Frobenius bimodules with additional properties, by considering those bimodules which preserve given dimension functions on $\Mod X$ and $\Mod Y$.  The main result is a description of how dimension preserving bimodules act on the indecomposable injectives and consequences of this behavior (Theorem \ref{decomp thm} and its corollaries).  In section 5, we study so-called ``right localizing" Frobenius bimodules,  which are a class of Frobenius bimodules that exhibit good local properties.  We show that under suitable hypotheses, right localizing $X,Y$-bimodules are precisely the Frobenius bimodules that ``come from" geometric data associated to $X$ and $Y$ (Theorem \ref{continuous thm}). In section 6, we consider Frobenius bimodules over noetherian schemes, relating them to the sheaf bimodules studied in \cite{Nyman thesis, Nyman, Van Sklyanin, Van P1}.  In particular, we show that a sheaf $X,Y$-bimodule which is locally free of finite rank on each side, such that the left and right duals agree, gives rise to a Frobenius bimodule in our sense (Proposition \ref{sheaf Frobenius prop}). Conversely, a right localizing Frobenius bimodule between schemes $X$ and $Y$ is given by tensoring by a suitable locally free sheaf $X,Y$-bimodule (Theorem \ref{right localizing scheme thm}). 

We study rank functions for Frobenius bimodules in section 7, in particular rank functions for right localizing bimodules.  Given a right localizing bimodule $_X\M_Y$, we show how to decompose $X$ and $Y$ into the disjoint union of weakly open subspaces such that the the restriction of $\M$ to these subspaces has constant rank (Theorem \ref{constant rank thm}).  Section 8 is devoted to proving a gluing theorem for certain right localizing bimodules. Under suitable hypotheses, we show that right localizing bimodules on weakly open covers of spaces which agree on the overlaps can be patched to give a right localizing bimodule between the original spaces (Theorem \ref{gluing thm}). In section 9 we study Frobenius $X,X$-bimodules, looking specifically at those bimodules which can be localized to give Frobenius $U,U$-bimodules for every weakly open subspace $U$ of $X$ and those bimodules which are ``commutative," in the sense that they commute with the center of $X$.  

In section 10 we study the category $\Frob(X,Y)$ of Frobenius $X,Y$-bimodules.  The main result in this section is a duality between $\Frob(X,Y)$ and $\Frob(Y,X)$ (Theorem \ref{duality thm}). Additionally, we consider categorical properties of the various special classes of Frobenius bimodules discussed earlier in the paper.  In section 11, we discuss in more depth the question of how to define a vector bundle on a noncommutative space $X$.  We come short of actually proposing a definition, instead contenting ourselves with discussing the pros and cons of three natural candidates for such a definition.  Finally, we collect some basic definitions and results on abelian categories that are used in the body of the paper in an appendix.

\subsubsection*{Acknowledgments} The author thanks David Arnold and Mark Sepanski for helpful feedback on this paper, and Adam Nyman for conversations on the ideas of this paper, especially for his help in understanding the material on sheaf bimodules in \cite{Nyman thesis, Nyman, Van Sklyanin, Van P1}, and for finding and helping to correct an error in the proof of Proposition \ref{sheaf Frobenius prop}(2).

\section{Noncommutative algebraic geometry} 
In this section we recall some of the basic machinery of noncommutative geometry as developed by Van den Bergh \cite{Van blowup}, Rosenberg \cite{Rosenberg local, Rosenberg book}, and Smith \cite{Smith integral, Smith subspaces}.  The terminology that we adopt is taken for the most part from Smith's papers \cite{Smith integral, Smith subspaces}.  Any or all of the above papers can and should be consulted for further information or details.

\subsection{Noncommutative spaces}
A \emph{noncommutative space} (or simply a \emph{space}) $X$ is a Grothendieck category $\Mod X$.  (Noncommutative spaces are called \emph{quasi-schemes} in \cite{Van blowup}.)  Recall that this means that $\Mod X$ is an Ab5 abelian category with a generator.  The objects of $\Mod X$ are called \emph{$X$-modules}. We write $\mod X$ for the full subcategory of $\Mod X$ consisting of noetherian $X$-modules, and call $X$ \emph{noetherian} if $\Mod X$ is a locally noetherian category; that is, if every $X$-module is the direct limit of its noetherian submodules.   An \emph{enriched space} is a pair $(X,\O_X)$, where $X$ is a space and $\O_X\in\Mod X$.  The distinguished $X$-module $\O_X$ is called the \emph{structure module}.

The following are three canonical examples of noncommutative spaces.
\begin{enumerate}
\item $\QCoh(\O_X)$, the category of quasicoherent sheaves on a quasicompact, quasiseparated scheme $X$ \cite[Appendix B.3]{Thomason Trobaugh}.  For consistency, we assume throughout the rest of this paper that all schemes are quasicompact and quasiseparated.
\item $\Mod R$, the category of right modules over a ring $R$.
\item $\Proj R$, the noncommutative projective scheme associated to an $\N$-graded $k$-algebra $R$ \cite{Artin Zhang, Verevkin}.  
\end{enumerate}
Each of the above spaces can be enriched in a natural way.

Being a Grothendieck category, $\Mod X$ is complete and cocomplete, and has injective hulls.  If $M$ is an $X$-module, then we write $E(M)$ for the injective hull of $M$.  If $\kappa$ is a cardinal, then we write $M^{(\kappa)}$ and $M^{\kappa}$ for the direct sum and direct product of $\kappa$ copies of $M$, respectively.

In \cite[Definitions 2.4 and 2.5]{Smith subspaces}, Smith gives definitions of weakly closed, closed, and weakly open subspaces of a space $X$.  We refer the reader to \cite{Smith subspaces} for most of the details, but we recall a few of the definitions and results that we will use most frequently in the paper.  

A \emph{weakly closed subspace} $Z$ of $X$ is a full subcategory $\Mod Z$ of $\Mod X$, which is closed under isomorphism, direct sums, and subquotients. A weakly closed subspace $Z$ is \emph{closed} if $\Mod Z$ is also closed under products.  Typically one writes $i:Z\rightarrow X$ to denote the inclusion of a weakly closed subspace, where $i_*:\Mod Z\rightarrow \Mod X$ denotes the inclusion functor.  A \emph{weakly open subspace} $U$ of $X$ is a full subcategory $\Mod U$ of $\Mod X$, closed under isomorphism and kernels, such that the inclusion functor $j_*:\Mod U\rightarrow\Mod X$ has an exact left adjoint $j^*$. The functor $j^*$ is referred to as \emph{restriction to $U$}, and if $M\in\Mod X$ we also write $M|_U$ for $j^*M$.  We will frequently use the fact that $j^*j_*\simeq\Id_U$.

There is a connection between weakly open and weakly closed subspaces of $X$, which we now describe. If $j:U\rightarrow X$ is the inclusion of a weakly open subspace of $X$, then $T(U)=\{M:M|_U=0\}$ is a localizing subcategory of $\Mod X$, and $\Mod U\simeq \Mod X/T(U)$ \cite[proof of Proposition 6.6]{Smith subspaces}.  Moreover, every weakly open subspace arises in this fashion:  If $T$ is a localizing subcategory of $\Mod X$, then there is a unique weakly open subspace $U$ of $X$ with $\Mod U\simeq \Mod X/T$ \cite[Propositions 6.5 and 6.6]{Smith subspaces}.  If $Z$ is a weakly closed subspace of $X$, let $\Mod_ZX$ denote the smallest localizing subcategory containing $Z$ (the category of $X$-modules \emph{supported} at $Z$).  Then the weakly open subspace $U$ with $\Mod U\simeq \Mod X/\Mod_ZX$ is the \emph{complement} to $Z$, written $X\setminus Z$ \cite[Definition 6.4]{Smith subspaces}.  If $U=X\setminus Z$ with $Z$ a closed subspace of $X$, then $U$ is called an \emph{open} subspace of $X$.  Thus every weakly open subspace is of the form $X\setminus Z$, and $Z$ is determined up to $\Mod_ZX$. Note that if $U$ and $V$ are weakly open subspaces of $X$ with $\Mod U\simeq \Mod X/T$ and $\Mod V\simeq \Mod X/S$ for localizing subcategories $S$ and $T$, then $U\subseteq V$ if and only if $S\subseteq T$.  

We recall a few more facts about weakly open subspaces:  If $\{U_i:i\in I\}$ is a collection of weakly open subspaces, then their \emph{union} is defined as follows. Write $\Mod U_i\simeq \Mod X/T_i$ for localizing subcategories $T_i$. Then $\bigcup_{i\in I}U_i$ is defined to be the weakly open subspace with $\Mod \bigcup_{i\in I}U_i\simeq \Mod X/\bigcap_{i\in I} T_i$ \cite[Definition 6.9]{Smith subspaces}. Since the intersection of localizing subcategories is localizing this definition makes sense. If $\bigcup_{i\in I}U_i=X$, then we say that the $U_i$ form a \emph{weakly open cover} of $X$ \cite[Definition 6.9]{Smith subspaces}.   Finally, if $U_1,U_2$ are weakly open subspaces with $U_i=X\setminus Z_i$ for $i=1,2$, we set $U_1\cap U_2=X\setminus Z_1\bullet Z_2$, where $\bullet$ is the Gabriel product \cite[Definition 6.14]{Smith subspaces}.  Since $\Mod_{Z_1\bullet Z_2}X=\Mod_{Z_2\bullet Z_1}X$, we have $U_1\cap U_2=U_2\cap U_1$.  

\subsection{Local and semilocal spaces} The notion of a local space was introduced by Rosenberg in \cite[Definition 3.1.1]{Rosenberg local}. An $X$-module $Q$ is called  \emph{quasifinal} if, given any $M\in\Mod X$, $Q$ is finitely subgenerated by $M$.  

If $Q$ is a quasifinal $X$-module and $S$ is a simple $X$-module, then $Q$ is isomorphic to a finite direct sum of copies of $S$.  Thus we may take $Q=S$, and we see additionally that $S$ must be the unique (up to isomorphism) simple $X$-module.

\begin{defn}  A space $X$ is \emph{local} if there exists a quasifinal $X$-module (which we always take to be simple).\end{defn}

\begin{lemma}  If $X$ is local with simple module $S$, then $E(S)$ is an injective cogenerator for $\Mod X$.\label{injective cogenerator lemma}\end{lemma}

\begin{proof} Let $M\in\Mod X$ and let $\sigma[M]$ denote the category subgenerated by $M$.  That is, $\sigma[M]$ is the smallest weakly closed subspace of $X$ which contains $M$, and every $X$-module in $\sigma[M]$ is isomorphic to a submodule of a quotient of a direct sum of copies of $M$.  Let $J$ denote the injective hull of $S$ in $\sigma[M]$; then 
 there are exact sequences $0\rightarrow A\rightarrow\oplus_{\a\in I} M\rightarrow B\rightarrow 0$ and $0\rightarrow J\rightarrow B\rightarrow C\rightarrow 0$.  Since $J$ is injective in $\sigma[M]$, the second sequence splits, and so there is an epic $\oplus_{\alpha\in I}M\rightarrow J$. So, for some index $\alpha$ the composition $M\xrightarrow{i_\alpha}\oplus_{\alpha\in I}M\rightarrow J$ is nonzero; i.e. $\Hom_X(M,J)\neq 0$.  Since $E(S)$ is also the injective hull of $J$ in $\Mod X$, we see that $\Hom_X(M,E(S))\neq 0$ for all $M\in \Mod X$. Thus $E(S)$ is an injective cogenerator.
\end{proof} 

It turns out that this property characterizes local spaces.

\begin{lemma} Let $X$ be a space with an indecomposable injective cogenerator $E$.  Then $X$ is local. \label{local lemma} 
\end{lemma}

\begin{proof}  If $S$ is a simple $X$-module, then the fact that $\Hom_X(S,E)\neq 0$ shows that $E\cong E(S)$, and that $S$ is the unique simple $X$-module up to isomorphism.  Now, given any $M\in\Mod X$, the image of $M$ under a nonzero $f\in\Hom_X(M,E)$ must contain $S$ as a submodule.  Thus $S$ is finitely subgenerated by each $X$-module $M$, showing that $S$ is quasifinal.
\end{proof}

We extend Rosenberg's definition of a local space to a semilocal space in the obvious way: We call a set $\{Q_i:i\in I\}$ of $X$-modules a \emph{quasifinal set} if, for all $M\in\Mod X$, there exists an index $i$ such that $Q_i$ is finitely subgenerated by $M$.  A space $X$ is \emph{semilocal} if it has a finite quasifinal set $\{Q_1,\dots, Q_n\}$. As above, we may assume that our quasifinal set is $\{S_1,\dots, S_n\}$, where the $S_i$ are a complete set of representatives for the isomorphism classes of the simple $X$-modules. Then we have the following generalization of the above lemmas.

\begin{lemma} A space $X$ is semilocal if and only if  there exists a finite set $\Sigma=\{S_1,\dots, S_t\}$ of simple $X$-modules such that $E=\oplus_{i=1}^tE(S_i)$ is an injective cogenerator for $\Mod X$.\label{semilocal lemma}
\end{lemma}

\begin{proof}  Suppose first that $X$ is semilocal. We take $\Sigma$ to be the set of all simple $X$-modules (up to isomorphism).  Given $M\in\Mod X$, there exists an index $i$ such that $S_i\in\sigma[M]$.  As in the proof of \ref{injective cogenerator lemma}, we have $\Hom_X(M,E(S_i))\neq 0$.  Thus, given any $M\in\Mod X$, we have $\Hom_X(M,E)\neq 0$, so that $E$ is an injective cogenerator for $\Mod X$.

Conversely, suppose that $E$ is an injective cogenerator.  Then as in the proof of Lemma \ref{local lemma}, we see that $\Sigma$ is a complete set of representatives for the isomorphism classes of simple $X$-modules.  Now, given $M\in\Mod X$, the fact that $\Hom_X(M,E)\neq 0$ implies that $\Hom_X(M,E(S_i))\neq 0$ for some $i$.  Thus the image of a nonzero $f\in\Hom_X(M,E(S_i))$ contains $S_i$ as a submodule, showing that $S_i$ is finitely subgenerated by $M$.  Hence $\Sigma$ is a finite quasifinal set, and $X$ is semilocal.
\end{proof}

\subsection{The injective spectrum}  
For a noncommutative space $X$, let $\Inj(X)$ denote the set of isomorphism classes of indecomposable injective $X$-modules.  We call $\Inj(X)$ the \emph{injective spectrum} of $X$, and will frequently think of it as an underlying ``point set" for $X$.  To reinforce this idea we shall write elements of $\Inj(X)$ using lowercase letters: $x$, $y$, etc.  We fix a representative for each isomorphism class in $\Inj(X)$, and denote the representative for $x$ by $E(x)$.  

If $E$ is an injective $X$-module, then $T(E)=\{M\in\Mod X:\Hom_X(M,E)=0\}$ is a localizing subcategory of $\Mod X$. Moreover, any localizing subcategory arises in this way;  indeed, given a localizing subcategory $T$ of $\Mod X$, let $\Sigma$ be a set of representatives for the $T$-torsionfree injective $X$-modules.  (That is, every $T$-torsionfree injective $X$-module is isomorphic to exactly one module in $\Sigma$.)  Then one checks readily that $T=T\bigl(\prod_{E\in\Sigma}E\bigr)$. If $x\in\Inj(X)$ we write  $T_x$ in place of $T(E(x))$. 

\begin{lemma} Let $X$ be a noetherian space, and let $E$ be an injective $X$-module. Let $\Sigma=\{x\in\Inj(X):\mbox{$E(x)$ is isomorphic to a summand of $E$}\}.$  Then $T(E)=\bigcap_{x\in\Sigma}T_x$.\label{localizing subcategory lemma}
\end{lemma}

\begin{proof} Given $x\in\Sigma$, we easily see that $\Hom_X(M,E)=0$ implies $\Hom_X(M,E(x))=0$, so that $T(E)\subseteq \bigcap_{x\in \Sigma}T_x$.  Conversely, suppose that  $\Hom_X(M,E)\neq 0$. Then $E(f(M))$ is isomorphic to a summand of $E$, where $f:M\rightarrow E$ is a nonzero morphism.   Since $X$ is noetherian $f(M)$ has a noetherian submodule, and hence $E(f(M))$ has an indecomposable summand. If $E(x)$ is isomorphic to an indecomposable summand of $E(f(M))$, then $x\in\Sigma$ and $\Hom_X(M,E(x))\neq 0$. 
\end{proof}

We write $X_x$ for the weakly open subspace of $\Mod X$ with $\Mod X_x\simeq \Mod X/T_x$, and we denote the inclusion by $j_x:X_x\rightarrow X$. Given $M\in\Mod X$, we write $M_x$ for $j^*_xM$ and call $M_x$ the \emph{stalk of $M$ at $x$}.  Then we have the following:

\begin{lemma} $X_x$ is a local space, with indecomposable injective cogenerator $E(x)_x$.
\end{lemma}

\begin{proof}  Since $E(x)$ is $T_x$-torsionfree, we see that $E(x)_x$ is an injective $X_x$-module, and since $j_{x*}E(x)_x\cong E(x)$, we see that $E(x)_x$ is indecomposable.  Finally, given $X_x$-module $N$, there is an $M\in\Mod X$ with $N=M_x$.  Then $\Hom_{X_x}(M_x,E(x)_x)=0$ if and only if $\Hom_X(M,E(x))=0$, if and only if $M\in T_x$, if and only if $M_x=0$.  So $E(x)_x$ is an indecomposable injective cogenerator for $\Mod X_x$, whence $X_x$ is local by Lemma \ref{local lemma}.
\end{proof}

The idea of studying $\Inj(X)$ (in the context of more general abelian categories) goes back to Gabriel \cite{Gabriel}, and in \cite{Gabriel} he introduced a natural topology on $\Inj(X)$ as follows.  If $M\in\mod X$, let $V(M)=\{x\in\Inj(X):\Hom_X(M,E(x))\neq 0\}$.  Then the \emph{Gabriel topology} on $\Inj(X)$ is obtained by taking $\{V(M):M\in\mod X\}$ as a basis for the closed sets. 

Let $U$ be a weakly open subspace of $\Mod X$, with $\Mod U\simeq\Mod X/T$.  Then there is a bijection between $\Inj(U)$ and $\{x\in\Inj(X):\mbox{$E(x)$ is $T$-torsionfree}\}$. Explicitly, given $u\in\Inj(U)$, there exists $x\in\Inj(X)$ with $E(u)\cong j^*E(x)$.  If $X$ is noetherian then this bijection is a homeomorphism:  If $N\in\mod U$, then there exists $M\in\mod X$ with $N=j^*M$.  Then $u\in V(N)$ if and only if $\Hom_U(N,E(u))\neq 0$, if and only if $\Hom_U(j^*M,j^*E(x))\neq 0$, if and only if $\Hom_X(M,j_*j^*E(x))\neq 0$, if and only if $\Hom_X(M,E(x))\neq 0$, if and only if $x\in V(M)$.  (Here we have used that $j_*j^*E(x)\cong E(x)$ whenever $E(x)$ is $T$-torsionfree.)

\begin{lemma} Let $X$ be a noetherian space, let $M$ be a noetherian $X$-module, and let $T(M)$ be the smallest localizing subcategory of $\Mod X$ which contains $M$.  Let $\mathfrak{U}$ be the basic open subset $V(M)^c$ of $\Inj(X)$.
\begin{enumerate}
\item $T(M)=\Mod_{\sigma[M]}X$, where $\sigma[M]$ is the category subgenerated by $M$.  
\item $T(M)=\bigcap_{x\in\mathfrak{U}}T_x$.
\item $\mathfrak{U}$ is homeomorphic to $\Inj(U)$, where $U$ is the weakly open subspace of $X$ with $\Mod U\simeq \Mod X/T(M)$.
\end{enumerate}\label{Gabriel topology lemma}
\end{lemma}

\begin{proof} (1) is immediate from the definitions.

(2) Let $E=\prod_{x\in\mathfrak{U}}E(x)$, so that $T(E)=\bigcap_{x\in\mathfrak{U}}T_x$.  Since $\Hom_X(M,E)=\prod_{x\in\mathfrak{U}}\Hom_X(M,E(x))=0$, we see that $M\in T(E)$ and hence $T(M)\subseteq T(E)$.  Conversely, suppose that $M\in T$ for some localizing subcategory $T$.  Then $T=\bigcap_{x\in\Sigma}T_x$ for some $\Sigma\subseteq \Inj(X)$ by Lemma \ref{localizing subcategory lemma}.  Since $M\in T_x$, we have $\Hom_X(M,E(x))=0$ for all $x\in\Sigma$; that is, $\Sigma\subseteq \mathfrak{U}$.  It follows that $T(E)\subseteq T$, so that $T(E)= T(M)$.

(3) Let $j:U\rightarrow X$ denote the inclusion.  Since $\Mod U\simeq \Mod X/T(M)$, there is a homeomorphism between $\Inj(U)$ and $\Sigma=\{x\in\Inj(X):\mbox{$E(x)$ is $T(M)$-torsionfree}\}$.  If $x\in V(M)$, then $\Hom_X(M,E(x))\neq 0$.  If $f$ is a nonzero morphism, then $f(M)$ is a nonzero submodule of $E(x)$ in $T(M)$, so that $E(x)$ is not $T(M)$-torsionfree.  This shows that $\Sigma\subseteq\mathfrak{U}$.  Conversely, let $x\in\mathfrak{U}$. If $E(x)$ is not $T(M)$-torsionfree, then $E(x)$ has a nonzero submodule $N$ with $N\in\sigma[M]$.  If we choose $N$ to be injective in $\sigma[M]$, then as in the proof of Lemma \ref{injective cogenerator lemma} we have $\Hom_X(M,N)\neq 0$.  This in turn implies that $\Hom_X(M,E(x))\neq 0$, contradicting the fact that $x\in\mathfrak{U}$. Thus $\Sigma=\mathfrak{U}$.
\end{proof}

If $U$ is a weakly open subspace of $X$, we will frequently treat $\Inj(U)$ as a subset of $\Inj(X)$ under the homeomorphism described above.

If $(X,\O_X)$ is an enriched space, then we define a sheaf of (not necessarily commutative) rings $\END(\O_X)$ on $\Inj(X)$ as follows.  Given a basic open subset $\mathfrak{U}$ of $\Inj(X)$, we let $U$ be the weakly open subspace associated to $\mathfrak{U}$ as in part (3) of the above lemma.   The rule $\mathfrak{U}\mapsto\End_U(\O_U)$ associates a ring to each basic open subset of $\Inj(X)$, where by definition $\O_U=\O_X|_U$.  If $\mathfrak{V}\subseteq\mathfrak{U}$ are basic open sets with associated weakly open subspaces $V$ and $U$ respectively, then it is easy to check that $V\subseteq U$.  If $j:\Mod V\rightarrow \Mod U$ denotes the inclusion, then there is a functorial isomorphism $j^*(\O_U)\cong \O_V$. Thus we have restriction homomorphisms $\rho_\mathfrak{V}^\mathfrak{U}:\End_U(\O_U)\rightarrow\End_V(j^*\O_U)\cong \End_V(\O_V)$ for $\mathfrak{V}\subseteq\mathfrak{U}$ basic open.  These data on basic open sets can then be used to construct a sheaf $\END(\O_X)$ on $\Inj(X)$ in a natural way.  It is easy to check that if $x\in\Inj(X)$, then $\END(\O_X)_x\cong\End_{X_x}(\O_{X,x})$.

\begin{remark} The set $\Inj(X)$ is studied in some detail in \cite{Pappacena injective}, but we caution the reader that ``weak Zariski topology" \cite[Definition 4.6]{Pappacena injective} is \emph{not} the same as the Gabriel topology, as claimed in \cite{Pappacena injective}.  (The former can be defined by taking $\{V(M):M\in\Mod X\}$ as a basis for the closed sets on $\Inj(X)$.)
\end{remark}

\subsection{Bimodules}  In \cite{Van blowup}, Van den Bergh defines the notion of a \emph{bimodule} between noncommutative spaces as follows.  If $X$ and $Y$ are spaces, then $\Lex(Y,X)$ denotes the category of additive left exact functors from $\Mod Y$ to $\Mod X$.  The category of \emph{weak $X,Y$-bimodules} is defined to be $\BIMOD(X,Y)=\Lex(Y,X)^{\opp}$.  If $\sh{M}$ is an object of $\BIMOD(X,Y)$, then the underlying left exact functor associated to $\sh{M}$ is denoted by $\HOM_Y(\sh{M},-)$.  If $\HOM_Y(\sh{M},-)$ has a left adjoint, then $\sh{M}$ is called an \emph{$X,Y$-bimodule}.  The full subcategory of $\BIMOD(X)$ whose objects are $X,Y$-bimodules will be denoted by $\Bimod(X,Y)$. If $X=Y$ then we will write $\BIMOD(X)$ for $\BIMOD(X,X)$, and similarly for $\Bimod(X)$.  If we write $_X\sh{M}_Y$ then it is understood that $\sh{M}$ is a weak $X,Y$-bimodule.
  
If $\sh{M}\in\Bimod(X,Y)$, then we fix a left adjoint to $\HOM_Y(\sh{M},-)$ and denote it by 
$-\otimes_X\sh{M}$ (in general $-\otimes_X\sh{M}$ is only determined up to natural equivalence). Given weak bimodules $_X\sh{M}_Y$ and $_Y\sh{N}_Z$, then we denote the composition of the underlying functors by $\sh{M}\otimes_Y\sh{N}$. Then $\sh{M}\otimes_Y\sh{N}$ is  weak $X,Z$-bimodule, and we have several nice formulas, such as 
\begin{equation}\HOM_Z(\sh{M}\otimes_Y\sh{N},-)=\HOM_Y(\sh{M},\HOM_Z(\sh{N},-))\end{equation}
and, when $\M$ and $\sh{N}$ are bimodules, 
\begin{equation}-\otimes_X(\sh{M}\otimes_Y\sh{N})\simeq(-\otimes_X\sh{M})\otimes_Y\sh{N}.\end{equation}
Here the $-\otimes_Y\sh{N}$ on the left is composition, and the one on the right is the left adjoint to $\HOM_Y(\sh{N},-)$.

\section{Frobenius bimodules}
Let $R$ and $S$ be rings, and let $_RM_S$ be an $R,S$-bimodule. Recall that the left and right duals of $M$ are defined by $(_RM)^*=\Hom_R(M,R)$ and $(M_S)^*=\Hom_S(M,S)$, respectively.  Note that each of $(_RM)^*$ and $(M_S)^*$ is naturally an $S,R$-bimodule.  A bimodule $_RM_S$ is called a \emph{Frobenius bimodule} if each of $_RM$ and $M_S$ is finitely generated projective, and there is an isomorphism of $S,R$-bimodules $(_RM)^*\cong(M_S)^*$ \cite[Definition 2.1]{Kadison}.  The isomorphism $(_RM)^*\cong(M_S)^*$  allows us to speak unambiguously of the \emph{dual} $M^*$ of $M$.

The following result gives a functorial description of Frobenius bimodules.

\begin{prop} Let $R$ and $S$ be rings.
\begin{enumerate}
\item {\cite[Proposition 2.4]{Kadison}} If $_RM_S$ is a Frobenius bimodule, then there are equivalences of functors $-\otimes_RM\simeq \Hom_R(M^*,-)$ and $\Hom_S(M,-)\simeq-\otimes_SM^*$.  Thus $-\otimes_RM$ is both a left and a right adjoint to $\Hom_S(M,-)$.
\item \cite[Theorem 2.1]{CGN} If $F:\Mod R\rightarrow \Mod S$ and $G:\Mod S\rightarrow \Mod R$ are additive functors such that each of $(F,G)$ and $(G,F)$ is an adjoint pair, then there exists a Frobenius bimodule $_RM_S$ such that $F\simeq-\otimes_RM$ and $G\simeq\Hom_S(M,-)$.
\end{enumerate}
\label{Frobenius bimodule prop}
\end{prop}

In light of Proposition \ref{Frobenius bimodule prop}, if $\A$ and $\B$ are abelian categories and $F:\A\rightarrow \B$ and $G:\B\rightarrow\A$ are functors, it is natural to call $F$ a \emph {Frobenius functor}, and $(F,G)$ a \emph{Frobenius pair}, if $G$ is both a left and right adjoint to $F$ \cite[Definitions 1.1 and 1.2]{CGN}.  Restricting to the case of noncommutative spaces and using Van den Bergh's bimodule notation, we are led immediately to the following definition:

\begin{defn} Let $X$ and $Y$ be spaces. An $X,Y$-bimodule $_X\M_Y$ is called a \emph{Frobenius bimodule} if $(-\otimes_X\M,\HOM_Y(\M,-))$ is a Frobenius pair.  
\end{defn}

We view the category $\Frob(X,Y)$ of Frobenius $X,Y$-bimodules as a full subcategory of $\BIMOD(X,Y)$, so that two Frobenius $X,Y$-bimodules $\M$ and $\sh{N}$ are isomorphic if and only if there is an equivalence of functors $\HOM_Y(\M,-)\simeq\HOM_Y(\sh{N},-)$ (equivalently, if there is an equivalence $-\otimes_X\M\simeq-\otimes_X\sh{N}$).  We discuss the categorical properties of $\Frob(X,Y)$ in section 10.

Since we will frequently make computations with the underlying functors associated to a Frobenius bimodule $_X\M_Y$, it is advantageous to introduce some abbreviated notation.  Thus we will often use the phrase ``let $_X\M_Y=(F,G)$ be a Frobenius bimodule" to mean that $F=-\otimes_X\M$ and $G=\HOM_Y(\M,-)$.  If $_X\M_Y=(F,G)$ is a Frobenius $X,Y$-bimodule, then its \emph{dual} is defined to be the Frobenius $Y,X$-bimodule $_Y\M^*_X=(G,F)$.  Thus we tautologically have the formulas $\HOM_Y(\M,-)=-\otimes_X\M^*$, $-\otimes_X\M=\HOM_Y(\M^*,-)$, and $\M^{**}=\M$. 

If $_X\M_Y=(F,G)$ is a Frobenius $X,Y$-bimodule, then the functors $F$ and $G$ enjoy many useful properties.  We record some of the more basic ones here, in a lemma that will be freely used in the sequel without explicit comment.  The reader is referred to the appendix for any unfamiliar terms.  

\begin{lemma}  Let $_X\M_Y=(F,G)$ be a Frobenius bimodule.  Then each of $F$ and $G$ commutes with direct and inverse limits, preserves finitely presented, finitely generated, finitely copresented, and finitely cogenerated modules, and projective and injective modules.
\end{lemma}

\begin{prop} Let $_X\M_Y$ be a Frobenius bimodule.
\begin{enumerate}
\item $\M$ is nonzero if and only if $\HOM_Y(\M,E)\neq 0$ for some (hence every) injective cogenerator $E$ for $\Mod Y$, if and only if $E'\otimes_X\M\neq 0$ for some (hence every) injective cogenerator $E'$ for $\Mod X$.
\item $-\otimes_X\M$ is faithful if and only if $\HOM_Y(\M,E)$ is a cogenerator for $\Mod X$ for some (hence every) injective cogenerator $E$ for $\Mod Y$.
\end{enumerate}
\label{injective cogenerator prop}
\end{prop}

\begin{proof}  (1) Suppose that $\M$ is nonzero.  Then there exists some $M\in\Mod X$ with $M\otimes_X\M\neq 0$.  This then implies that $\Hom_Y(M\otimes_X\M,E)\neq 0$ for any injective cogenerator $E$ for $\Mod Y$, and the adjoint isomorphism gives $\Hom_X(M,\HOM_Y(\M,E))\neq 0$. Conversely, if $M=\HOM_X(\M,E)\neq 0$, then clearly $\M\neq 0$.  The last statement follows by considering the dual bimodule $\M^*$.

(2) Suppose that  $-\otimes_X\M$ is faithful, so that $M\otimes_X\M\neq 0$ whenever $M\neq 0$.  Since $\Hom_Y(M\otimes_X\M,E)\neq 0$ for all nonzero $M$, we have $\Hom_X(M,\HOM_Y(\M,E))\neq 0$ for all nonzero $M$.  Since $\HOM_Y(\M,E)$ is injective, this says that $\HOM_Y(\M,E)$ is an injective cogenerator.  For the converse, the fact that $\Hom_X(M,\HOM_Y(\M,E))\neq 0$ for all nonzero $M$ shows that $\Hom_Y(M\otimes_X\M,E)\neq 0$ for all nonzero $M$.  Thus $-\otimes_X\M$ takes nonzero $X$-modules to nonzero $Y$-modules and hence is faithful.
\end{proof}

Suppose that $_X\M_Y=(F,G)$ is a Frobenius bimodule, and let $Z$ be a weakly closed subspace of $Y$.  Then we define $F^{-1}Z$  by declaring $\Mod F^{-1}Z=\{M\in\Mod X:F(M)\in\Mod Z\}$.  

\begin{lemma} Let $_X\M_Y=(F,G)$ be a Frobenius bimodule. 
\begin{enumerate}
\item If $Z$ is a weakly closed (respectively closed) subspace of $Y$, then $F^{-1}Z$ is a weakly closed (respectively closed) subspace of $X$.
\item If $T$ is a localizing subcategory of $\Mod Y$, then $F^{-1}T$ is a localizing subcategory of $\Mod X$.
\item $F^{-1}(\bigcap_{i\in I}Z_i)=\bigcap_{i\in I}F^{-1}Z_i$.
\end{enumerate}
\end{lemma}

\begin{proof}  Parts (1) and (2) follow from the fact that $F$ is exact and commutes with direct sums and products, and (3) is straightforward.
\end{proof}

The above lemma enables us to define $F^{-1}U$ for a weakly open subspace $U$ of $Y$, as follows.  Write $\Mod U\simeq\Mod Y/T$ for some localizing subcategory $T$, and define $F^{-1}U$ to be the weakly open subspace of $X$ with $\Mod F^{-1}U\simeq \Mod X/F^{-1}T$. Since $F^{-1}(\bigcap_{i\in I}T_i)=\bigcap_{i\in I}F^{-1}T_i$, we see that $F^{-1}(\bigcup_{i\in I}U_i)=\bigcup_{i\in I}F^{-1}U_i$.  In particular if $\{U_i:i\in I\}$ is a weakly open cover of $Y$, then $\{F^{-1}U_i:i\in I\}$ is a weakly open cover of $X$.  Finally, note that $F^{-1}$ preserves containments: if $U\subseteq V$ are weakly open subspaces of $\Mod Y$, then $F^{-1}U\subseteq F^{-1}V$.

The definition of $F^{-1}U$ enables us to study a Frobenius bimodule locally.

\begin{lemma} Let $\M=(F,G)$ be a Frobenius bimodule between spaces $X$ and $Y$, and let $V$ and $U$ be weakly open subspaces of $X$ and $Y$ respectively, with inclusions $j_V$ and $j_U$.   If $F^{-1}U\subseteq V$, then there exists an exact functor $\bar F:\Mod V\rightarrow \Mod U$ having a right adjoint such that $\bar Fj^*_V=j_U^*F$.  Moreover, if $V=F^{-1}U$ is nonzero, then $\bar F$ is faithful.\label{restriction lemma}
\end{lemma}

\begin{proof}  Write $\Mod U\simeq\Mod Y/T$ and $\Mod V\simeq \Mod X/S$ for localizing subcategories $T$ and $S$.  Since $F^{-1}U\subseteq V$ we have $S\subseteq F^{-1}T$.  Since $j^*_UF$ is an exact functor which vanishes on $F^{-1}T$, and consequently on $S$, the universal property of quotient categories \cite[Corollary II.1.2]{Gabriel} ensures the existence of an exact functor $\bar F$ satisfying $\bar F j_V^*=j^*_UF$. Since each of $j^*_U$, $j^*_V$, and $F$ commute with sums, so too must $\bar F$.  Then $\bar F$ has a right adjoint by the Adjoint Functor Theorem. For the final statement, suppose that $V=F^{-1}U$ is nonzero. If $M\in\Mod F^{-1}U$, then there exists $N\in\Mod X$ with $N=j^*_{F^{-1}U}M$, and we have $\bar Fj^*_{F^{-1}U}N=j^*_UF(N)=0$, so that $F(N)\in T$.  Since $N\in F^{-1}T$ we must have $M=j^*_{F^{-1}U}N=0$. Since $\bar F$ takes nonzero objects to nonzero objects it is faithful.  
\end{proof}

The most important case of the above lemma is when $V=F^{-1}U$; in this case we write $F|_U$ for $\bar F$ and denote its right adjoint by $G|_U$.  The corresponding $F^{-1}U,U$-bimodule will be denoted by $\M|_U$ and called the \emph{restriction of $\M$ to $U$}.  Note that the bimodule $\M|_U$ is \emph{faithfully flat}, in the sense that $F|_U$ is an exact and faithful functor.  In general $\M|_U$ need not be a Frobenius bimodule, as the following example illustrates.  

\begin{example}  If $X$ and $Z$ are spaces, then their \emph{disjoint union} $X\sqcup Z$ is defined by setting $\Mod X\sqcup Z=\Mod X\times \Mod Z$; thus objects in $\Mod X\sqcup Z$ are pairs $(M,N)$ with $M\in\Mod X$ and $N\in\Mod Z$, with homomorphisms taken componentwise.

Let $X$ be a space, and let $Y=X\sqcup X$ be the disjoint union of two copies of $X$.    We write $\Delta:\Mod X\rightarrow \Mod Y$ for the functor defined on modules by $\Delta(M)=(M,M)$. The direct sum and direct product functors $\oplus:\Mod Y\rightarrow\Mod X$ and $\Pi:\Mod Y\rightarrow \Mod X$, defined on modules by $\oplus((M,N))=M\oplus N$ and $\Pi((M,N))=M\prod N$, are left and right adjoint to $\Delta$, respectively.  Since there is a natural equivalence of functors $\oplus\simeq\Pi$, we see that $_X\M_Y=(\Delta,\Pi)$ is a Frobenius bimodule.

Now, let $V$ be a weakly open subspace of $X$ such that the inclusion functor $j_*:\Mod V\rightarrow \Mod X$ is not exact ($j_*$ is only left exact in general), and write $\Mod V\simeq\Mod X/T$ for a localizing subcategory $T$.  Then $T\times 0=\{(M,0):M\in T\}$ is a localizing subcategory of $\Mod Y$, and its corresponding weakly open subspace is $U=V\sqcup X$.  Now, $\Delta^{-1}(T\times 0)=0$, and so the restriction of $\M$ to $U$ is an $X,U$-bimodule.  An easy computation shows that $\Delta|_U$ is given by $\Delta|_U(M)=(j^*M,M)$ and $\Pi|_U$ is given by $\Pi|_U(M,N)=j_*M\prod N$.  Since $j_*$ is not exact, neither is $\Pi|_U$ and hence $\M|_U$ is not Frobenius.\qed
\label{Delta example}
\end{example}

The following result gives a sufficient condition for $\M|_U$ to be Frobenius.  

\begin{prop} Let $\M=(F,G)$ be a Frobenius $X,Y$-bimodule, and let $U$ be a weakly open subspace of $Y$, with $\Mod U\simeq \Mod Y/T$. If $T\subseteq G^{-1}F^{-1}T$, then $\M|_U$ is a Frobenius $F^{-1}U,U$-bimodule.
\label{local Frobenius prop}
\end{prop}

\begin{proof} To ease notation, let $V=F^{-1}U$, and denote the inclusions by $j:U\rightarrow Y$ and $k:V\rightarrow X$.  Since $T\subseteq G^{-1}F^{-1}T$, we can apply Lemma \ref{restriction lemma} to $\M^*$ and $V$ to obtain a $U,V$-bimodule $(\bar G,\bar F)$ such that $\bar Gj^*=k^*G$.  We will show that $F|_U\simeq \bar F$ and $G|_U\simeq \bar G$.  

Note that each of $(k^*G,Fk_*)$ and $(\bar Gj^*,j_*\bar F)$ is an adjoint pair; since $k^*G=\bar Gj^*$ and adjoints are unique up to natural equivalence we have that $Fk_*\simeq j_*\bar F$.  Composing on the left by $j^*$ and using $j^*j_*\simeq \Id_U$, we see that $j^*Fk_*\simeq \bar F$.  But $j^*F=F|_Uk^*$, and $k^*k_*\simeq\Id_V$, so that $F|_U\simeq\bar F$.  Similarly each of $(j^*F,Gj_*)$ and $(F|_Uk^*,k_*G|_U)$ is an adjoint pair, so that $k_*G|_U\simeq Gj_*$. Composing on the left with $k^*$ gives $G|_U\simeq k^*Gj_*$, and using $k^*G=\bar Gj^*$ gives $G|_U\simeq \bar Gj^*j_*\simeq \bar G$.
\end{proof}

Using the Hom and Tensor notation for bimodules, the formulas relating $\M$ to $\M|_U$ can be restated as saying that whenever $T\subseteq G^{-1}F^{-1}T$, there are functorial isomorphisms for all $M\in\Mod X$ and $N\in\Mod Y$:
\begin{equation}
\begin{split}
M|_{F^{-1}U}\otimes_{F^{-1}U}\M|_U&\cong (M\otimes_X\M)|_U\\
\HOM_U(\M|_U,N|_U)&\cong \HOM_Y(\M,N)|_{F^{-1}U},
\end{split}\label{eq:projection}
\end{equation}
We refer to these as the \emph{projection formulas} for $\M$.

We can prove analogous results to the above starting with the functor $G$ instead of $F$; specifically, if $V$ is a weakly open subspace of $X$, then $_X\M_Y$ induces a $V,G^{-1}V$-bimodule which we can denote by $_V|\M$.  Furthermore, if $\Mod V\simeq \Mod X/S$ and $S\subseteq F^{-1}G^{-1}S$, then $_V|\M$ is a Frobenius $V,G^{-1}V$-bimodule.  A moments consideration shows that, in this case, we have the formula $(_V|\M)^*=\M^*|_V$. 

One technical difficulty that we shall encounter below is that, given a Frobenius bimodule $\M=(F,G)$ between spaces $X$ and $Y$, neither $F$ nor $G$ need be faithful.  The above result enables us to bypass this difficulty by passing to suitable open subspaces of $X$ and $Y$, as we now show.  

\begin{defn} Let $\M=(F,G)$ be a Frobenius $X,Y$-bimodule.  We define the \emph{kernels} of $F$ and $G$ to be $\ker F=\{M\in\Mod X:F(M)=0\}$ and $\ker G=\{N\in\Mod Y:G(N)=0\}$.  We define the \emph{supports} of $F$ and $G$ to be $\Supp(F)=\{x\in\Inj(X):\mbox{$F(M)\neq 0$ for all nonzero $M\leq E(x)$}\}$ and $\Supp(G)=\{y\in\Inj(Y):\mbox{$G(N)\neq 0$ for all nonzero $N\leq E(y)$}\}$.  
\end{defn}

Note that, since $F$ and $G$ are exact and commute with sums and products, $\ker F$ and $\ker G$ are closed subspaces of $X$ and $Y$, respectively.   We set $X_F=X\setminus \ker F$ and $Y_G=Y\setminus \ker G$, so that $X_F$ and $Y_G$ are open subspaces of $X$ and $Y$, respectively.  Finally, note that $\ker F$ and $\ker G$ are in fact localizing subcategories of $\Mod X$ and $\Mod Y$, respectively, so that $\Mod X_F\simeq \Mod X/\ker F$ and $\Mod Y_G=\Mod X/\ker G$.

\begin{prop} Let $\M=(F,G)$ be a nonzero Frobenius bimodule between spaces $X$ and $Y$.
\begin{enumerate}
\item $\ker GF=\ker F$ and $\ker FG=\ker G$.
\item $X_F=F^{-1}Y_G$ and $Y_G=G^{-1}X_F$.
\item If $X$ and $Y$ are noetherian, then $\Inj(X_F)=\Supp(F)$ and $\Inj(Y_G)=\Supp(G)$.
\item Writing $\M|_{Y_G}=(\bar F,\bar G)$, each of $\bar F$ and $\bar G$ is faithful.
\end{enumerate}\label{faithful prop}
\end{prop}

\begin{proof} (1) Clearly $\ker F\subseteq\ker GF$.  If $M\in\Mod X$ with $GF(M)=0$, then $\Hom_X(M,GF(M))=0$.  By the adjoint isomorphism $\Hom_Y(F(M),F(M))=0$.  Thus $F(M)=0$ and $M\in\ker F$.  The proof that $\ker G=\ker FG$ is analogous.

(2) By definition, $F^{-1}Y_G$ is the weakly open subspace of $X$ with $\Mod F^{-1}Y_G\simeq \Mod X/F^{-1}\ker G$.  But $F^{-1}\ker G=\{M:F(M)\in\ker G\}=\{M:GF(M)=0\}=\ker GF$.  By part (1), $\ker GF=\ker F$, so that $\Mod F^{-1}Y_G\simeq \Mod X/\ker F$.  Thus $F^{-1}Y_G=X_F$.  Similarly we see that $G^{-1}X_F=Y_G$.
 
(3) We identify $\Inj(X_F)$ with $\{x\in\Inj(X):\mbox{$E(x)$ is $\ker F$-torsionfree}\}$.  But $E(x)$ is $\ker F$-torsionfree if and only if $F(M)\neq 0$ for all nonzero submodules $M$ of $E(x)$, if and only if $x\in\Supp(F)$.  Similarly $\Inj(Y_G)=\Supp(G)$.

(4)  This follows from Lemma \ref{restriction lemma} applied to each of $\M$ and $\M^*$.
\end{proof}

\begin{remark} (1) The fact that each of $X_F$ and $Y_G$ is open, and not just weakly open, is a strong condition.  Suppose that $X$ is affine, say $X\simeq \Mod R$.  Then there is a bijection between the closed subspaces of $X$ and the two-sided ideals of $R$ \cite{Rosenberg book}.  In particular, $\ker F\simeq \Mod R/I$ for some two-sided ideal $I$ of $R$.  Since $\ker F$ is also closed under extensions, we have that $R/I^2\in\Mod R/I$, since there is an exact sequence $0\rightarrow I/I^2\rightarrow R/I^2\rightarrow R/I\rightarrow 0$, and $I/I^2$ and $R/I$ are both in $\Mod R/I$.  Thus $I$ kills $R/I^2$, from which it follows that $I^2=I$.  Hence we conclude that in the affine case, $\ker F\simeq \Mod R/I$, where $I$ is an idempotent two-sided ideal of $R$.

(2) More generally, if $X$ is a scheme admitting an ample line bundle, then \cite[Theorem 4.1]{Smith subspaces} shows that every closed subspace of $\QCoh(\O_X)$ is of the form $\QCoh(\O_Z)$ for some closed subscheme $Z$ of $X$.  Thus in this case it is possible to prove as above that $\ker F=\QCoh(\O_X/\sh{I})$, where $\sh{I}$ is a sheaf of ideals of $\O_X$ satisfying $\sh{I}^2=\sh{I}$.

(3) If $\ker F$ is closed under injective envelopes, then $F(M)=0$ for some nonzero submodule of $E(x)$ if and only if $F(E(x))=0$; consequently in this case we have that $\Supp(F)=\{x\in\Inj(X):F(E(x))\neq 0\}$.  This applies in particular if $X$ is a noetherian scheme, because every localizing subcategory of $\QCoh(\O_X)$ is closed under injective envelopes in this case \cite[Proposition VI.2.4]{Gabriel}.  
\end{remark}

\begin{lemma} Let $_X\M_Y$ be a Frobenius bimodule, and let $S$ and $T$ be simple $X$- and $Y$-modules, respectively. 
\begin{enumerate}
\item $S\otimes_X\M$ has an essential socle.
\item If $S\otimes_X\M\neq 0$, then $E(S)\otimes_X\M$ is an injective hull for $S\otimes_X\M$; thus $E(S)\otimes_X\M\cong \oplus_{i=1}^t E(T_i)$ for simple $Y$-modules $T_i$.
\item If $S\otimes_X\M\neq 0$ and $E(T)$ is a summand of $E(S)\otimes_X\M$, then $E(S)$ is a summand of $\HOM_Y(\M,E(T))$. 
\item If $E(T)$ is a summand of $E(S)\otimes_X\M$ then $\HOM_Y(\M,T)\neq 0$.
\end{enumerate}
\label{injective hull lemma}
\end{lemma}

\begin{proof} (1) Since $S\otimes_X\M$ is finitely copresented, hence finitely cogenerated, this follows from \cite[Theorem 10.4(2)]{AF}. (Their proof works in any Grothendieck category.)

(2)  Let $f:S\otimes_X\M\rightarrow E$ be monic with $E$ injective.  Then there exists $\tilde f:S\rightarrow \HOM_Y(\M,E)$, which is necessarily monic because $S$ is simple.  Thus there is a unique monic $\tilde g:E(S)\rightarrow \HOM_Y(\M,E)$ such that $\tilde g\circ i=\tilde f$, where $i:S\rightarrow E(S)$ is the canonical map.  It follows  that there is a unique $g:E(S)\otimes_X\M\rightarrow E$ with $g(i\otimes_X\M)=f$.  The last statement follows from (1) because the socle of $S\otimes_X\M$ is a finite direct sum of simple $Y$-modules \cite[Proposition 9.7]{AF}.

(3) Since $S\otimes_X\M\neq 0$, $E(S)\otimes_X\M$ is an injective hull of $S\otimes_X\M$ by part (2).  Thus there is a nonzero $f\in\Hom_Y(S\otimes_X\M,E(T))$, so that $\Hom_X(S,\HOM_Y(\M,E(T)))\neq 0$. Since $S$ is simple any nonzero morphism must be monic, so that $E(S)$ is a summand of $\HOM_Y(\M,E(T))$.

(4) Since $E(T)$ is a summand of $E(S)\otimes_X\M$, we have $\Hom_Y(T,E(S)\otimes_X\M)\neq 0$, so that $\Hom_X(\HOM_Y(\M,T),E(S))\neq 0$.  Thus $\HOM_Y(\M,T)$ is nonzero. 
\end{proof}

We close this section with some results on Frobenius bimodules between (semi)local spaces.

\begin{prop} Let $\M$ be a nonzero Frobenius bimodule between spaces $X$ and $Y$, and assume that $X$ is local.  
\begin{enumerate}
\item $-\otimes_X\M$ is faithful.
\item If $\HOM_Y(\M,-)$ is faithful, then $Y$ is semilocal.
\end{enumerate}
\end{prop}

\begin{proof} (1) Let $S$ be the unique (up to isomorphism) simple $X$-module.  Write $\M=(F,G)$ to ease notation, and note that by definition $\ker G=G^{-1}0$.  Since $0\subseteq F^{-1}G^{-1}0$, we can apply Proposition \ref{local Frobenius prop} to $\M^*$ to conclude that there is a Frobenius $X,Y_G$-bimodule $\bar \M=(\bar F,\bar G)$, and it follows from Lemma \ref{restriction lemma} that $\bar G$ is faithful.  Moreover, since the localizing subcategory of $\Mod X$ used is just the zero subcategory, we see that $\bar F=j^*F$, where $j:Y_G\rightarrow Y$ denotes the inclusion.  Now, let $T$ be a simple $Y_G$-module.  Then $\bar G(T)\neq 0$, as $\bar G$ is faithful. It follows from part (2) of Lemma \ref{injective hull lemma} that $\bar G(E(T))\cong E(S)^{(n)}$ for some positive integer $n$, and then part (4) of Lemma \ref{injective hull lemma} implies that $\bar F(S)=j^*F(S)$ is nonzero, so that $F(S)\neq 0$. The proof is concluded by noting that, since $S\in\sigma[M]$ for every $X$-module $M$ and $F(S)\neq 0$, we must have $F(M)\neq 0$ for all $X$-modules $M$.  Thus $F$ is faithful.

(2) Keeping the above notation, Proposition \ref{injective cogenerator prop}(2) shows that $F(E(S))$ is an injective cogenerator for $\Mod Y$.  But Proposition \ref{injective hull lemma}(2) shows that $F(E(S))$ is isomorphic to a finite direct sum of injective hulls of simple $Y$-modules.  By Lemma \ref{semilocal lemma}, $Y$ is semilocal.
\end{proof}

\section{Dimension preserving bimodules}  
In this section we consider Frobenius bimodules between spaces equipped with dimension functions.  We shall assume that all dimension functions are exact, finitely partitive, ordinal valued, and commute with direct limits. If $d$ is a dimension function on $\Mod X$ and $\a$ is an ordinal,  then an $X$-module $M$ is called \emph{$\alpha$-homogeneous} (with respect to $d$) if $d(N)=d(M)=\alpha$ for every nonzero submodule $N$ of $M$, and $M$ is called \emph{$\a$-critical} if $d(M)=\a$ and $d(M/N)<\a$ for all proper quotients $M/N$ of $M$.  

We note that if $X$ is noetherian, then Krull dimension in the sense of Gabriel \cite[p. 382]{Gabriel} is such a dimension function. Moreover, the Krull-Schmidt Theorem holds for injective modules when $X$ is noetherian:  If $E\in\Mod X$ is injective, then there exist indecomposable injectives $\{E(x_i):i\in I\}$ and cardinals $\kappa_i$ such that $E\cong\oplus_{i\in I}E(x_i)^{(\kappa_i)}$, and this decomposition is unique up to isomorphism.

\begin{defn} Let $X$ and $Y$ be spaces with dimension functions $d$ and $\d$, respectively. We say that a functor $F:\Mod X\rightarrow\Mod Y$ is \emph{dimension preserving} if $\delta(F(M))=d(M)$ whenever $F(M)\neq 0$.  A Frobenius bimodule $_X\M_Y$ is called \emph{dimension preserving} if each of $-\otimes_X\M$ and $\HOM_Y(\M,-)$ preserves dimension.
\end{defn}

We remark in passing that if $_X\M_Y$ is dimension preserving, then so too is its dual $_Y\M^*_X$. In general a Frobenius bimodule need not be dimension preserving for given dimension functions $d$ and $\d$. Indeed, if one uses two different dimension functions $d$ and $\d$ on a single space $X$, then even the identity functor need not be dimension preserving.  In subsequent sections we shall primarily apply the following results in the case where $X$ and $Y$ are noetherian spaces, and $d$ and $\d$ are both taken to be Krull dimension.  However, there is no advantage to specializing to this case immediately, and so we develop the ideas in a more formal framework.    

Unless stated to the contrary, we assume throughout this section that $X$ and $Y$ are spaces equipped with fixed dimension functions $d$ and $\d$, respectively.  When we say that a Frobenius bimodule is dimension preserving, it is understood that this means with respect to the given dimension functions $d$ and $\d$.

\begin{lemma}  Let $_X\M_Y$ be a dimension preserving Frobenius bimodule and let $M$ be an $\a$-homogeneous $Y$-module.  If $\HOM_Y(\M,M)$ is nonzero then it is an $\a$-homogeneous $X$-module.\label{homogeneous lemma}
\end{lemma}

\begin{proof}  Let $N$ be an $X$-module with $d(N)<\alpha$. Then since $\M$ is dimension preserving $\d(N\otimes_X\M)<\alpha$ also.  Since $M$ is $\a$-homogeneous, we have $\Hom_Y(N\otimes_X\M,M)=0$, so that $\Hom_X(N,\HOM_Y(\M,M))=0$ by the adjoint isomorphism.  In particular if $K$ is a nonzero submodule of $\HOM_Y(\M,M)$, then $d(K)\geq\a$. Since  $d(\HOM_Y(\M,M))=\a$ we conclude that $d(K)=\a$ and so $\HOM_Y(\M,M)$ is $\a$-homogeneous.
\end{proof}

We need to fix some more notation.  Given an ordinal $\a$, we let $S_\a=\{M\in\Mod X:d(M)<\a\}$, and similarly $T_\a=\{N\in\Mod Y:\d(N)<\a\}$.  Each $S_\a$ is a localizing subcategory of $\Mod X$, and similarly each $T_\a$ is a localizing subcategory of $\Mod Y$.  We write $V_\a$ for the weakly open subspace of $\Mod X$ with $\Mod V_\a\simeq \Mod X/S_\a$, and write $k_\a:V_\a\rightarrow X$ for the inclusion.  Similarly we let $U_\a$ be the weakly open subspace of $\Mod Y$ with $\Mod U_\a\simeq \Mod Y/T_\a$, and we denote the inclusion by $j_\a:U_\a\rightarrow Y$.  Finally, we denote by $\Inj_\a(X)$ the set of those $x\in\Inj(X)$ such that the critical dimension of $E(x)$ is $\a$. (The \emph{critical dimension} of $E(x)$ is the dimension of a critical submodule of $E(x)$.)

If $_X\M_Y=(F,G)$ is a dimension preserving Frobenius bimodule, then $S_\a\subseteq F^{-1}T_\a$ and $T_\a\subseteq G^{-1}S_\a$; combining these gives $T_\a\subseteq G^{-1}F^{-1}T_\a$.  In particular Proposition \ref{local Frobenius prop} shows that there is a Frobenius $V_\a,U_\a$-bimodule $\M_\a=(F_\a,G_\a)$ satisfying the following formulas:
\begin{equation}\begin{split}
k^*_\a(-)\otimes_{V_\a}\M_\a&=j^*_\a(-\otimes_X\M)\\
\HOM_{U_\a}(\M_\a,j^*_\a(-))&=k^* _\a(\HOM_Y(\M,-)).\label{alpha formulas}
\end{split}
\end{equation}
We can induce dimension functions on $V_\a$ and $U_\a$ be setting $d(k^*_\a M)=d(M)-\a$ and $\delta(j^*_\a N)=\d(N)-\a$ for $M\in\Mod X$ and $N\in\Mod Y$, making $\M_\a$ a dimension preserving $V_\a,U_\a$-bimodule. (Here we abuse notation and write $d$ and $\d$ for the induced dimension functions.  This should not cause confusion.)  We can (and shall) identify $\Inj(V_\a)$ with $\bigcup_{\b\geq\a}\Inj_\b(X)$, and similarly we identify $\Inj(U_\a)$ with $\bigcup_{\b\geq\a}\Inj_\b(Y)$.  If $x\in\Inj_\b(X)$ for $\b\geq\a$, then $x\in\Inj_{\b-\a}(V_\a)$, and similarly if $y\in\Inj_\b(Y)$ for $\b\geq\a$, then $y\in\Inj_{\b-\a}(U_\a)$.

\begin{lemma}  Let $_X\M_Y=(F,G)$ be a dimension preserving Frobenius bimodule, and let $x\in\Inj_\a(X)$.  Then $x\in\Supp(F)$ if and only if $F(M)\neq 0$ for some (hence every) critical submodule of $E(x)$.\label{nonzero crit lemma}
\end{lemma}

\begin{proof}  If $x\in\Supp(F)$, then clearly $F(M)\neq 0$.  Suppose now that $F(M)\neq 0$, and let $N\leq E(x)$ be a submodule with $F(N)=0$.  If $N\neq 0$, then $N\cap M$ is a nonzero submodule of $M$ with $F(N\cap M)=0$.  This implies that $F(M)=F(M/N\cap M)$.  But $\d(F(M))=\a$ since $F(M)$ is nonzero and $\M$ is dimension preserving, while $\d(F(M/N\cap M))<\a$ since $d(M/N\cap M)<\a$.  Thus $N\cap M=0$ and $F(N)$ is nonzero for every $N\leq E(x)$.
\end{proof}

\begin{thm}  Let $_X\M_Y=(F,G)$ be a nonzero dimension preserving Frobenius bimodule between noetherian spaces $X$ and $Y$, and let $x\in\Supp(F)\cap \Inj_\a(X)$.  Then there exists a positive integer $n_x$ depending on $x$ such that $E(x)\otimes_X\M\cong\oplus_{i=1}^{n_x} E(y_i)^{(n_i)}$ with each $y_i\in\Supp(G)\cap \Inj_\a(Y)$. \label{decomp thm} \end{thm}

\begin{proof}  Suppose first that $\a=0$, so that $E(x)$ is the injective hull of a simple $X$-module $S$.  Then $E(x)\otimes_X\M$ is the injective hull of $S\otimes_X\M$ by Lemma \ref{injective hull lemma}(2), and $S\otimes_X\M$ is noetherian and artinian. The first follows because $S\otimes_X\M$ is finitely generated and $X$ is a noetherian space. For the latter, note that $S\otimes_X\M$ is $0$-homogeneous by Lemma \ref{homogeneous lemma}, and every $0$-homogeneous module is artinian \cite[Lemma 3.10(a)]{Pappacena injective}. Thus $S\otimes_X\M$ is a finite-length $X$-module; it follows that if $E(x)\otimes_X\M\neq 0$, then it is isomorphic to a finite direct sum of injective hulls of $0$-critical $Y$-modules.  Hence $E(x)\otimes_X\M\cong \oplus_{i=1}^{n_x}E(y_i)^{(n_i)}$, with each $E(y_i)\in\Inj_0(Y)$. Moreover each $y_i\in\Supp(G)$ by Lemma \ref{injective hull lemma}(4). This proves the result for $\a=0$.

Now suppose $\a>0$, and let $E(x)\in\Supp(F)\cap\Inj_\a(X)$.  Since $E(x)\otimes_X\M$ is injective, we can write 
\begin{equation}
E(x)\otimes_X\M\cong \bigl(\bigoplus_{y\in J} E(y)^{(\kappa_y)}\bigr)\oplus E,\label{decomp eq}\end{equation}
where each $\kappa_y$ is a cardinal, $J\subseteq\Inj_\a(Y)$, and each indecomposable summand of $E$ has critical dimension different from $\a$.  

If $M$ is a critical submodule of $E(x)$, then $M\otimes_X\M\neq 0$, so that $M\otimes_X\M$ is $\a$-homogeneous.  In particular $j^*_\a(M\otimes_X\M)\neq 0$, so that $k_\a^*M\otimes_{V_\a}\M_\a$ is nonzero.  Since $k_\a^*M$ is a simple $V_\a$-module and $k^*_\a E(x)$ is its injective hull, we see that $x\in\Supp(F_\a)\cap\Inj_0(V_\a)$.

It follows by the first paragraph, applied to the $V_\a,U_\a$-bimodule $\M_\a$, that
\begin{equation}k^*_\a E(x)\otimes_{V_\a}\M_\a\cong\bigoplus_{i=1}^{n_x} j^*_\a E(y_i)^{(n_i)},\end{equation}
where each $E(y_i)\in\Inj_0(U_\a)=\Inj_\a(Y)$. Appealing again to Lemma \ref{injective hull lemma}(4), we have that $y_i\in\Supp(G_\a)$.  If $N$ is an $\a$-critical submodule of $E(y_i)$, then $\HOM_{U_\a}(\M_\a,j^*_\a N)\neq 0$ implies that $\HOM_Y(\M,N)\neq 0$, so that $y_i\in\Supp(G)$ by Lemma \ref{nonzero crit lemma}.

Now, the formula $k^*_\a E(x)\otimes_{V_\a}\M_\a\cong j_\a^*(E(x)\otimes_X\M)$ and equation \eqref{decomp eq} give that  
\begin{equation}\bigoplus_{i=1}^{n_x} j^*_\a E(y_i)^{(n_i)}\cong \bigl(\bigoplus_{y\in J} j^*_\a E(y)^{(\kappa_y)}\bigr)\oplus j^*_\a E.\end{equation}
Using the uniqueness of decomposition of indecomposable injectives in $\Mod U_\a$, we see that 
\[E(x)\otimes_X\M\cong \bigl(\bigoplus_{i=1}^{n_x} E(y_i)^{(n_i)}\bigr)\oplus E,\]
where each $E(y_i)$ has critical dimension $\a$ and $\d(E)<\a$.  

We finish the proof by showing that $E=0$. The adjoint isomorphism gives $\Hom_Y(E,E(x)\otimes_X\M)\cong\Hom_X(\HOM_Y(\M,E),E(x))$. Since $d(\HOM_Y(\M,E))=\delta(E)<\alpha$, and every nonzero submodule of $E(x)$ has dimension $\geq \a$, it follows that $\Hom_X(E,E(x)\otimes_X\M)=0$.  Since $E$ is a summand of $E(x)\otimes_X\M$, we must have $E=0$.
\end{proof}  

\begin{cor} Let $\M$ be a dimension preserving Frobenius bimodule between noetherian spaces $X$ and $Y$. If $M$ is an $\a$-critical $X$-module with injective hull $E(x)$ and $M\otimes_X\M\neq 0$, then $E(x)\otimes_X\M$ is an injective hull for $M\otimes_X\M$.
\label{essential containment cor}
\end{cor}

\begin{proof} Arguing as in the proof of the previous theorem, we see that $k^*_\a M\otimes_{V_\a}\M_\a\neq 0$.  It follows from Lemma \ref{injective hull lemma}(2) that $k^*_\a E(x)\otimes_{V_\a}\M_\a$ is an injective hull for $k^*_\a M\otimes_{V_\a}\M_\a$.  Using formula \eqref{alpha formulas} we conclude that $j^*_\a(E(x)\otimes_X\M)$ is an injective hull for $j^*_\a(M\otimes_X\M)$.  Since $j_{\a*}$ preserves essential containments and $E(x)\otimes_X\M$ is $T_\a$-torsionfree, we have that $E(x)\otimes_X\M$ is an injective hull of $j_{\a*}j_\a^*(M\otimes_X\M)$.  Finally, we note that $M\otimes_X\M\leq j_{\a*}j_\a^*(M\otimes_X\M)$ because $M\otimes_X\M$ is $T_\a$-torsionfree.
\end{proof}  

\begin{cor}  Let $\M$ be a dimension preserving Frobenius bimodule between noetherian spaces $X$ and $Y$, and let $x\in\Supp(F)\cap\Inj_\a(X)$. If $E(y)$ is a isomorphic to a summand of $E(x)\otimes_X\M$, then $y\in\Supp(G)$ and $E(x)$ is isomorphic to a summand of $\HOM_Y(\M,E(y))$.\label{dual summand cor}
\end{cor}

\begin{proof} By Theorem \ref{decomp thm}, $y\in\Inj_\a(Y)=\Inj_0(U_\a)$.  Let $N$ be a critical submodule of $E(y)$, so that $j^*_\a N$ is a simple $U_\a$-module. Now, $k^*_\a M$ is a simple $V_\a$-module and $x\in\Supp(F_\a)\cap\Inj_0(V_\a)$.  By Lemma \ref{injective hull lemma}(3) we conclude that $k^*_\a E(x)$ is a summand of $\HOM_{U_\a}(\M_\a,j^*_\a E(y))=k^*_\a\HOM_Y(\M,E(y))$.  Applying $k_{\a*}$ and using the fact that $E(x)$ and $\HOM_Y(\M,E(y))$ are $S_\a$-torsionfree shows that $E(x)$ is isomorphic to a summand of $\HOM_Y(\M,E(y))$.  Also, $\HOM_{U_\a}(\M_\a,j^*_\a N)=k^*_\a\HOM_Y(\M,N)$ is nonzero by Lemma \ref{injective hull lemma}(4).  By Lemma \ref{nonzero crit lemma}, $y\in\Supp(G)$. 
\end{proof}

\begin{cor} Let $\M=(F,G)$ be a dimension preserving bimodule between noetherian spaces $X$ and $Y$, and assume that $\ker F$ is closed under injective envelopes.  Then $F$ preserves essential containments.
\label{strong essential cor}
\end{cor}

\begin{proof}  Since $\ker F$ is closed under injective envelopes, we see that if $F(M)=0$ for some submodule of $E(x)$, then $F(E(x))=0$.  This, combined with Corollary \ref{essential containment cor}, shows that $F(E(M))$ is an injective hull of $F(M)$ for all critical $X$-modules $M$.  Now, if $M$ is an arbitrary noetherian $X$-module, then $M$ contains an essential submodule which is a finite direct sum of critical $X$-modules, say $M_1\oplus\dots\oplus M_t\leq M$, and $E(M)=E(M_1)\oplus \dots \oplus E(M_t)$.
The result now follows from the fact that $F(E(M))$ is an injective hull for $F(M_1)\oplus\dots \oplus F(M_t)$.
\end{proof}

We close this section by showing that for noetherian spaces $X$ and $Y$ which are close to being commutative, every Frobenius bimodule $_X\M_Y$ preserves Krull dimension.  The precise condition that we impose is the following:

\begin{enumerate} 
\item[$(*)$]If $\Hom_X(E(x_1),E(x_2))\neq 0$ with $E(x_1)\in\Inj_\a(X)$, $E(x_2)\in\Inj_\beta(X)$, then $\a\geq \beta$.
\end{enumerate}
Condition $(*)$ is satisfied with respect to Krull dimension if $\Mod X\simeq \Mod R$, where $R$ is a commutative ring, or more generally a (two-sided) FBN ring. Indeed, let $M$ and $N$ be critical $R$-modules with $\Hom_R(E(M),E(N))\neq 0$.  If $f(M)\neq 0$ for some $f$, then we have that $\Kdim f(M)<\Kdim M$, so that $\Kdim N<\Kdim M$.  If, on the other hand, we have that $f(M)=0$ for all $f\in\Hom_X(E(M),E(N))$ then, letting $U$ be a prime submodule of $E(M)$ with $f(U)=0$, we can find $K\leq E(M)$ and $V\leq E(N)$ prime such that the exact sequence $0\rightarrow U\rightarrow K\rightarrow V\rightarrow 0$ satisfies the hypotheses of ``Jategaonkar's Main Lemma" \cite[Theorem 11.1]{Goodearl Warfield}. Since $R$ satisfies the strong second layer condition, we conclude that there is a link $\ass E(N)\leadsto\ass E(M)$. By \cite[Corollary 12.6 and Theorem 13.13]{Goodearl Warfield}, we conclude that $\Kdim M=\Kdim N$ in this case. We shall see in Lemma \ref{scheme star lemma} that condition $(*)$ also holds for dimension of support on $\QCoh(\O_X)$, where $X$ is a noetherian scheme.

\begin{thm} Let $X$ and $Y$ be noetherian spaces satisfying condition $(*)$ with respect to Krull dimension.  Then every Frobenius bimodule $_X\M_Y$ preserves Krull dimension.\label{Krull dim thm}
\end{thm}

\begin{proof}

The definition of dimension preserving shows that the zero bimodule preserves every dimension function, so we may assume that $\M$ is nonzero.  Suppose first that $M\in\Mod X$ is simple, and let $N\in\Mod Y$ be an $\a$-critical quotient module of $M\otimes_X\M$. Then $\Hom_Y(E(M)\otimes_X\M,E(N))\neq 0$.  Writing $E(M)\otimes_X\M=\oplus_{i=1}^tE(S_i)$ for simple $Y$-modules $S_i$, we must have $\Hom_Y(E(S_i),E(N))\neq 0$ for some $S_i$.  Condition ($*$) in $\Mod Y$ then says that $N$ must be $0$-critical; i.e. $N$ is simple.  It follows that $M\otimes_X\M$ is $0$-dimensional.  A similar argument interchanging $X$ and $Y$ and using $\M^*$ shows that if $N\in\Mod Y$ is simple, then $\HOM_Y(\M,N)$ is $0$-dimensional.  Since every $0$-dimensional $X$-module is a direct limit of simple $X$-modules (and similarly for $Y$-modules), we see that each of $-\otimes_X\M$ and $\HOM_Y(\M,-)$ take $0$-dimensional modules to $0$-dimensional modules.

We now proceed by transfinite induction, assuming that $-\otimes_X\M$ and $\HOM_Y(\M,-)$ each preserve Krull dimension for all $X$- and $Y$-modules of Krull dimension strictly less than $\a$.  In particular we have that $S_\a\subseteq F^{-1}T_\a$ and $T_\a\subseteq G^{-1}F^{-1}T_\a$, where $\M=(F,G)$.  We show first that condition ($*$) passes down to $U_\a$ and $V_\a$.

Any nonzero indecomposable injective in $U_\a$ is of the form $j^*_\a E(y)$ for some $E(y)\in\Inj(Y)$ of critical dimension $\geq \a$.  So, if $\Hom_{U_\a}(j^*_\a E(y_1),j^*_\a E(y_2))\neq 0$, then since $E(y_i)\cong j_{\a*}j^*_\a E(y_i)$ for $i=1,2$, we have $\Hom_Y(E(y_1),E(y_2))\neq 0$. Condition $(*)$ in $\Mod Y$ then implies that the critical dimension of $E(y_2)$ is at most the critical dimension of $E(y_1)$; it follows that the critical dimension of $j^*_\a E(y_2)$ is at most the critical dimension of $j^*_\a E(y_1)$.  The argument for $V_\a$ is similar.

By Proposition \ref{local Frobenius prop}, there is a Frobenius $V_\a,U_\a$-bimodule $\M_\a$ such that formulas \eqref{alpha formulas} hold.  Now, if $M$ is an $\a$-critical $X$-module, then $k^*_\a M$ is a simple $V_\a$-module, and so $k^*_\a M\otimes_{V_\a}\M_\a\cong j^*_\a(M\otimes_X\M)$ is either zero or $0$-dimensional.  This says precisely that $\Kdim M\otimes_X\M\leq\a$, and similarly we see that $\Kdim\HOM_Y(\M,N)\leq \a$ for any $\a$-critical $Y$-module $N$.  If $M\otimes_X\M\neq 0$, let $\Kdim M\otimes_X\M=\b$, and let $N$ be a critical submodule of $M\otimes_X\M$.  Then $\Hom_Y(N,M\otimes_X\M)\neq 0$, so that $\Hom_X(\HOM_Y(\M,N),M)\neq 0$.  Since $\Kdim\HOM_Y(\M,N)\leq \b$ and $M$ is $\a$-critical, we conclude that $\b\geq\a$.  

Thus $\b=\a$ and $-\otimes_X\M$ preserves the Krull dimension of any critical $X$-module.   Since any noetherian $X$-module has a critical composition series and $-\otimes_X\M$ is exact, we see that $-\otimes_X\M$ preserves the dimension of any noetherian $X$-module, and since any $X$-module is the direct limit of its noetherian submodules, we see that $-\otimes_X\M$ preserves Krull dimension.  In a similar way we see that $\HOM_Y(\M,-)$ preserves Krull dimension. 
\end{proof}

\section{Right localizing bimodules}
We wish to study a Frobenius bimodule $_X\M_Y$ locally; that is, by studying one or both of the restrictions $\M|_U$ and $_V|\M$ for weakly open subspaces $V$ and $U$ of $X$ and $Y$, respectively.  As mentioned in section 3 above, the difficulty is that in general, the bimodules $\M|_U$ and $_V|\M$ need not be Frobenius.  The following definitions impose the sufficient conditions of Proposition \ref{local Frobenius prop} to ensure that $\M|_U$ is a Frobenius  $F^{-1}U,U$-bimodule for all weakly open subspaces $U$ of $Y$ (respectively, that $_V|\M$ is a Frobenius $V,G^{-1}V$-bimodule for all weakly open subspaces $V$ of $X$).

\begin{defn} Let $\M$ be a Frobenius bimodule between spaces $X$ and $Y$.  Then $\M$ is called \emph{right localizing} if $T\subseteq G^{-1}F^{-1}T$ for every localizing subcategory $T$ of $\Mod Y$. Similarly $\M$ is \emph{left localizing} if $S\subseteq F^{-1}G^{-1}S$ for every localizing subcategory $S$ of $\Mod X$.
\end{defn}

For spaces $X$ and $Y$, a category equivalence between $\Mod X$ and $\Mod Y$ is an easy example of a Frobenius bimodule that is both left and right localizing.  Also, if $_X\M_Y=(F_1,G_1)$ and $_Y\sh{N}_Z=(F_2,G_2)$ are either left or right localizing, then an easy computation shows that $\M\otimes_Y\sh{N}=(F_2F_1,G_1G_2)$ is, as well. The following example shows that the two notions are in general distinct.

\begin{example}  Let $X$ be a space, and let $Y=X\sqcup X$ be the disjoint union of two copies of $X$. If $_X\M_Y=(\Delta,\Pi)$ is the Frobenius bimodule of Example \ref{Delta example}, then we claim that $\M$ is left localizing but not right localizing.  

Given a localizing subcategory $T$ of $\Mod X$, write $T=T(E)$ for some injective $X$-module $E$.  Then 
\begin{multline*}\Delta^{-1}\Pi^{-1}T(E)=\{M:\Hom_X(\Pi\Delta(M),E)=0\}\\
=\{M:\Hom_X(M,\Pi\Delta(E))=0\}=T(\Pi\Delta(E))=T(E\textstyle\prod E)=T(E),\end{multline*}
showing that $\M$ is left localizing.

On the other hand, fix an injective $X$-module $E$, and consider the injective $Y$-module $E'=(E,0)$. Since $\Delta\Pi(E')=(E,E)$, we see that $\Pi^{-1}\Delta^{-1}T(E')=T((E,E))$.  This shows that $\M$ is not right localizing: $(0,E)$ is in $T(E')$ but not in $T((E,E))$.\qed\label{not symmetric example}
\end{example}

Nevertheless it is clear that $\M$ is left localizing if and only if $\M^*$ is right localizing, and so we restrict our attention to the latter type of bimodule.  
 
\begin{lemma} If $_X\M_Y=(F,G)$ is a right localizing Frobenius bimodule and $U$ is a weakly open subspace of $Y$, then $\M|_U$ is a right localizing $F^{-1}U,U$-bimodule.
\end{lemma}

\begin{proof}  Write $\Mod U\simeq \Mod Y/T$ for a localizing subcategory $T$ of $\Mod Y$.  Now, the localizing subcategories of $\Mod Y/T$ are of the form $S/T$, where $S$ is a localizing subcategory of $\Mod Y$ containing $T$.  The result then follows from the projection formulas \eqref{eq:projection} and the fact that $\M$ is right localizing.
\end{proof}

\begin{lemma} Let $Y$ be a noetherian space.  Then a Frobenius $X,Y$-bimodule $\M=(F,G)$ is right localizing if and only if $T_y\subseteq G^{-1}F^{-1}T_y$ for every $y\in\Inj(Y)$.
\label{localizing test lemma}
\end{lemma}

\begin{proof} One direction is clear.  Given a localizing subcategory $T$ of $\Mod Y$, we can write $T=\bigcap_{y\in\Sigma}T_y$ for some $\Sigma\subseteq \Inj(Y)$ by Lemma \ref{localizing subcategory lemma}.  Then 
\[T=\bigcap_{y\in\Sigma}T_y\subseteq \bigcap_{y\in\Sigma}G^{-1}F^{-1}T_y=G^{-1}F^{-1}\bigl(\bigcap_{y\in \Sigma}T_y\bigr)=G^{-1}F^{-1}T.\]
\end{proof}

If $_X\M_Y=(F,G)$ is a right localizing Frobenius bimodule, then an easy computation shows that $\M^*\otimes_X\M=(FG,FG)$ is both a left and right localizing Frobenius bimodule on $Y$.  

\begin{lemma} If $_X\M_Y$ is a right localizing Frobenius bimodule, then $\M^*\otimes_X\M$ preserves every dimension function on $\Mod Y$.\label{dim preserving lemma}
\end{lemma}

\begin{proof} Let $\d$ be a dimension function on $\Mod Y$. Writing $\M=(F,G)$, we must show that $FG$ preserves $\d$.  Retaining the notation of the previous section, we have that $T_\a$ is a localizing subcategory of $\Mod Y$ for every ordinal $\a$.  Since $\M$ is right localizing $T_\a\subseteq G^{-1}F^{-1}T_\a$ for all $\a$; this says precisely that $\d(FG(M))=\d(M)$ for every $Y$-module $M$ with $FG(M)\neq 0$.  
\end{proof}

\begin{prop} Let $_X\M_Y=(F,G)$ be a nonzero right localizing Frobenius bimodule, and suppose that $Y$ is noetherian.  Then for each $y\in\Supp(G)$, there exists a positive integer $n_y$ such that $FG(E(y))\cong E(y)^{(n_y)}$.\label{FG prop}
\end{prop}

\begin{proof}  Let the critical dimension of $E(y)$ be $\a$. Note that $G(M)\neq 0$ implies $FG(M)\neq 0$ for $Y$-modules $M$, since $\ker FG=\ker G$.  Thus $y\in\Supp(FG)$.  Since $\M^*\otimes_X\M=(FG,FG)$ preserves Krull dimension by Lemma \ref{dim preserving lemma}, we have by Theorem \ref{decomp thm} that $FG(E(y))\cong \oplus_{i=1}^{n_y}E(y_i)$, with each $y_i\in\Inj_\a(Y)$. This implies that $G^{-1}F^{-1}T_y=\bigcap_{i=1}^{n_y}T_{y_i}$. Since $\M$ is right localizing, we have $T_y\subseteq T_{y_i}$ for all $i$.  Fix $i$, and let $M_i$ be a critical submodule of $E(y_i)$. Then $M_i\not\in T_{y_i}$, so that $M_i\not\in T_y$;  that is, $\Hom_Y(M_i,E(y))\neq 0$.  Since $M_i$ is critical of dimension $\a$, $M_i$ is isomorphic to a submodule of $E(y)$. Hence $E(y)= E(y_i)$ for all $i$, showing that $FG(E(y))\cong E(y)^{(n_y)}$.
\end{proof}

\begin{lemma} Let $X$ and $Y$ be noetherian spaces, and let $_X\M_Y=(F,G)$ be a right localizing Frobenius bimodule.  Assume that $\ker G$ is closed under injective envelopes, and that $F$ is faithful. Then $G$ preserves essential containments.
\end{lemma}

\begin{proof}  Let $N\leq M$ be $Y$-modules with $N$ essential in $M$.  If $G(N)$ is not an essential submodule of $G(M)$, then there exists $K\in\Mod X$ such that $G(N)\oplus K$ is a submodule of $G(M)$.  Applying $F$ gives that $FG(N)\oplus F(K)$ is a submodule of $FG(M)$. Since $\ker G=\ker FG$ and $FG$ preserves Krull dimension, we have that $FG$ preserves essential containments by Corollary \ref{strong essential cor}, so that $FG(N)$ is an essential submodule of $FG(M)$.  It follows that $F(K)=0$ and, since $F$ is faithful, that $K=0$.
\end{proof}

\begin{prop} Let $X$ and $Y$ be noetherian spaces  and let $_X\M_Y=(F,G)$ be a right localizing Frobenius bimodule. Assume that $\ker G$ is closed under injective envelopes and that $F$ is faithful. If $U$ and $V$ are weakly open subspaces of $Y$, then $F^{-1}(U\cap V)=F^{-1}U\cap F^{-1}V$.\label{intersection prop}
\end{prop}

\begin{proof} If we write $\Mod U\simeq \Mod Y/S$ and $\Mod V\simeq\Mod Y/T$ for localizing subcategories $S$ and $T$ of $\Mod Y$, then we have $\Mod U\cap V\simeq \Mod Y/\Mod_{S\bullet T}Y$.  By definition, $\Mod F^{-1}U\simeq \Mod X/F^{-1}S$, $\Mod F^{-1}V\simeq \Mod X/F^{-1}T$, and $\Mod F^{-1}(U\cap V)\simeq\Mod X/F^{-1}(\Mod_{S\bullet T}Y)$.  Since $\Mod F^{-1}U\cap F^{-1}V\simeq \Mod X/\Mod_{F^{-1}S\bullet F^{-1}T}X$, we must verify that 
\[F^{-1}(\Mod_{S\bullet T}Y)=\Mod_{F^{-1}S\bullet F^{-1}T}X.\]

If $M\in\Mod_{F^{-1}S\bullet F^{-1}T}X$, then $M$ has a filtration with successive slices in either $F^{-1}S$ or $F^{-1}T$.  Since $F$ is exact, we see that $F(M)$ has a filtration with successive slices in either $S$ or $T$.  Thus  $\Mod_{F^{-1}S\bullet F^{-1}T}X\subseteq F^{-1}(\Mod_{S\bullet T}Y)$.  For the reverse containment, let $M\in F^{-1}(\Mod_{S\bullet T}Y)$ with $M$ noetherian.  Then $F(M)$ is also noetherian, and so contains an essential submodule which is a finite direct sum of uniform $Y$-modules.  Each of these $Y$-modules in turn contains a (necessarily essential) submodule in either $S$ or $T$; thus $F(N)$ contains an essential submodule of the form $K\oplus L$, with $K\in S$ and $L\in T$.  Now, $G$ preserves essential containments by the previous lemma, and so $G(K)\oplus G(L)$ is an essential submodule of $GF(M)$.  Since $F$ is faithful there is a monic $M\rightarrow GF(M)$; in particular either $M\cap G(K)$ or $M\cap G(L)$ is nonzero. Since $\M$ is right localizing $FG(K)$ and $FG(L)$ are in $S$ and $T$, respectively, and so either $F(M\cap G(K))$ is a nonzero module in $S$, or $F(M\cap G(L))$ is a nonzero module in $T$.  

We have shown that every noetherian module in $F^{-1}(\Mod_{S\bullet T}Y)$ contains a nonzero submodule in either $F^{-1}S$ or $F^{-1}T$.  From this it follows readily that $M\in\Mod_{F^{-1}S\bullet F^{-1}T}X$.  Since every $X$-module is the direct limit of its noetherian submodules this proves the reverse containment $F^{-1}(\Mod_{S\bullet T}Y)\subseteq \Mod_{F^{-1}S\bullet F^{-1}T}X$.
\end{proof}

Our main results in this section show that, over noetherian spaces, right localizing Frobenius bimodules come from geometric data.  The following is a precise formulation of this idea. 

\begin{thm}  If $_X\M_Y=(F,G)$ is a nonzero right localizing Frobenius bimodule between noetherian spaces $X$ and $Y$, then there is a surjective continuous function $f:\Supp(F)\rightarrow \Supp(G)$ such that $F(E(x))\cong E(f(x))^{(m_x)}$ for some positive integer $m_x$. \label{continuous thm}
\end{thm}

\begin{proof} Write $\M|_{Y_G}=(\bar F,\bar G)$.  By Proposition \ref{faithful prop}, $\bar F$ and $\bar G$ are faithful, and we may identify $\Inj(X_F)$ and $\Inj(Y_G)$ with $\Supp(F)$ and $\Supp(G)$, respectively.  

If $\eta:\Id_{X_F}\rightarrow \bar G\bar F$ denotes the unit of the adjoint pair $(\bar F,\bar G)$, then $\eta_M$ is monic for all $X_F$-modules $M$.  In particular, $E(x)$ is isomorphic to a summand of $\bar G\bar F(E(x))$ for all $x\in\Inj(X_F)$.  Hence there exists $y\in\Inj(Y_G)$ such that $E(x)$ is isomorphic to a summand of $\bar G(E(y))$. Applying $\bar F$ then gives that $\bar F(E(x))$ is isomorphic to a summand of $\bar F\bar G(E(y))$.  By Proposition \ref{FG prop}, applied to $\M|_{Y_G}$, we have that $\bar F\bar G E(y)\cong E(y)^{(n_y)}$ for some positive integer $n_y$.  So, $\bar F(E(x))$ is isomorphic to a summand of $E(y)^{(n_y)}$. Writing $y=f(x)$, we see that $\bar F(E(x))\cong E(f(x))^{(m_x)}$ for some positive integer $m_x$.  Given $y\in\Inj(Y_G)$, let $E(x)$ be a summand of $\bar G(E(y))$.  Then $\bar F(E(x))$ is a summand of $\bar F\bar G(E(y))\cong E(y)^{(n_y)}$, showing that $y=f(x)$.  Thus $f$ is surjective.  

To see that $f$ is continuous, recall that   $\{V(M):M\in\mod Y_G\}$ give a basis for the closed sets in $\Inj(Y_G)$.  Now, $x\in f^{-1}(V(M))$ if and only if $\Hom_{Y_G}(M,E(f(x)))\neq 0$, if and only if $\Hom_{Y_G}(M,\bar F(E(x)))\neq 0$, if and only if $\Hom_{X_F}(\bar G(M),E(x))\neq 0$, if and only if $x\in V(\bar G(M))$.   Since $\bar G(M)$ is noetherian,  $f^{-1}(V(M))=V(\bar G(M))$ is closed, proving the continuity of $f$.
\end{proof}

Under the homeomorphisms between $\Inj(X_F)$ and $\Supp(F)$ and $\Inj(Y_G)$ and $\Supp(G)$, the above result can also be phrased in terms of the existence of a surjective continuous map $f:\Inj(X_F)\rightarrow\Inj(Y_G)$.  If $Y$ is in addition an enriched space, then $f$ can be extended to a morphism of ringed spaces.  We begin with a lemma.

\begin{lemma}  Assume the hypotheses of Theorem \ref{continuous thm}, and suppose that $F$ and $G$ are both faithful.  Let $\mathfrak{U}\subseteq\Inj(Y)$ be a basic open subset, and let $U$ be the weakly open subspace of $Y$ corresponding to $\mathfrak{U}$. Then $F^{-1}U$ is the weakly open subspace of $X$ corresponding to $f^{-1}\mathfrak{U}$.
\end{lemma}

\begin{proof} Write $\mathfrak{U}=V(M)^c$, so that $\Mod U\simeq \Mod Y/T(M)$ by Lemma \ref{Gabriel topology lemma}. The proof of the continuity of $f$ in Theorem \ref{continuous thm} shows that $f^{-1}V(M)=V(G(M))$, and so $f^{-1}\mathfrak{U}=V(G(M))^c$.  If $V$ is the weakly open subspace of $X$ which corresponds to $f^{-1}\mathfrak{U}$, then $\Mod V\simeq \Mod X/T(G(M))$.  Thus it suffices to show that $T(G(M))=F^{-1}T(M)$.

By Lemma \ref{Gabriel topology lemma}, we have $T(M)=\bigcap_{y\in\mathfrak{U}}T_y$, and so $F^{-1}T(M)=\bigcap_{y\in \mathfrak{U}}F^{-1}T_y$.  Similarly we have $T(G(M))=\bigcap_{x\in f^{-1}\mathfrak{U}}T_x$.  If $E(x)$ is a isomorphic to a summand of $G(E(y))$, then $F(E(x))$ is a isomorphic to a summand of $FG(E(y))\cong E(y)^{(n_y)}$; thus $E(x)$ is a summand of $G(E(y))$ if and only if $y=f(x)$. It follows that $F^{-1}T_y=\bigcap_{f(x)=y}T_x$.  Thus $F^{-1}T(M)=\bigcap_{y\in \mathfrak{U}}\bigcap_{f(x)=y}T_x=\bigcap_{x\in f^{-1}\mathfrak{U}}T_x=T(G(M))$.
\end{proof}

\begin{prop} Assume the hypotheses of Theorem \ref{continuous thm}, and suppose further  that $Y$ is enriched, with structure module $\O_Y$.   Then the continuous function $f:\Inj(X_F)\rightarrow\Inj(Y_G)$ can be extended to a morphism of ringed spaces $f:(\Inj(X_F),\END(\bar G(\O_{Y_G})))\rightarrow (\Inj(Y_G),\END(\O_{Y_G}))$, where $\O_{Y_G}=\O_Y|_{Y_G}$.\label{morphism prop}
\end{prop}

\begin{proof} Changing notation, we assume without loss of generality that $F$ and $G$ are faithful, and that $f:\Inj(X)\rightarrow\Inj(Y)$ is surjective and continuous.  Recall that the sheaf $\END(\O_Y)$ is obtained by the rule $\mathfrak{U}\mapsto\End_U(\O_U)$ on basic open subsets $\mathfrak{U}$ of $\Inj(Y)$.  Similarly the pushforward sheaf $f_*\END(G(\O_Y))$ is obtained by the rule $\mathfrak{U}\mapsto \End_{F^{-1}U}(G(\O_Y)|_{F^{-1}U})$.  (Here we have used the previous lemma to know that the weakly open subspace of $X$ associated to $f^{-1}\mathfrak{U}$ is $F^{-1}U$.)  Using the functorial isomorphisms $G(\O_Y)|_{F^{-1}U}\cong G|_U(\O_U)$, we shall identify $\End_{F^{-1}U}(G(\O_Y)|_{F^{-1}U})$ with $\End_U(G|_U(\O_U))$.

Given a basic open subset $\mathfrak{U}$, we define $\vphi(\mathfrak{U}):\End_U(\O_U)\rightarrow\End_U(G|_U(\O_U))$ by  $\vphi(\mathfrak{U})(f)=G|_U(f)$. If $\rho_\mathfrak{V}^\mathfrak{U}$ and $\sigma_\mathfrak{V}^\mathfrak{U}$ denote the restriction homomorphisms for $\mathfrak{U}\mapsto \End_U(\O_U)$ and $\mathfrak{U}\mapsto \End_U(G|_U(\O_U))$ respectively, then it is straightforward to check that $\sigma_\mathfrak{V}^\mathfrak{U}\vphi(\mathfrak{U})=\vphi(\mathfrak{U})\rho_\mathfrak{V}^\mathfrak{U}$ whenever $\mathfrak{V}\subseteq\mathfrak{U}$ are basic open subsets of $\Inj(Y)$.  

Passing from these data on basic open sets to sheaves, we have constructed a morphism $\vphi:\END(\O_Y)\rightarrow f_*\END(G(\O_Y))$.  Thus $(f,\vphi):(\Inj(X),\END(G(\O_Y)))\rightarrow (\Inj(Y),\END(\O_Y))$ is a morphism of ringed spaces.
\end{proof}

The nest result says roughly that over certain spaces, Frobenius bimodules which come from geometric data are necessarily right localizing, and is a partial converse to Theorem \ref{continuous thm}.

\begin{prop} Let $_X\M_Y=(F,G)$ be a nonzero dimension preserving Frobenius bimodule between noetherian spaces $X$ and $Y$, and assume that $\ker G$ is closed under injective envelopes. If there exists a function $f:\Supp(F)\rightarrow \Supp(G)$ such that $F(E(x))\cong E(f(x))^{(m_x)}$ for some positive integer $m_x$, then $\M$ is right localizing, and $f$ is necessarily surjective and continuous. \label{automatically continuous} \end{prop}

\begin{proof} We must show that $T_y\subseteq G^{-1}F^{-1}T_y$ for all $y\in\Inj(Y)$, by Lemma \ref{localizing test lemma}. Since $\ker G$ is closed under injective envelopes, $y\in\Supp(G)$ if and only if $G(E(y))\neq 0$.  If $G(E(y))=0$, then $G^{-1}F^{-1}T_y=\Mod Y$, so that $T_y\subseteq G^{-1}F^{-1}T_y$ trivially.  if $G(E(y))\neq 0$, then $y\in\Supp(G)$ and so we may write $G(E(y))\cong \oplus_{i=1}^nE(x_i)$ with each $x_i\in\Supp(F)$.  By Corollary \ref{dual summand cor}, $E(y)$ is a summand of $F(E(x_i))$ for each $i$, so that $y=f(x_i)$ for all $i$.  In particular there is a positive integer $n_y$ such that $FG(E(y))\cong E(y)^{(n_y)}$.  

Now, $M\in T_y$ if and only if $\Hom_Y(M,E(y))=0$, if and only if $\Hom_Y(M,FG(E(y))=0$, if and only if $\Hom_Y(FG(M),E(y))=0$.  Thus $T_y=G^{-1}F^{-1}T_y$.  The fact that $f$ is surjective and continuous follows as in the final paragraph of Theorem \ref{continuous thm}.
\end{proof}

These results in turn enable us to characterize when certain right localizing Frobenius bimodules are left localizing.

\begin{prop} Let $X$ and $Y$ be noetherian spaces and let $\M$ be a nonzero right localizing Frobenius $X,Y$-bimodule.  If $\M$ is left localizing, then $f:\Supp(F)\rightarrow\Supp(G)$ is a homeomorphism.  If $\ker G$ is closed under injective envelopes and $\M$ is dimension preserving, then $\M$ is left localizing whenever $f$ is injective.\label{homeo prop}
\end{prop}

\begin{proof} 
Suppose that $\M$ is left localizing.  Then applying Theorem \ref{continuous thm} to the dual bundle $\M^*$ shows that there is a surjective, continuous function $g:\Supp(G)\rightarrow \Supp(F)$ such that $E(y)\otimes_Y\M^*\cong E(g(y))^{(n_y)}$ for some positive integer $n_y$.  Since $\M$ is also right localizing, we see that 
\[E(y)\otimes_Y\M^*\otimes_X\M\cong E(g(y))^{(n_y)}\otimes_X\M\cong E(fg(y))^{(m_xn_y)}.\]
Now Proposition \ref{FG prop} and the Krull Schmidt Theorem for indecomposable injectives show that $y=fg(y)$ for all $y\in\Supp(G)$.  In a similar way we see that $gf(x)=x$ for all $x\in\Supp(F)$.  Thus $g=f^{-1}$, and since $g$ is continuous $f$ is a homeomorphism.

Conversely, suppose $\ker G$ is closed under injective envelopes and that $f$ is injective. Since $f$ is automatically  surjective, it is a bijection.  Denote the inverse of $f$ by $g$.  If $y\in\Supp G$, say $y=f(x)$, then the fact that $E(f(x))$ is the only indecomposable summand of $F(E(x))$ shows that $E(x)$ is the only indecomposable summand of $G(E(y))$.  Thus $G(E(y))\cong E(g(y))^{(n_y)}$ for all $y\in\Supp G$, and Proposition \ref{automatically continuous} shows that $g$ is continuous.  Thus $\M$ is left localizing.
\end{proof}

\section{Frobenius bimodules between schemes}
In this section we consider Frobenius bimodules between $\QCoh(\O_X)$ and $\QCoh(\O_Y)$ where $X$ and $Y$ are separated noetherian schemes.  We begin by recalling some important facts from \cite[Chapitre VI]{Gabriel}.  If $X$ is a noetherian scheme, then there is a bijection between the underlying point set of $X$ and $\Inj(\QCoh(\O_X))$, defined by sending an indecomposable injective quasicoherent $\O_X$-module $E(x)$ to the generic point of the scheme-theoretic support of $E(x)$.  If $\Inj(\QCoh(\O_X))$ is endowed with its Gabriel topology, then this bijection becomes a homeomorphism, and we shall typically identify the underlying point set of $X$ with $\Inj(\QCoh(\O_X))$ under this homeomorphism.  Also, if $U$ is an open subscheme of $X$, then we can identify $\QCoh(\O_U)$ with a weakly open subspace of $\QCoh(\O_X)$; specifically, we have $\QCoh(\O_U)\simeq\QCoh(\O_X)/T(U)$, where $T(U)=\{\sh{F}\in\QCoh(\O_X):\sh{F}|_U=0\}$.

\begin{lemma} If $X$ is a noetherian scheme, then $\QCoh(\O_X)$ satisfies condition $(*)$ with respect to Krull dimension.\label{scheme star lemma}
\end{lemma}

\begin{proof}  Let $x_1,x_2\in X$, and suppose that $\vphi\in\Hom_X(E(x_1),E(x_2))$ is a nonzero morphism. The image of $\vphi$ is a nonzero subsheaf $\sh{F}$ of $E(x_2)$.  Since the support of $E(x_2)$ is an integral subscheme of with generic point $x_2$, we must have $\sh{F}_{x_2}\neq 0$. So $\vphi_{x_2}(E(x_1)_{x_2})\neq 0$, showing that $x_2\in\Supp(E(x_1))$.  It follows that $\Supp(E(x_2))\subseteq\Supp(E(x_1))$, so that the dimension of support of $E(x_1)$ is at least equal to the dimension of support of $E(x_2)$.  Since dimension of support agrees with Krull dimension for noetherian schemes, $\QCoh(\O_X)$ satisfies condition $(*)$ with respect to Krull dimension.
\end{proof}

Because we will work with actual sheaves in this section, we shall drop the notation for bimodules that we adopted above, and refer exclusively to the underlying functors, i.e. we will write a Frobenius bimodule as $(F,G)$. Also, we will refer to a Frobenius bimodule $(F,G)$ between $\QCoh(\O_X)$ and $\QCoh(\O_Y)$ as a Frobenius $X,Y$-bimodule. 

We recall the definition of a sheaf bimodule given in \cite{Nyman thesis, Nyman, Van Sklyanin, Van P1}.  Since we will work exclusively with coherent sheaf bimodules in this section, we restrict our attention to this situation. We use the definition in \cite{Nyman thesis}. (The original definition, due to Van den Bergh \cite[Definition 2.3]{Van Sklyanin}, is what we call a ``finite sheaf bimodule" here.)

\begin{defn}[{\cite[Definition 3.4]{Nyman thesis}}] Let $X$ and $Y$ be schemes, with fiber product $X\times Y$.  Then a \emph{(coherent) sheaf $X,Y$-bimodule} is a coherent $\O_{X\times Y}$-module $\sh{E}$, such that each of the morphisms $pr_1|_{\Supp \sh{E}}$ and $pr_2|_{\Supp \sh{E}}$ is affine.  Here $pr_1:X\times Y\rightarrow X$ and $pr_2:X\times Y\rightarrow Y$ denote the canonical projection morphisms.  We say that $\E$ is \emph{finite} if $pr_1|_{\Supp \sh{E}}$ and $pr_2|_{\Supp \sh{E}}$ are finite morphisms. 
\end{defn}

If $W$ is a scheme with morphisms $\a:W\rightarrow X$ and $\b:W\rightarrow Y$, then the universal property of $X\times Y$ gives a morphism $(\a,\b):W\rightarrow X\times Y$.  If $\M$ is a coherent $\O_W$-module, then we write $_\a\M_\b$ for $(\a,\b)_*\M$.  If each of the morphisms $\a$ and $\b$ is affine, then clearly $\sh{E}={_\a\M_\b}$ is a sheaf $X,Y$-bimodule.  Moreover, any sheaf $X,Y$-bimodule arises in this way, by taking $W$ to be the scheme-theoretic support of $\E$, with $\a$ and $\b$ the inclusion maps. It is clear that $\E$ is finite if and only if $\a$ and $\b$ are finite morphisms.

A sheaf $X,Y$-bimodule $\E$ defines a right exact functor $-\otimes_{\O_X}\E:\QCoh(\O_X)\rightarrow \QCoh(\O_Y)$ by the rule $\sh{F}\otimes_{\O_X}\E=pr_{2*}(pr_1^*\sh{F}\otimes_{\O_{X\times Y}}\E)$ \cite[p. 442]{Van Sklyanin}.  If $\E={_\a\M_\b}$, then one can check that $\sh{F}\otimes_{\O_X}\E=\b_*(\a^*\sh{F}\otimes_{\O_W}\M)$.  $\E$ is called \emph{locally free on the left (right)} if $pr_{1*}\E$ ($pr_{2*}\E$) is a locally free $\O_X$-module ($\O_Y$-module).   If $\E={_\a\M_\b}$ is locally free of finite rank on each side, then the left and right duals to $\E$ are the locally free sheaf $Y,X$-bimodules defined by the formulas
\begin{equation}
\begin{split}
(_X\E)^*&={_\b[\a^\uparrow\HOM_{\O_X}(\a_*\M,\O_X)]_\a}\\
(\E_Y)^*&={_\b[\b^\uparrow\HOM_{\O_Y}(\b_*\M,\O_Y)]_\a}.
\end{split}
\end{equation}
The definitions of $(_X\E)^*$ and $(\E_Y)^*$ are from \cite[Definition 3.9]{Nyman}, and we refer the reader to \cite[Section 3.1]{Nyman} for the definitions and basic properties of the functors $\a^\uparrow$ and $\b^\uparrow$.

The following proposition is a scheme-theoretic analogue of the characterization of Frobenius bimodules between rings given in section 3.

\begin{prop} Let $X$ and $Y$ be noetherian schemes, and let $\E$ be a sheaf $X,Y$-bimodule. 
\begin{enumerate}
\item If $\E$ is locally free of finite rank on each side and there is an isomorphism of sheaf $Y,X$-bimodules $(_X\E)^*\cong(\E_Y)^*$, then $-\otimes_{\O_X}\E$ is a Frobenius functor.
\item Assume that $X$ and $Y$ are smooth of the same dimension and that $\E$ is finite. If $-\otimes_X\E$ is a Frobenius functor, then $\E$ is locally free of finite rank on each side and $(_X\E)^*\cong(\E_Y)^*$ as sheaf $Y,X$-bimodules. 
\end{enumerate}
\label{sheaf Frobenius prop}
\end{prop}

\begin{proof} Let $W=\Supp(\E)$ and write $\E={_\a\M_\b}$ as above.

(1) By \cite[Proposition 3.14]{Nyman}, there are adjoint pairs $(-\otimes_{\O_X}\E,-\otimes_{\O_Y}(\E_Y)^*)$ and $(-\otimes_{\O_Y}(\E_Y)^*,-\otimes_{\O_X}((\E_Y)^*_X)^*)$.  Thus it suffices to prove that there is an isomorphism of sheaf $X,Y$-bimodules $\E\cong ((\E_Y)^*_X)^*$. By definition, we have
\begin{equation}
\begin{split}
((\E_Y)^*_X)^*&={_\a[\a^\uparrow\HOM_{\O_X}(\a_*(\M_Y)^*,\O_X)]_\b}\\
&={_\a[\a^\uparrow\HOM_{\O_X}(\a_*(\b^\uparrow\HOM_{\O_Y}(\b_*\M,\O_Y)),\O_X)]_\b}.
\end{split}\label{eq:dual}\end{equation}
By hypothesis, there is an isomorphism of sheaf $Y,X$-bimodules
\[_\b[\b^\uparrow\HOM_{\O_Y}(\b_*\M,\O_Y)]_\a\cong {_\b[\a^\uparrow\HOM_{\O_X}(\a_*\M,\O_X)]_\a},\]
which in turn implies that there is an isomorphism of coherent $\O_W$-modules
\begin{equation}\b^\uparrow\HOM_{\O_Y}(\b_*\M,\O_Y)\cong \a^\uparrow\HOM_{\O_X}(\a_*\M,\O_X).\label{eq:W iso}\end{equation}
Substituting \eqref{eq:W iso} into \eqref{eq:dual} gives
\[((\E_Y)^*_X)^*={_\a[\a^\uparrow\HOM_{\O_X}(\a_*(\a^\uparrow\HOM_{\O_X}(\a_*\M,\O_X)),\O_X)]_\b.}\]
Now, $\a_*\a^\uparrow\simeq\Id$ by \cite[Lemma 3.1]{Nyman}.  Thus we have
\[((\E_Y)^*_X)^*\cong{_\a[\a^\uparrow\HOM_{\O_X}(\HOM_{\O_X}(\a_*\M,\O_X),\O_X)]_\b}.\]
Now, $\a_*\M$ is a locally free $\O_X$-module of finite rank, so that there is an $\O_X$-module isomorphism $\HOM_{\O_X}(\HOM_{\O_X}(\a_*\M,\O_X),\O_X)\cong\a_*\M$.  Hence
\[((\E_Y)^*_X)^*\cong{_\a[\a^\uparrow\a_*\M]_\b}\cong {_\a\M_\b}\cong \E.\]

(2) Recall that $-\otimes_{\O_X}\E=pr_{2*}(pr^*_1(-)\otimes_{\O_{X\times Y}}\E)$.  It follows from \cite[Lemma 3.15]{Nyman thesis} and \cite[Proposition 2.2.7]{Van Sklyanin} that $pr_{2*}$ is exact and faithful when restricted to the image of $pr^*_1(-)\otimes_{\O_{X\times Y}}\E$.  Thus $pr^*_1(-)\otimes_{\O_{X\times Y}}\E$ is exact and takes coherent $\O_X$-modules to coherent $\O_{X\times Y}$-modules.  Let $U$ be an open affine subset of $Y$ and let $V$ be an open affine subset of $X$, so that $V\times U$ is an open affine subset of $X\times Y$.  If we write $j:V\rightarrow X$ for the inclusion, then $j$ is affine since $X$ is separated, so that $j^*$ is exact.  Since $V\times U$ is affine, taking sections over $V\times U$ is exact, and so the functor $\QCoh(\O_V)\rightarrow \Mod\O_{X\times Y}(V\times U)$ given by $\sh{F}\mapsto (pr^*_1j^*\sh{F}\otimes_{\O_{X\times Y}}\E)(V\times U)$ is exact and preserves noetherian objects.

Now, $(pr^*_1j^*\sh{F}\otimes_{\O_{X\times Y}}\E)(V\times U)=\sh{F}(V)\otimes_{\O_X(V)}M$, where $M$ is the $\O_X(V),\O_Y(U)$-bimodule $\E(V\times U)$.  Since $\QCoh(\O_V)\simeq\Mod \O_X(V)$, we conclude that the functor $-\otimes_{\O_X(V)}M:\Mod \O_X(V)\rightarrow \Mod\O_{X\times Y}(V\times U)$ is exact and preserves noetherian modules.  From this we conclude that $_{\O_X(V)}M$ is flat and finitely generated, and since $\O_X(V)$ is noetherian, that $_{\O_X(V)}M$ is finitely generated projective.  Hence $\E$ is locally free of finite rank on the left.  Since $X$ and $Y$ are smooth of the same dimension and $\E$ is finite, $\E$ is also locally free of finite rank on the right \cite[Proposition 3.1.6]{Van P1}.

Since $\E$ is locally free of finite rank on each side, we know that the left and right adjoints to $-\otimes_{\O_Y}\E$ are $-\otimes_{\O_Y}(_X\E)^*$ and $-\otimes_{\O_Y}(\E_Y)^*$, respectively. Thus there is an equivalence of functors $-\otimes_{\O_Y}(_X\E)^*\simeq -\otimes_{\O_Y}(\E_Y)^*$, and since the functor determines the sheaf bimodule up to isomorphism \cite[Lemma 3.1.1]{Van P1}, we have $(_X\E)^*\cong(\E_Y)^*$ as sheaf $Y,X$-bimodules
\end{proof}

\begin{cor} Let $X$ and $Y$ be smooth schemes of the same dimension and let $\E$ be a finite sheaf $X,Y$-bimodule, locally free of finite rank on each side.  Then $-\otimes_{\O_X}\E$ is a Frobenius $X,Y$-bimodule if and only if $\omega_X\otimes_{\O_X}\E\cong\E\otimes_{\O_Y}\omega_Y$ as sheaf $X,Y$-bimodules, where $\omega_X$ and $\omega_Y$ are the dualizing sheaves on $X$ and $Y$ respectively.
\end{cor}

\begin{proof}  According to \cite[Lemma 3.1.8]{Van P1}, there is an isomorphism of sheaf $X,Y$-bimodules $((\E_Y)^*_X)^*\cong \omega_X^{-1}\otimes_{\O_X}\E\otimes_{\O_Y}\omega_Y$.  Since $-\otimes_X\E$ is Frobenius if and only if $((\E_Y)^*_X)^*\cong\E$ as sheaf $X,Y$-bimodules, the result follows.
\end{proof}

An obvious question at this point is to what extent part (2) of Proposition \ref{Frobenius bimodule prop} carries over to the scheme-theoretic setting.  That is, if $(F,G)$ is a Frobenius bimodule between noetherian schemes $X$ and $Y$, is there a sheaf $X,Y$-bimodule $\E$ such that $F\simeq -\otimes_{\O_X}\E$?  We provide an affirmative answer in the case where $(F,G)$ is right localizing and $F$ and $G$ are faithful.  We begin with a pair of lemmas.  

\begin{lemma} Let $X$ and $Y$ be noetherian schemes, and let $(F,G)$ be a right localizing Frobenius $X,Y$-bimodule with $F$ and $G$ faithful.  If $U$ is an open subscheme of $Y$, then $F^{-1}\QCoh(\O_U)=\QCoh(\O_{f^{-1}U})$, where $f:X\rightarrow Y$ is the continuous function of Theorem \ref{continuous thm}.
\end{lemma}

\begin{proof} In the above notation, we need to show that $T(f^{-1}U)=F^{-1}T(U)$.  Let $E=\bigoplus_{y\in U}E(y)$.  Then $\Hom_Y(\sh{G},E)\cong\bigoplus_{y\in U}\Hom_Y(\sh{G},E(y))$ for all coherent $\O_Y$-modules $\sh{G}$, and we see that when $\sh{G}$ is coherent, $\sh{G}\in T(U)$ if and only if $\Hom_Y(\sh{G},E)=0$.  Since $Y$ is noetherian every $\sh{G}\in\QCoh(\O_Y)$ is the direct limit of its coherent subsheaves; from this it follows that in fact $T(U)=\{\sh{G}:\Hom_Y(\sh{G},E)=0\}$.  Similarly, if we let $E'=\bigoplus_{x\in f^{-1}U} E(x)$, then $T(f^{-1}U)=\{\sh{F}:\Hom_X(\sh{F},E')=0$.

Given $z\in X$, we have that $E(z)$ is a summand of $G(E)$ if and only if $E(f(z))$ is a summand of $E$, if and only if $z\in f^{-1}U$, if and only if $E(z)$ is a summand of $E'$.  Since 
\[F^{-1}T(U)=\{\sh{F}:\Hom_Y(F(\sh{F}),E)=0\}=\{\sh{F}:\Hom_X(\sh{F},G(E))=0\}\]
we see that $F^{-1}T(U)=T(f^{-1}U)$ as claimed.
\end{proof}

\begin{lemma} Let $(F,G)$ be a right localizing Frobenius bimodule between noetherian schemes $X$ and $Y$ and assume that $F$ and $G$ are faithful. If $U$ is an open affine subscheme of $Y$ then $f^{-1}U$ is an open affine subscheme of $X$.\label{affine lemma}
\end{lemma}

\begin{proof} By the previous lemma $(F|_U,G|_U)$ is a Frobenius $f^{-1}U,U$-bimodule, and $F|_U$ is faithful.  Thus the adjoint triple $(G|_U,F|_U,G|_U)$ forms an affine map in the sense of \cite{Rosenberg book}.  Since $U$ is affine, it follows from \cite[Proposition 6.4.1]{Rosenberg book} that $f^{-1}U$ is affine.
\end{proof}

\begin{thm} Let $X$ and $Y$ be noetherian schemes, and let $(F,G)$ be a right localizing Frobenius $X,Y$-bimodule with $F$ and $G$ faithful.  Then $F\simeq -\otimes_{\O_X}\E$, where $\E$ is a sheaf $X,Y$-bimodule, locally free of finite rank on each side, such that $(_X\E)^*\cong(\E_Y)^*$.\label{right localizing scheme thm}
\end{thm}

\begin{proof} Let $f:X\rightarrow Y$ be he continuous map of Theorem \ref{continuous thm}, and let $\{U_i:i\in I\}$ be an affine open cover of $Y$; then Lemma \ref{affine lemma} shows that $\{f^{-1}U_i:i\in I\}$ is an affine open cover of $X$.  We set $W=\bigcup_{i\in I}f^{-1}U_i\times U_i$ and view $W$ as an open subscheme of $X\times Y$. Let $\a:W\rightarrow X$ and $\beta:W\rightarrow Y$ denote the canonical maps.  We shall construct a coherent $\O_W$-module $\M$, and then set $\E={_\a\M_\b}$. 

We construct $\M$ as follows. Given $i\in I$, we set $R_i=\O_X(f^{-1}U_i)$ and $S_i=\O_Y(U_i)$.  Then there are category equivalences $\QCoh(\O_{f^{-1}U_i})\simeq \Mod R_i$, $\QCoh(\O_{U_i})\simeq \Mod S_i$, and $\QCoh(\O_{f^{-1}U_i\times U_i})\simeq\Mod R_i\otimes S_i$.  Now, for each $i\in I$ we have a Frobenius pair $(F|_{U_i},G|_{U_i})$ which we can view as a Frobenius pair between $\Mod R_i$ and $\Mod S_i$.  Thus by Proposition \ref{Frobenius bimodule prop} there is a Frobenius $R_i,S_i$-bimodule $M_i$ such that $F|_{U_i}\simeq-\otimes_RM_i$.  Viewing $M_i$ as an $R_i\otimes S_i$-module, let $\tilde M_i$ be the coherent $\O_{f^{-1}U_i\times U_i}$-module associated to $M_i$.  We show that the sheaves $\{\tilde M_i:i\in I\}$ can be glued to give a coherent sheaf $\M$ on $W$.

To ease notation, let $W_i=f^{-1}U_i\times U_i$.  Given indices $i$ and $j$, we must show that $\tilde M_i|_{W_i\cap W_j}\cong\tilde M_j|_{W_i\cap W_j}$.  Given a point $w$ in $W_i\cap W_j$, we can find an open affine subset $V$ of $U_i\cap U_j$ such that $w\in f^{-1}V\times V$, since $\{f^{-1}V\times V\}$, as $V$ ranges over the open affine subsets of $Y$, gives a basis for the topology of $W$.  Thus, it suffices to show that $\tilde M_i|_{f^{-1}V\times V}\cong \tilde M_j|_{f^{-1}V\times V}$ for all open affine $V\subseteq U_i\cap U_j$.  Since $V$ is affine,  $\tilde M_i|_{f^{-1}V\times V}$ and $\tilde M_j|_{f^{-1}V\times V}$ are sheaves associated to  $\O_X(f^{-1}V),\O_Y(V)$-bimodules $N_i$ and $N_j$, respectively.  Also, under the category equivalences $\QCoh(\O_{f^{-1}V})\simeq \Mod\O_X(f^{-1}V)$ and $\QCoh(\O_{V})\simeq \Mod\O_Y(V)$, we see that the functor $(F|_{U_i})|_V$ is naturally equivalent to $-\otimes_{\O_X(f^{-1}V)}N_i$, and the functor $(F|_{U_j})|_V$ is naturally equivalent to $-\otimes_{\O_X(f^{-1}V)}N_j$.  But $(F|_{U_i})|_V\simeq (F|_{U_j})|_V\simeq F|_V$, and so there is an equivalence of functors  $-\otimes_{\O_X(f^{-1}V)}N_i\simeq -\otimes_{\O_X(f^{-1}V)}N_j$.  From this it follows that $N_i\cong N_j$ as $\O_X(f^{-1}V),\O_Y(V)$-bimodules, so that $\tilde M_i|_{f^{-1}V\times V}\cong\tilde M_j|_{f^{-1}V\times V}$ as claimed.

In order to show that $\E={_\a\M_\b}$ is a sheaf $X,Y$-bimodule we must show that the maps $\a$ and $\b$ are affine.  But this is clear: for the open affine cover $\{f^{-1}U_i:i\in I\}$ of $X$, we have $\alpha^{-1}(f^{-1}U_i)=W_i$, which is affine.  Similarly $\b^{-1}(U_i)=W_i$ is affine.  

We next show that $F\simeq -\otimes_{\O_X}\E$, and to prove this it suffices to prove that $F|_{U_i}\simeq (-\otimes_{\O_X}\E)|_{U_i}$ for all $i$.  Given $\sh{F}\in\QCoh(\O_X)$, we have $F|_{U_i}(\sh{F}|_{f^{-1}U_i})\cong F(\sh{F})|_{U_i}$ by the projection formulas \eqref{eq:projection}. Under the category equivalences $\QCoh(\O_{f^{-1}U_i})\simeq\Mod R_i$ and $\QCoh(\O_{U_i})\simeq\Mod S_i$, $F(\sh{F})|_{U_i}$ is sent to $\sh{F}(f^{-1}U_i)\otimes_{R_i}M_i$.  On the other hand, by definition  $\sh{F}\otimes_{\O_X}\E=\b_*(\a^*\sh{F}\otimes_{\O_W}\M)$, and under the above category equivalences, $\b_*(\a^*\sh{F}\otimes_{\O_W}\M)|_{U_i}$ is sent to 
\[\b_*(\a^*\sh{F}\otimes_{\O_W}\M)(U_i)\cong (\a^*\sh{F}\otimes_{\O_W}\M)(W_i)\cong\sh{F}(f^{-1}U_i)\otimes_{R_i}M_i.\]
Thus $F\simeq -\otimes_{\O_X}\E$. 

Finally, we show that $\E$ has the stated properties.  Since $M_i=\M(W_i)$ is a Frobenius bimodule, it is finitely-generated projective on the left and on the right.  This shows that $\a_*\M$ and $\b_*\M$ are each locally free of finite rank, so that $\E$ is locally free of finite rank on each side.  By definition we have 
$(_X\E)^*={_\b[\a^\uparrow\HOM_{\O_X}(\a_*\M,\O_X)]_\a}$ and
 $(\E_Y)^*={_\b[\b^\uparrow\HOM_{\O_Y}(\b_*\M,\O_Y)]_\a}.$
Thus it suffices to show that there is an isomorphism of $\O_W$-modules
\begin{equation}
\a^\uparrow\HOM_{\O_X}(\a_*\M,\O_X)\cong {\b^\uparrow\HOM_{\O_Y}(\b_*\M,\O_Y)}.\label{dual sheaf iso}
\end{equation}
Taking sections of the left hand side over $W_i=\a^{-1}(f^{-1}U_i)$ gives 
\begin{equation}\begin{split}
\a^\uparrow\HOM_{\O_X}(\a_*\M,\O_X)(W_i)&=\HOM_{\O_X}(\a_*\M,\O_X)(f^{-1}U_i)\\
&\cong \Hom_{\O_X(f^{-1}U_i)}(\a_*\M(f^{-1}U_i),\O_X(f^{-1}U_i))\\
&\cong \Hom_{R_i}(M_i,R_i)\\
&\cong (_{R_i}M_i)^*.\end{split}\end{equation}
Similarly $\b^\uparrow\HOM_{\O_Y}(\b_*\M,\O_Y)(W_i)\cong ({M_i}_{S_i})^*$.  Since there is bimodule isomorphism $(_{R_i}M_i)^*\cong ({M_i}_{S_i})^*$ for all $i$, we see that the isomorphism \eqref{dual sheaf iso} holds.
\end{proof}

If $R$ and $S$ are rings and $M$ is a Frobenius $R,S$-bimodule, then the Endomorphism Ring Theorem \cite[Theorem 2.5]{Kadison} asserts that $E=\End_R(M^*_R)$ is a Frobenius ring extension of $S$.  The definition of Frobenius ring extension is given in \cite[Definition 1.1]{Kadison}, and is equivalent to the classical notion of Frobenius algebra when $S$ is commutative.  The following is a scheme-theoretic version of the Endomorphism Ring Theorem for right localizing bimodules.

\begin{prop} Let $(F,G)$ be a right localizing Frobenius bimodule between $X$ and $Y$, and assume that $F$ and $G$ are faithful.  Let $\A=\END_{\O_X}(pr_{2*}\E^*)$, where $F\simeq-\otimes_{\O_X}\E$.  Then $f$ extends to a morphism of ringed spaces $f:(X,\A)\rightarrow (Y,\O_Y)$, and $f_*\A$ is a sheaf of Frobenius $\O_Y$-algebras.
\end{prop}

\begin{proof} Inspecting the definition of $\END(\O_Y)$ given in section 2, we see that $\END(\O_Y)$ is determined on the basic open subsets of $Y$ by \[U\mapsto\Hom_{\O_Y(U)}(\O_Y(U),\O_Y(U))\cong\O_Y(U),\]
so that $\END(\O_Y)\cong\O_Y$ as sheaves of rings on $Y$. Now, $G(\O_Y)=\O_Y\otimes_{\O_Y}\E^*\cong pr_{2*}\E^*$, and the definition of $\END(G(\O_Y))$ is determined on a basic open subset $f^{-1}U$ of $X$ by 
\[f^{-1}U\mapsto \End_{\O_X(f^{-1}U)}(pr_{2*}\E^*(f^{-1}U)).\]
Thus $\END(G(\O_Y))\cong\END_{\O_X}(pr_{2*}\E^*)$ as sheaves of rings on $X$, and Proposition \ref{morphism prop} shows that $f$ induces a morphism of ringed spaces $f:(X,\A)\rightarrow (Y,\O_Y)$.

Let $\M^*$ be the $\O_W$-module with $_\b\M^*_\a=\E^*$.  If $U$ is an open affine subset of $Y$, then $\M^*(f^{-1}U\times U)$ is the dual to the Frobenius bimodule $\M(f^{-1}U\times U)$.  If we let $R=\O_X(f^{-1}U)$, $S=\O_Y(U)$, and $M=\M(f^{-1}U\times U)$, then we have that $\M^*(f^{-1}U\times U)\cong M^*$ as $S,R$-bimodules.  In particular, $f_*\A(U)=\A(f^{-1}U)\cong\Hom_R(M^*_R)$ is a Frobenius $S$-algebra, by the affine Endomorphism Ring Theorem.  Thus $f_*\A$ is a sheaf of Frobenius $\O_Y$-algebras.
\end{proof}

\section{Rank functions}
In this section we introduce left and right rank functions associated to a Frobenius bimodule $_X\M_Y$.  In the case that $\M$ is right localizing, we show that these rank functions behave well, giving a kind of additivity principle, and a decomposition of $X$ and $Y$ into disjoint weakly open subspaces where $\M$ has constant rank.  We assume throughout this section that $X$ and $Y$ are noetherian spaces.

\begin{defn} Let $_X\M_Y$ be a Frobenius bimodule between noetherian spaces $X$ and $Y$.  Given $x\in\Inj(X)$ and $y\in\Inj(Y)$, we define the \emph{$y$-component of the right rank of $\M$ at $x$} to be $\rrk_\M(x,y)=\kappa$, where $\kappa$ is the multiplicity with which $E(y)$ occurs as a summand of $E(x)\otimes_X\M$.  Similarly, we define the \emph{$x$-component of the left rank of $\M$ at $y$} to be $\lrk_\M(y,x)=\nu$, where $\nu$ is the multiplicity with which $E(x)$ occurs as a summand of $\HOM_Y(\M,E(y))$. Finally, we define the \emph{total right rank of $\M$ at $x$} to be $\rho_\M(x)=\sum_{y\in\Inj(Y)}\rrk_\M(x,y)$, and the \emph{total left rank of $\M$ at $y$} to be $\lambda_\M(y)=\sum_{x\in\Inj(X)}\lrk_\M(y,x)$.  
\end{defn}

The Krull-Schmidt Theorem for indecomposable injective modules shows that these notions are all well-defined.  While all of the ranks defined above may be \emph{a priori} infinite cardinals, we shall see below that they are finite in many cases of interest.  It is clear from the definitions that $\rrk_\M(x,y)=\lrk_{\M^*}(x,y)$, $\lrk_\M(y,x)=\rrk_{\M^*}(y,x)$, $\rho_\M(x)=\lambda_{\M^*}(x)$, and $\lambda_\M(y)=\rho_{\M^*}(y)$. 

\begin{lemma} Let $_X\M_Y$ is a Frobenius bimodule, and let $y,y_1,y_2\in\Inj(Y)$.  Then $\rrk_{\M^*\otimes_X\M}(y_1,y_2)=\lrk_{\M^*\otimes_X\M}(y_1,y_2)$ and $\rho_{\M^*\otimes_X\M}(y)=\lambda_{\M^*\otimes_X\M}(y)$.
\end{lemma}

\begin{proof}  These follow immediately from the definitions and the fact that $\M^*\otimes_X\M=(FG,FG)$ is self-adjoint.
\end{proof}

\begin{prop} Let $_X\M_Y=(F,G)$ be a right localizing Frobenius bimodule, and let $f:\Supp(F)\rightarrow\Supp(G)$ be the continuous function of Theorem \ref{continuous thm}.  Then, for all $y\in\Supp(G)$, we have 
\begin{equation}
\lambda_{\M^*\otimes_X\M}(y)=\sum_{y=f(x)}\lrk_\M(y,x)\rrk_\M(x,y).\label{eq:additivity}
\end{equation}
Additionally, $\rho_\M(x)$ and  $\lambda_\M(y)$ are finite for all $x\in\Supp(F)$ and $y\in\Supp(G)$.
\end{prop} 

\begin{proof} The formula clearly holds when $\M=0$, so suppose that $\M\neq 0$, and let $x\in\Supp(F)$.  Since $X$ and $Y$ are noetherian, Theorem \ref{continuous thm} shows that $F(E(x))\cong E(f(x))^{(m_x)}$ for some positive integer $m_x$.  It follows from the definitions that $m_x=\rrk_\M(x,f(x))$, and that $\rrk_\M(x,y)=0$ if $y\neq f(x)$.  This shows in particular that $\rho_\M(x)=\rrk_\M(x,f(x))$ is finite for all $x\in\Supp(F)$.

Suppose now that $y\in\Supp(G)$, so that $FG(E(y))\cong E(y)^{(n_y)}$ for some positive integer $n_y$ be Proposition \ref{FG prop}.  If we write $G(E(y))\cong\oplus_{x\in\Inj(X)} E(x)^{(\lrk_\M(y,x))}$, then $\lrk_\M(y,x)\neq 0$ only for $x\in\Supp(F)$, and in this case we have that $F(E(x))$ is a summand of $FG(E(y))\cong E(y)^{(n_y)}$.  Thus $\lrk_\M(y,x)\neq 0$ if and only if $y=f(x)$, and we have
\[E(y)^{(n_y)}\cong FG(E(y)) \cong \bigoplus_{y=f(x)}F(E(x))^{(\lrk_\M(y,x))}\cong \bigoplus_{y=f(x)}E(y)^{(\lrk_\M(y,x)\rrk_\M(x,y))}.\]
Thus $n_y=\sum_{y=f(x)}\lrk_\M(y,x)\rrk_\M(x,y)$.  Since $n_y=\lambda_{\M^*\otimes_X\M}(y)$, we have formula \eqref{eq:additivity}.  Finally, note that $\lambda_\M(y)=\sum_{x\in\Inj(X)}\lrk_\M(y,x)=\sum_{y=f(x)}\lrk_\M(y,x)$ is finite because $\sum_{y=f(x)}\lrk_\M(y,x)\rrk_\M(x,y)$ is finite.
\end{proof}

Our next result is an analogue of the well-known fact that a sheaf $\sh{L}$ on a scheme $X$ is invertible if and only if it is locally free of rank $1$.  We introduce the following notation: If $X$ is a space, then we denote by $\I(X)$ the full subcategory of $\Mod X$ consisting of the injective $X$-modules.

\begin{prop} Let $_X\M_Y=(F,G)$ be a Frobenius bimodule between noetherian spaces $X$ and $Y$, and assume that $F$ and $G$ are faithful.  Then $F:\Mod X\rightarrow\Mod Y$ is an equivalence of categories if and only if $\rho_\M(x)=1$ for all $x\in\Inj(X)$ and $\lambda_\M(y)\in\Inj(Y)$.
\end{prop}

\begin{proof} Clearly if $F$ is a category equivalence, then $\rho_\M(x)=1$ for all $x$, and since $G$ is the inverse equivalence to $F$ we see that $\lambda_\M(y)=1$ for all $y$ as well.  

Conversely, suppose the stated conditions hold, and let $f:\Inj(X)\rightarrow\Inj(Y)$ and $g:\Inj(Y)\rightarrow\Inj(X)$ be the functions defined by $F(E(x))\cong E(f(x))$ and $G(E(y))\cong E(g(y))$.  Since $F$ and $G$ are faithful the unit maps $E(x)\rightarrow E(gf(x))$ and $E(y)\rightarrow E(fg(y))$ are monic for all $x$ and $y$, and hence isomorphisms; in particular we see that $f$ and $g$ are inverse bijections.  Now, given any injective $X$-module $E$, we can write $E\cong\oplus_{x\in\Inj(X)}E(x)^{(\kappa_x)}$ for cardinals $\kappa_x$, and since $F$ and $G$ commute with sums we see that $GF(E)\cong E$ for all injective $X$-modules $E$.  Similarly we have $FG(E')\cong E'$ for all injective $Y$-modules $E'$.  Moreover, if $E_1$ and $E_2$ are injective $X$-modules, then $\Hom_Y(F(E_1),F(E_2))\cong \Hom_X(E_1,GF(E_2))\cong\Hom_X(E_1,E_2)$.  We have shown that $F:\I(X)\rightarrow\I(Y)$ is fully faithful, and that every object in $\I(Y)$ is isomorphic to $F(E)$ for some $E\in\I(X)$; thus $F$ is an equivalence of categories between $\I(X)$ and $\I(Y)$.  Since $F$ is exact, \cite[Proposition I.9.14]{Gabriel} implies that $F:\Mod X\rightarrow \Mod Y$ is an equivalence of categories.
\end{proof}

The hypothesis that $F$ and $G$ are faithful in the above proposition is a necessary one, as the following example illustrates.    

\begin{example} Let $R=\begin{pmatrix} k&0\\k&k\end{pmatrix}$ be the ring of $2\times 2$ lower triangular matrices over a field $k$, and let $X=\Mod R$.  There are up to isomorphism two simple right $R$-modules, namely $S_1=(k\ 0)$ and $S_2=(0\ k)=(k\ k)/(k\ 0)$.  Note that $E(S_1)=(k\ k)$ and $E(S_2)=S_2$.  Now, let $I$ be the ideal consisting of those matrices whose first column is possibly nonzero.  Note that $R/I\cong k$ as rings, and $R/I$ is a projective left $R$-module, isomorphic to the second column of $R$.  Thus we may view $R/I$ as an $R,k$-bimodule, and an easy computation of the duals shows that $R/I$ is Frobenius.  Note that the functor $-\otimes_RR/I$ is the same as $M\mapsto M/MI$ for a right $R$-module $M$; in particular we see that $E(S_1)\otimes_RR/I\cong E(S_2)\otimes_RR/I\cong k$ as right $k$-modules.  On the other hand, $\Hom_k(R/I,k)\cong E(S_2)$ as right $R$-modules.  Since $\Inj(\Mod R)=\{E(S_1),E(S_2)\}$ and $\Inj(\Mod k)=\{k\}$, we see that the left and right ranks of $R/I$ are identically equal to $1$; however $R/I$ clearly does not induce an equivalence of categories.
\label{triangular example 1}\qed\end{example}

Of course, $-\otimes_RR/I$ is not faithful in this example; since $I=\ann(S_1)$ we have $S_1\otimes_RR/I=0$.  The difficulty that this example poses is that localizing subcategories of $\Mod R$ need not be closed under injective envelopes; in particular if $T_2=\{M\in\Mod R:\Hom_R(M,E(S_2))=0\}$, then $S_1\in T_2$ but $E(S_1)\not\in T_2$. We shall return to this idea several more times below, and see that the failure of localizing subcategories to be closed under injective envelopes is an obstacle to carrying out many local constructions of interest.

If $\E$ is a locally free sheaf of finite rank on a scheme $X$, then it is possible to decompose $X$ as the disjoint union of open subschemes $X\cong\bigsqcup_{\lambda\in\Lambda}U_\lambda$, such that $\E|_{U_\lambda}$ has constant rank.  The following is a noncommutative version of this decomposition.  

\begin{thm} Let $_X\M_Y=(F,G)$ be a right localizing Frobenius bimodule with  $F$ and $G$  both faithful.  Let $\Lambda=\{\lambda_\M(y):y\in\Inj(Y)\}$.  Then there exist weakly open subspaces $\{U_\lambda:\lambda\in\Lambda\}$ of $Y$ such that the following hold:
\begin{enumerate}
\item $Y=\bigsqcup_{\lambda\in\Lambda}U_\lambda$ is the disjoint union of the $U_\lambda$. (That is, $Y=\bigcup_{\lambda\in\Lambda}U_\lambda$ and $U_\lambda\cap U_\mu=\emptyset$ for $\lambda\neq \mu$.)
\item $X=\bigsqcup_{\lambda\in\Lambda}F^{-1}U_\lambda$ is the disjoint union of the $F^{-1}U_\lambda$.
\item $\M|_{U_\lambda}$ is a Frobenius $F^{-1}U_\lambda,U_\lambda$-bimodule of constant total left rank $\lambda$ for all $\lambda\in\Lambda$.
\end{enumerate}\label{constant rank thm}
\end{thm}

\begin{proof}   Given $y\in\Inj(Y)$, we let $\O_y$ denote the largest critical submodule of $E(y)$. (Here we use Krull dimension on $\Mod Y$.)  Then we define $T_\lambda$ to be the smallest localizing subcategory of $\Mod Y$ containing $\{\O_y:\lambda_\M(y)\neq\lambda\}$. Thus a $Y$-module $M$ is in $T_\lambda$ if and only if the total left rank of the injective hull of any critical subquotient of $M$ is different from $\lambda$.  We then let $U_\lambda$ be the weakly open subspace of $\Mod Y$ with $\Mod U_\lambda\simeq\Mod Y/T_\lambda$.

We show that $\bigcup_{\lambda\in\Lambda}U_\lambda=Y$ and that $U_\lambda\cap U_\mu=\emptyset$ for $\lambda\neq \mu$.  For the first, we need to show that $\bigcap_{\lambda\in \Lambda} T_\lambda=0$.  If $M\in\Mod Y$ is nonzero, then let $N$ be a critical submodule of $M$, with injective hull $E(y)$.  If $\lambda_\M(y)=\lambda$, then $M\not\in T_\lambda$.  Hence the only $Y$-module contained in each $T_\lambda$ is $0$.  For the second claim, we need to show that $\Mod_{T_\lambda\bullet T_\mu}Y=\Mod Y$ for all $\lambda\neq \mu$.  This is equivalent to showing that every noetherian $Y$-module has a filtration with slices in either $T_\lambda$ or $T_\mu$.  Since every noetherian $Y$-module has a critical composition series, it suffices to show that any critical $Y$-module is in either $T_\lambda$ or $T_\mu$.  But this is clear: If $N$ is a critical $Y$-module with injective hull $E(y)$, either $\lambda_\M(y)\neq\lambda$ and $N\in T_\lambda$, or $\lambda_\M(y)\neq \mu$ and $N\in T_\mu$.  Hence we have a decomposition $Y=\bigsqcup_{\lambda\in\Lambda}U_\lambda$.

Since $\{U_\lambda:\lambda\in\Lambda\}$ is a weakly open cover for $Y$, we have by the remarks preceeding Lemma \ref{restriction lemma} that $\{F^{-1}U_\lambda:\lambda\in\Lambda\}$ is a weakly open cover for $X$.  Also, if $\lambda\neq \mu$, then $F^{-1}U_\lambda\cap F^{-1}U_\mu=F^{-1}(U_\lambda\cap U_\mu)=F^{-1}\emptyset=\emptyset$, by Proposition \ref{intersection prop}.  Thus $X=\bigsqcup_{\lambda\in\Lambda}F^{-1}U_\lambda$.

We finish by showing that $\M|_{U_\lambda}$ has the stated properties.  Write $j_\lambda:U_\lambda\rightarrow Y$ and $k_\lambda:F^{-1}U_\lambda\rightarrow X$ for the inclusions, and set $\M|_{U_\lambda}=(F_\lambda,G_\lambda)$.  Then any indecomposable injective $U_\lambda$-module is isomorphic to $j^*_\lambda E(y)$ for some $T_\lambda$-torsionfree $E(y)$.  Note that $E(y)$ is $T_\lambda$-torsionfree if and only if $\lambda_\M(y)=\lambda$, and that $G_\lambda(j^*_\lambda E(y))=k^*_\lambda G(E(y))$.  If we write $G(E(y))\cong\oplus_{i=1}^tE(x_i)^{(\lrk_\M(y,x_i))}$, then $G_\lambda(j^*_\lambda E(y))\cong \oplus_{i=1}^tk^*_\lambda E(x_i)^{(\lrk_\M(y,x_i))}$.  We claim that each $E(x_i)$ is $F^{-1}T_\lambda$-torsionfree.  Since $F(E(x_i))\cong E(y)^{(\rrk_\M(x,y))}$, we see that, if $M$ is a nonzero submodule of $E(x_i)$ in $F^{-1}T_\lambda$, then $F(M)$ is a nonzero submodule of $E(y)^{(\rrk_\M(x,y))}$ in $T_\lambda$, a contradiction.  Since each $k^*_\lambda E(x_i)$ is an indecomposable injective $F^{-1}U_\lambda$-module, we see that $\lambda_{\M|_{U_\lambda}}(y)=\sum_{i=1}^t\lrk_\M(y,x_i)=\lambda_\M(y)=\lambda$.  Thus $\M|_{U_\lambda}$ has constant total left rank $\lambda$.
\end{proof}

The decomposition of a scheme as $X=\bigsqcup_{\lambda\in\Lambda}U_\lambda$ as a union of pairwise disjoint open subschemes implies that there is a category equivalence $\QCoh(\O_X)\simeq\bigoplus_{\lambda\in\Lambda}\QCoh(\O_{U_\lambda})$.  In the noncommutative case, this need no longer be true:  It can happen that $X$ is the disjoint union of two weakly open (even open) subspaces $U_1$ and $U_2$ without $\Mod X$ being equivalent to $\Mod U_1\times\Mod U_2$.  The following concrete example illustrates this.  

\begin{example} Again let $R$ be the ring of $2\times 2$ lower triangular matrices.  We retain the notation of Example \ref{triangular example 1}.  For $i=1,2$, let $T_i=\{M\in\Mod R:\Hom_R(M,E(S_i))=0\}$, and let $U_i$ be the weakly open subspace of $X=\Mod R$ with $\Mod U_i\simeq \Mod R/T_i$.  Note that each $U_i$ is in fact open, since $U_i$ is the open complement to $\Mod R/I_i$, where $I_i$ is the two-sided ideal of $R$ consisting of those matrices whose $(i,i)$-entry is $0$ (for $i=1,2$).

We show that $U_1$ and $U_2$ are disjoint, and that their union is $X$.  For the first claim, we must show that $\Mod_{T_1\bullet T_2}X=\Mod X$.  Since $S_2\in T_1$ and $S_1\in T_2$ and every noetherian $R$-module has finite length this is clear.   Similarly, we have that $T_1\cap T_2=0$, since $E(S_1)\oplus E(S_2)$ is an injective cogenerator for $\Mod R$.  

However, we claim that $\Mod R$ is \emph{not} equivalent to $\Mod U_1\times\Mod U_2$.  To see this, note that we have $\Mod U_1\simeq\Mod U_2\simeq\Mod k$, but $\Mod R\not\simeq\Mod k\times\Mod k$: the exact sequence $0\rightarrow S_1\rightarrow E(S_1)\rightarrow S_2\rightarrow 0$ gives a nonsplit extension of $S_1$ by $S_2$ in $\Mod R$, while the category $\Mod k\times\Mod k$ is semisimple.\qed
\end{example}

As noted above, the difficulty in Example \ref{triangular example 1} is that $T_2$ is not closed under injective envelopes.  It turns out that this is precisely the obstacle to $\Mod X$ being equivalent to $\Mod U_1\times\Mod U_2$, which we now prove in greater generality. 
 
\begin{prop} Let $X$ be a noetherian space, and let $\{U_i:i\in I\}$ be a collection of pairwise disjoint weakly open subspaces of $X$ whose union is $X$.  Write $\Mod U_i\simeq\Mod X/T_i$ for localizing subcategories $T_i$ of $\Mod X$.  If each $T_i$ is closed under injective envelopes, then there is a category equivalence $\Mod X\simeq\bigoplus_{i\in I}\Mod U_i$, given by $M\mapsto(M|_{U_i})_{i\in I}$ and $f\mapsto (f|_{U_i})_{i\in I}$. \label{category decomposition prop}
\end{prop}

\begin{proof} Let $j_i:U_i\rightarrow X$ denote the inclusion.  Since $X$ is noetherian $j_{i*}$ commutes with direct sums for all $i$.  We first show that, given $x\in\Inj(X)$, there is a unique $i\in I$ with $j^*_iE(x)\neq 0$.  If we denote the torsion functor for $T_i$ by $\tau_i$, then $j^*_iE(x)=0$ if and only if $\tau_iE(x)=E(x)$. (Here we use that each $T_i$ is closed under injective envelopes.) Since $\bigcap_{i\in I}T_i=0$ by hypothesis, we see that $j^*_iE(x)\neq 0$ for some $i\in I$.  If $k\neq i$, then $U_i\cap U_k=\emptyset$ implies that $\Mod_{T_i\bullet T_k}X=\Mod X$.  Consequently, if $\tau_iE(x)=0$, then we must have $\tau_kE(x)\neq 0$, which implies that $\tau_kE(x)=E(x)$.  Hence $j^*_iE(x)\neq 0$ for a unique $i\in I$.

Given $E\in\I(X)$, there exist cardinals $\{\kappa_x:x\in\Inj(X)\}$ such that $E\cong\bigoplus_{x\in\Inj(X)}E(x)^{(\kappa_x)}$.  Given $i\in I$, let $\Sigma_i=\{x\in\Inj(X):j^*_iE(x)\neq 0\}$. Since each of $j^*_i$ and $j_{i*}$ commute with sums and $j_{i*}j^*_iE(x)\cong E(x)$ for $x\in\Sigma_i$, we see that $j_{i*}j^*_iE\cong\bigoplus_{x\in\Sigma_i}E(x)^{(\kappa_x)}$, and since $j^*_iE(x)$ nonzero for a single $i\in I$, we see that $\{\Sigma_i:i\in I\}$ actually partitions $\Inj(X)$.  Combining these we see that 
\[\bigoplus_{i\in I}j_{i*}j^*_iE\cong\bigoplus_{i\in I}\bigoplus_{x\in\Sigma_i}E(x)^{(\kappa_x)}\cong\bigoplus_{x\in\Inj(X)}E(x)^{(\kappa_x)}\cong E.\]
Moreover, given an injective $U_i$-module $Q_i$, there exists a $T_i$-torsionfree injective  $X$-module $E_i$ such that $j^*_iE_i\cong Q_i$, and $E_i$ is necessarily $T_k$-torsion for $k\neq i$.  From this it follows that, given $(Q_i)_{i\in I}$ in $\bigoplus_{i\in I}\I(U_i)$, there exists $E\in\I(X)$ with $(j^*_iE)_{i\in I}\cong (Q_i)_{i\in I}$, namely $E=\bigoplus_{i\in I}E_i$. 

Finally, given $E_i,E_2\in\I(X)$, we have
\begin{equation}
\begin{split}
\Hom_X(E_1,E_2)&\cong \Hom_X(\bigoplus_{i\in I}j_{i*}j^*_iE_1,E_2)\\
&\cong\prod_{i\in I}\Hom_X(j_{i*}j^*_iE_1,E_2)\\
&\cong\prod_{i\in I}\Hom_{U_i}(j^*_iE_1,j^*_iE_2)\\
&\cong\Hom_{\bigoplus U_i}((j^*_iE_1)_{i\in I},(j^*_iE_2)_{i\in I})
\end{split}
\end{equation}
Hence the functor $E\mapsto (j^*_iE)_{i\in I}$ from $\I(X)$ to $\bigoplus_{i\in I}\I(U_i)$ is fully faithful. Since every object of $\bigoplus_{i\in I}\I(U_i)$ is isomorphic to $(j^*_iE)_{i\in I}$ for some $E\in\I(X)$, we conclude that $E\mapsto (j^*_iE)_{i\in I}$ gives an equivalence of categories $\I(X)\simeq \bigoplus_{i\in I}\I(U_i)$.  Since $E\mapsto (j^*_iE)_{i\in I}$ is exact, it follows from \cite[Proposition I.9.14]{Gabriel} that it determines an equivalence of categories $\Mod X\simeq \bigoplus_{i\in I}\Mod U_i$.
\end{proof}

The above proposition leads to the following corollary to Theorem \ref{constant rank thm}.

\begin{cor} Assume the hypotheses and notation of Theorem \ref{constant rank thm}. Assume that $\M$ is dimension preserving and that each $T_\lambda$ is closed under injective envelopes.  Then we have decompositions $\Mod Y\simeq\bigoplus_{\lambda\in\Lambda}\Mod U_\lambda$, $\Mod X\simeq\bigoplus_{\lambda\in\Lambda}\Mod F^{-1}U_\lambda$, and $\M\cong \bigoplus_{\lambda\in\Lambda}\M|_{U_\lambda}$.\label{constant rank cor}
\end{cor}

\begin{proof} The only thing that needs to be established is that each $F^{-1}T_\lambda$ is closed under injective envelopes.  Since $\M$ is dimension preserving and $F$ is faithful this follows by Corollary \ref{strong essential cor}.
\end{proof}

In general the set $\Lambda$ of allowable total left ranks need not be finite, even in the situation of Corollary \ref{constant rank cor}.  For example, we can take $\Mod X\simeq\oplus_{n\in \N}\Mod k$ for a field $k$, and we can define a Frobenius $X,X$-bimodule $\M$ componentwise by $\M_n=\Id_k^{(n)}$.  However, in the presence of certain finiteness assumptions on $\Mod Y$, it is possible to prove that $|\Lambda|$ is finite.

\begin{prop}  Assume the hypotheses of Theorem \ref{constant rank thm}, and assume that each $T_\lambda$ is closed under injective envelopes.  Assume that there exists an injective $Y$-module $E$ such that $E$ is isomorphic to a finite direct sum of indecomposable injective $Y$-modules, and such that $E$ subgenerates $\Mod Y$.  Then $\Lambda$ is finite.
\end{prop}

\begin{proof} Write $E\cong\oplus_{i=1}^tE(y_i)^{(n_i)}$ for positive integers $n_i$.  For each $i$, the proof of Proposition \ref{category decomposition prop} shows that there is a single index $\lambda_i$ such that $E(y_i)|_{U_{\lambda_i}}\neq 0$.  Consequently $E|_{U_\lambda}\neq 0$ if and only if $\lambda\in\{\lambda_1,\dots,\lambda_t\}$.  Since $(-)|_{U_\lambda}$ is exact and commutes with sums for all $\lambda$, and $E$ subgenerates $\Mod Y$, we see that $N|_{U_\lambda}\neq 0$ if and only if $\lambda\in\{\lambda_1,\dots,\lambda_t\}$ for all $Y$-modules $N$.  It follows that $\Mod Y\simeq\oplus_{i=1}^t\Mod U_{\lambda_i}$ and $|\Lambda|$ is finite.
\end{proof}

\begin{remark}  The existence of an injective $Y$-module satisfying the above hypotheses is ensured in the following cases: $Y$ is affine, $Y$ is a scheme admitting an ample line bundle, or $Y$ is an integral space in the sense of \cite[Definition 3.1]{Smith integral}.

For the first, note that if $\Mod Y\simeq\Mod R$, then $R$ is a right noetherian ring, and $E(R)$ satisfies the hypotheses of the proposition.  For the second, note that $\{\O_Y(-n):n\in \N\}$ is a set of generators for $\QCoh(\O_Y)$, where $\O_Y(1)$ denotes an ample line bundle on $Y$.  Since $\O_Y(-n)$ is isomorphic to a subsheaf of $\O_Y$ for all $n\in\N$, we see that $\O_Y$ is a noetherian subgenerator for $\QCoh(\O_Y)$, and so $E(\O_Y)$ satisfies the hypotheses of the proposition.  Finally, part of the definition of ``integral space"  in \cite {Smith integral} is the existence of an indecomposable injective $Y$-module $E$ which subgenerates $\Mod Y$. 
\end{remark}

\section{Gluing right localizing bimodules}
Since right localizing Frobenius bimodules have good local properties, one might reasonably expect to be able to glue right localizing bundles on weakly open covers.  More precisely, we can ask the following question:  Suppose $X$ and $Y$ are spaces with weakly open covers $\{V_i:i\in I\}$ and $\{U_i:i\in I\}$ respectively, and that there are right localizing Frobenius $V_i,U_i$-bimodules $\M_i$ for $i\in I$ such that $\M_i|_{U_i\cap U_j}\cong \M_j|_{U_i\cap U_j}$ for all $i,j$.  Does there exist a right localizing Frobenius $X,Y$-bimodule $\M$ such that $\M|_{U_i}\cong\M_i$?  We shall see below that the answer is ``yes" in some cases, provided that the index set $I$ is finite.  However, the next example reveals that without additional hypotheses, the answer is ``no" even in the case where $Y=X$ and $|I|=2$.

\begin{example} Let $R$, $U_1$, and $U_2$ be as in Example \ref{triangular example 1}.  We define Frobenius $U_i,U_i$-bimodules $\M_i$ for $i=1,2$ by letting $\M_1=\Id_{U_1}$ and $\M_2=0$.  Clearly $\M_1$ and $\M_2$ are both right localizing, and since $U_1\cap U_2=\emptyset$ there is no overlap condition to verify.  We claim that there can be no Frobenius $X,X$-bimodule with $\M|_{U_i}=\M_i$.

Suppose such an $\M$ existed.  Then there would be nonnegative integers $m,n$ such that $E(S_1)\otimes_X\M\cong E(S_1)^{(m)}\oplus E(S_2)^{(n)}$.  If $j_1:U_1\rightarrow X$ denotes the inclusion, then $j^*_1(E(S_1)\otimes_X\M)\cong j^*_1E(S_1)\otimes_{U_1}\M_1=j^*_1E(S_1)$, so that $m=1$.  On the other hand, if $j_2:U_2\rightarrow X$ denotes the inclusion, then $j^*_2(E(S_1)\otimes_X\M)\cong j^*_2E(S_1)\otimes_{U_2}\M_2=0$.  But $j^*_2E(S_1)$ is a nonzero summand of $j^*_2(E(S_1)\otimes_X\M)$.\qed
\end{example}

As in the previous section, the obstacle to things working smoothly is the fact that $T_2$ is not closed under injective envelopes.  The purpose of this section is to prove that, under suitable hypotheses, this is the only obstacle to gluing over finite covers.  We need to fix some notation and hypotheses, which will remain in force for the rest of the section.

\begin{notation} $X$ and $Y$ will be noetherian spaces, and $\{V_1,V_2\}$ and $\{U_1,U_2\}$ will be weakly open covers for $X$ and $Y$, respectively.  We write $\Mod V_i\simeq \Mod X/S_i$ for localizing categories $S_i$, and similarly we write $\Mod Y\simeq \Mod Y/T_i$ for localizing subcategories $T_i$, for $i=1,2$.  For ease of notation, we set $S_{12}=\Mod_{S_1\bullet S_2}X$ and $T_{12}=\Mod_{T_1\bullet T_2}Y$, so that $\Mod V_1\cap V_2\simeq\Mod X/S_{12}$ and $U_1\cap U_2\simeq \Mod Y/T_{12}$.  We assume that $S_i$ and $T_i$ are closed under injective envelopes for $i=1,2$.

Also, we let $k_i:V_i\rightarrow X$ and $j_i:U_i\rightarrow Y$ denote the inclusions for $i=1,2$, and we write $k_{12}:V_1\cap V_2\rightarrow X$ and $j_{12}:U_1\cap U_2\rightarrow Y$ for the inclusions.  By \cite[Lemma 6.12]{Smith subspaces}, there are weakly open immersions $\b_i:V_1\cap V_2\rightarrow V_i$ and $\a_i:U_1\cap U_2\rightarrow U_i$ such that $k_{12}=k_i\b_i$ and $j_{12}=j_i\a_i$ for $i=1,2$.\label{notation}
\end{notation}

The following is our main gluing result for right localizing bimodules.

\begin{thm}  Keep the above notation, and let $\M_i=(F_i,G_i)$ be a right localizing Frobenius $V_i,U_i$-bimodule, with $F_i$ and $G_i$ faithful, for $i=1,2$.  Assume further that $F_1^{-1}(U_1\cap U_2)=F_2^{-1}(U_1\cap U_2)= V_1\cap V_2$.  If $\M_1|_{U_1\cap U_2}=\M_2|_{U_1\cap U_2}$, then there exists a right localizing Frobenius $X,Y$-bimodule $\M=(F,G)$ such that each of $F$ and $G$ are faithful, and $\M_i\cong\M|_{U_i}$ for $i=1,2$.\label{gluing thm}
\end{thm}

Before proving Theorem \ref{gluing thm}, we establish two technical lemmas.

\begin{lemma} $S_{12}$ is closed under injective envelopes.  \label{closed injective lemma}
\end{lemma}

\begin{proof} Any $X$-module $M$ in $S_{12}$  has a filtration with successive slices in either $S_1$ or $S_2$.  Suppose that $M$ is noetherian and $M\in S_{12}$.  We prove that $E(M)\in S_{12}$ by induction on the smallest length of a filtration on $M$ with slices in $S_1$ or $S_2$.  Denote this length by $l$.

If $l=1$, then either $M\in S_1$ or $M\in S_2$; since each of $S_1$ and $S_2$ is closed under injective envelopes we have in either case that $E(M)\in S_{12}$.  Now suppose $l>1$, and let $S_1$ be the first term in a filtration of $M$ of length $l$.  Then either $M_1\in S_1$ or $M_1\in S_2$.  Either way we have $E(M_1)\in S_{12}$.  Also, since $M/M_1$ has a filtration of length strictly less than $l$, we have $E(M/M_1)\in S_{12}$ by induction.  Now the exact sequence
\[0\rightarrow E(M_1)\rightarrow E(M)\rightarrow E(M/M_1)\]
shows that $E(M)\in S_{12}$ also.  Since $X$ is noetherian every $X$-module is the direct limit of its noetherian submodules, and since $S_{12}$ is closed under direct limits it follows readily that $E(M)\in S_{12}$ for all $M\in S_{12}$.
\end{proof}

\begin{lemma} Keep the above notation. Then every injective $X$-module $E$ is isomorphic to a direct sum $E_1\oplus E_2\oplus Q$ such that the following hold:
\begin{enumerate}
\item $E_1$ and $E_2$ are in $S_{12}$, and $Q$ is $S_{12}$-torsionfree.
\item $E_1$ is $S_1$-torsionfree and $S_2$-torsion, and $E_2$ is $S_2$-torsionfree and $S_1$-torsion.
\end{enumerate}
Moreover, any such decomposition is unique up to isomorphism.
\label{injective decomposition lemma}
\end{lemma}

\begin{proof}  Denote the torsion functors associated to $S_1$, $S_2$, and $S_{12}$ by $\sigma_1$, $\sigma_2$, and $\sigma_{12}$ respectively.  Let $E$ be an injective $X$-module.  Since $S_{12}$ is closed under injective envelopes, $\sigma_{12}E$ is injective, and necessarily isomorphic to a summand of $E$.  Thus $E\cong \sigma_{12}E\oplus Q$, where $Q$ is $S_{12}$-torsionfree.  

Let $E_1=\sigma_2(\sigma_{12}E)$ and $E_2=\sigma_1(\sigma_{12}E)$.  Since $\{V_1,V_2\}$ is a weakly open cover for $X$, $S_1\cap S_2=0$.  From this it follows that $E_1$ is $S_1$-torsionfree and $S_2$-torsion, and similarly $E_2$ is $S_2$-torsionfree and $S_1$-torsion.  Also, since each of $S_1$ and $S_2$ are closed under injective envelopes, each of $E_1$ and $E_2$ is an injective summand of $\sigma_{12}E$.  Now, the fact that $\sigma_{12}E$ is $S_{12}$-torsion and $S_{12}=\Mod_{S_1\bullet S_2}X$ show that in fact $\sigma_{12}E\cong E_1\oplus E_2$.  Thus $E$ has a decomposition as $E_1\oplus E_2\oplus Q$ with the stated properties.  The uniqueness statement is clear.
\end{proof}

Note that analogous statements to Lemmas \ref{closed injective lemma} and \ref{injective decomposition lemma} hold, \emph{mutatis mutandis}, over $Y$.

\begin{proof}[Proof of Theorem \ref{gluing thm}]  We divide the proof into several steps for the convenience of the reader.  Recall that $\I(X)$ denotes the full subcategory of injective $X$-modules.  We shall identify $\I(V_i)$ with the full subcategory of $\I(X)$ consisting of the $S_i$-torsionfree injective $X$-modules, for $i=1,2$, and similarly for $\I(V_{12})$ and the $S_{12}$-torsionfree injective $X$-modules.  We also make similar identifications for $\I(U_1)$, $\I(U_2)$, and $\I(U_{12})$.  Finally, we write $\M_1|_{U_1\cap U_2}=\M_2|_{U_1\cap U_2}=(F_{12},G_{12})$.

\subsubsection*{Step 1}  Let $Q\in\I(V_{12})$.  Then, viewing $Q$ as an $S_{12}$-torsionfree injective $X$-module as above,  we establish natural isomorphisms
\begin{equation}
j_{12*}F_{12}k_{12}^*(Q)\cong j_{1*}F_1k^*_1(Q)\cong j_{2*}F_2k^*_2(Q).\label{eq:torsionfree injective}
\end{equation}
Since $F_{12}\b^*_1\simeq\a_1^*F_1$, $j_{12*}=j_{1*}\a_{1*}$ and $k_{12}^*=\b^*_1k^*_1$, there are natural isomorphisms
\[j_{12*}F_{12}k^*_{12}(Q)\cong j_{1*}\a_{1*}F_{12}\b^*_1k^*_1(Q)\cong j_{1*}\a_{1*}\a_1^*F_1k^*_1(Q)\]
for all $Q\in\Mod X$.  Since $Q$ is $S_{12}$-torsionfree, $F_1k_1^*(Q)$ is $T_{12}/T_1$-torsionfree.  Since there is a natural isomorphism $\a_{1*}\a^*_1E\cong E$ whenever $E\in\I(U_1)$ is $T_{12}/T_1$-torsionfree, we obtain the first isomorphism in \eqref{eq:torsionfree injective}.  The second isomorphism is similar.

Similarly, given $Q'\in\I(U_{12})$, there are natural isomorphisms
\[k_{12*}G_{12}j_{12}^*(Q')\cong k_{1*}G_1j^*_1(Q')\cong k_{2*}G_2j^*_2(Q').\]

\subsubsection*{Step 2:}  We define a functor $F:\I(X)\rightarrow \I(Y)$ as follows.  Given $E$ in $\I(X)$, we fix an isomorphism $E\cong E_1\oplus E_2\oplus Q$ as in Lemma \ref{injective decomposition lemma}.  Then we define $F$ on objects by \[F(E)=j_{1*}F_1k^*_1(E_1)\oplus j_{2*}F_2k^*_2(E_2)\oplus j_{12*}F_{12}k^*_{12}(Q).\]
We need to define $F$ on morphisms.  To do so, let $E'\in\I(X)$ with decomposition $E'\cong E'_1\oplus E'_2\oplus Q'$.  Then the conditions (1)--(3) of Lemma \ref{injective decomposition lemma} imply that 
\begin{multline}
\Hom_X(E,E')\cong \Hom_X(E_1,E'_1)\oplus\Hom_X(E_2,E'_2)\\
\oplus\Hom_X(Q,E'_1)\oplus\Hom_X(Q,E'_2)\oplus\Hom_X(Q,Q').\label{Hom decomposition equation}\end{multline} 
Given $f\in\Hom_X(E,E')$, we can decompose $f$ as $f_1+\dots+f_5$, corresponding to the five terms (in order) given in \eqref{Hom decomposition equation}.  We show how to determine $F(f_i)$ in each of the corresponding terms in the decomposition of $\Hom_X(F(E),F(E'))$.

Since $E_1$ is $S_1$-torsionfree, we have $E_1\cong k_{1*}k_1^*E_1$, so that 
$\Hom_X(E_1,E'_1)\cong \Hom_{V_1}(k^*_1E_1,k^*_1E'_1)$.  So, we define $F(f_1)$ to be the image of $f_1$ under the following maps:
\begin{multline}
\Hom_X(E_1,E_1')\xrightarrow{\cong} \Hom_{V_1}(k^*_1E_1,k^*_1E_1')\xrightarrow{F_1}\\
\Hom_{U_1}(F_1k_1^*(E_1),F_1k_1^*(E'_1))
\xrightarrow{j_{1*}}\Hom_Y(j_{1*}F_1k_1^*(E_1),j_{1*}F_1k_1^*(E'_1)).\end{multline}
Similarly, we define $F(f_2)$ and $F(f_5)$, using $k_2$ and $F_2$ and $k_{12}$ and $F_{12}$ respectively.

We proceed to define $F(f_3)$, the definition for $F(f_4)$ being analogous.  By Step 1, we have a natural isomorphism $j_{12*}F_{12}k^*_{12}(Q)\cong j_{1*}F_1k_1^*(Q)$.  Since $E_1'$ and $Q$ are both $T_1$-torsionfree, we have as above that $\Hom_X(Q,E_1')\cong\Hom_{V_1}(k^*_1Q,k^*_1E'_1)$.  Thus we define $F(f_3)$ as the image of $f_3$ under the maps
\begin{multline}
\Hom_X(Q,E_1')\xrightarrow{\cong} \Hom_{V_1}(k^*_1Q,k^*_1E_1')
\xrightarrow{F_1}\Hom_{U_1}(F_1k^*_1(Q),F_1k_1^*(E'_1))\\
\xrightarrow{j_{1*}}\Hom_Y(j_{1*}F_1k^*_1(Q),j_{1*}F_1k^*_1(E_1'))
\xrightarrow{\cong}
\Hom_Y(j_{12*}F_{12}k^*_{12}(Q),j_{1*}F_{1}k^*_{1}(E_1')).\end{multline}

The verification that $F(fg)=F(f)F(g)$ for $f:E'\rightarrow E''$, $g:E\rightarrow E'$ is routine but tedious, using the various natural isomorphisms defined above, and therefore left to the industrious reader.

 Similar constructions, using the functors $G_1$, $G_2$, and $G_{12}$, give a functor $G:\I(Y)\rightarrow\I(X)$.  In particular, given $E'\in\I(Y)$, we can write $E'\cong E_1'\oplus E_2'\oplus Q'$ by Lemma \ref{injective decomposition lemma}, and we have \[G(E')=k_{1*}G_1j^*_1(E'_1)\oplus k_{2*}G_2j^*_2(E'_2)\oplus k_{12*}G_{12}j^*_{12}(Q').\]

\subsubsection*{Step 3:} We verify that $F:\I(X)\rightarrow\I(Y)$ is both a left and right adjoint to $G$.  Choose $E\in\I(X)$ and $E'\in\I(Y)$ with decompositions $E\cong E_1\oplus E_2\oplus Q$ and  $E'\cong E'_1\oplus E'_2\oplus Q'$ as in Lemma \ref{injective decomposition lemma}. Then we have as in \eqref{Hom decomposition equation} that

\begin{multline}
\Hom_X(E,G(E'))\cong \Hom_X(E_1,k_{1*}G_1j^*_1(E_1'))\oplus\Hom_X(E_2,k_{2*}G_2j^*_2(E_2'))\\
\oplus\Hom_X(Q,k_{1*}G_1j^*_1(E_1'))\oplus\Hom_X(Q,k_{2*}G_2j^*_2(E_2'))\\
\oplus\Hom_X(Q,k_{12*}G_{12}j^*_{12}(Q')).\end{multline}
We treat each term individually.  For the first,
\begin{equation}
\begin{split}
\Hom_X(E_1,k_{1*}G_1j^*_1(E_1'))&\cong \Hom_{V_1}(k^*_1E_1,G_1j^*_1(E_1'))\\
&\cong \Hom_{U_1}(F_1k^*_1(E_1),j^*_1E_1')\\
&\cong \Hom_Y(j_{1*}F_1k^*_1(E_1),j_{1*}j^*_1E_1')\\
&\cong \Hom_Y(j_{1*}F_1k^*_1(E_1),E_1'), 
\end{split}
\end{equation}
where we have used that natural isomorphism $I\cong j_{1*}j^*_1I$ for $T_1$-torsionfree injective $Y$-modules $I$. Similar calculations give the isomorphisms \[\Hom_X(E_2,k_{2*}G_2j^*_2(E_2'))\cong \Hom_X(j_{2*}F_2k^*_2(E_2),E_2')\] 
and 
\[\Hom_X(Q,k_{12*}G_{12}j^*_{12}(Q'))\cong \Hom_X(j_{12*}F_{12}k^*_{12}(Q),Q').\]

For the remaining terms, we have
\begin{equation}
\begin{split}
\Hom_X(Q,k_{1*}G_1j^*_1(E_1'))&\cong \Hom_{V_1}(k_1^*Q,G_1j^*_1E_1')\\
&\cong \Hom_{U_1}(F_1k^*_1(Q),j^*_1E_1')\\
&\cong \Hom_Y(j_{1*}F_1k^*_1(Q),j_{1*}j^*_1E'_1)\\
&\cong \Hom_Y(j_{12*}F_{12}k^*_{12}(Q),j_{1*}j_1^*E'_1)\\
&\cong \Hom_Y(j_{12*}F_{12}k^*_{12}(Q),E'_1),
\end{split}
\end{equation}
where we have used \eqref{eq:torsionfree injective} and the natural isomorphism $I\cong j_{1*}j^*_1I$ for $T_1$-torsionfree injective $Y$-modules $I$.  A similar calculation yields the formula $\Hom_X(Q,k_{2*}G_2j^*_2(E_2'))\cong  \Hom_Y(j_{12*}F_{12}k^*_{12}(Q),E'_2)$.

Combining the above five isomorphisms gives the desired isomorphism $\Hom_X(E,G(E'))\cong\Hom_Y(F(E),E')$ for all $E\in\I(X)$ and $E'\in\I(Y)$.  A similar calculation shows that $\Hom_Y(E',F(E))\cong \Hom_X(G(E'),E)$ for all $E\in\I(X)$ and $E'\in\I(Y)$.  Thus, $G$ is both a left and right adjoint to $F$, viewed as functors between $\I(X)$ and $\I(Y)$.

\subsubsection*{Step 4:}  We extend $F$ and $G$ to functors between $\Mod X$ and $\Mod Y$ in the usual way. Since we need the explicit description of this extension on modules, we briefly recall it.  Given $M\in\Mod X$, we write  $M$ as the kernel of a map of injective $X$-modules: $0\rightarrow M\xrightarrow{f}E\rightarrow E'$. Then $F(M)$ is defined to be the kernel of the morphism $F(f):F(E)\rightarrow F(E')$.  (This requires choosing one such kernel for each morphism.)  Similarly we define $G:\Mod Y\rightarrow \Mod X$.

We now verify that $G$ is both a left and right adjoint to $F$. Given $M\in\Mod X$ and $N\in\Mod Y$, we choose injective resolutions $0\rightarrow M\xrightarrow{i} E_1\xrightarrow{\vphi} E_2$ and $0\rightarrow N\xrightarrow{j} E_1'\xrightarrow{\psi} E_2'$.  Given $f\in\Hom_Y(F(M),N)$, we can find $\a\in\Hom_Y(F(E_1),E_1')$ and $\b\in\Hom_Y(F(E_2),E_2')$ such that the following diagram is commutative, with exact rows:
\begin{equation}
\begin{CD}
0@>>>F(M) @>F(i) >> F(E_1) @>F(\vphi) >>F(E_2)\\
& & @VVf V @VV\alpha V @VV\beta V\\
0@>>>N @>j >> E_1' @>\psi >>E_2'.
\end{CD}
\end{equation}

Using the fact that $G$ is a right adjoint to $F$ between $\I(X)$ and $\I(Y)$, we can find $\a':\Hom_X(E_1,G(E_1'))$ and $\b'\in\Hom_X(E_2,G(E_2'))$ so that the following diagram is commutative, with exact rows:
\begin{equation}
\begin{CD}
0@>>>M @>i >> E_1 @>\vphi >>E_2\\
& & & & @VV\alpha' V @VV\beta' V\\
0@>>>G(N) @>G(j) >> G(E_1') @>G(\psi) >>G(E_2').
\end{CD}
\end{equation}
Now, $G(\psi)\circ\alpha'\circ i=\beta'\circ\vphi\circ i=0$, so that there exists a unique map $g:M\rightarrow G(N)$ such that $j\circ g =\alpha'\circ i$.  We have therefore constructed a map $\Phi:\Hom_Y(F(M),N)\rightarrow \Hom_X(M,G(N))$.  In a similar way we can construct a map $\Psi:\Hom_Y(M,G(N))\rightarrow \Hom_Y(F(M),N)$.  We leave to the reader the verification that $\Phi\circ \Psi$ and $\Psi\circ\Phi$ are the identity maps.  It follows that $(F,G)$ is an adjoint pair between $\Mod X$ and $\Mod Y$.  Interchanging $F$ and $G$ shows that $(G,F)$ is also an adjoint pair.  Thus $_X\M_Y=(F,G)$ is a Frobenius $X,Y$-bimodule.  

\subsubsection*{Step 5:}  We conclude the proof by showing that each of $F$ and $G$ is faithful, that $\M$ is right localizing, and that $\M|_{U_i}\cong \M_i$ for $i=1,2$.  To see that $F$ is faithful, fix $M\in\Mod X$ and let $E(x)$ be an indecomposable injective summand of $E(M)$.  Then either $E(x)$ is $S_1$ torsion and $S_2$ torsionfree, or $S_2$ torsion and $S_1$-torsionfree, or $S_{12}$-torsionfree, by Lemma \ref{injective decomposition lemma}.  Using the faithfulness of $G_1$ and $G_2$ and the fact that they agree on the overlap, we see that there exists $E(y)\in\Inj(Y)$, which is either $T_1$ torsion and $T_2$ torsionfree, or $T_2$ torsion and $T_1$-torsionfree, or $T_{12}$-torsionfree, such that $E(x)$ is a summand of $k_{1*}G_1j^*_1(E(y))$, or $k_{2*}G_1j^*_2(E(y))$, or $k_{12*}G_{12}j^*_{12}(E(y))$ respectively.  In any case we have that $E(x)$ is isomorphic to a summand of $G(E(y))$.  Since $\Hom_X(M,G(E(y)))\neq 0$, we have that $\Hom_Y(F(M),E(y))\neq 0$ and so $F$ takes nonzero $X$-modules to nonzero $Y$-modules.  Hence $F$ is faithful, and a similar proof shows that $G$ is faithful.

To show that $\M$ is right localizing, it suffices by Lemma \ref{localizing test lemma} to show that $G^{-1}F^{-1}T_y=T_y$ for all $y\in\Inj(Y)$.  (We have equality instead of containment because $F$ and $G$ are faithful.)  Given $y\in \Inj(Y)$, there are three possibilities for $E(y)$:  it is either $T_{12}$-torsionfree, $T_1$-torsionfree and $T_2$-torsion, or $T_2$-torsionfree and $T_1$-torsion.  We prove that $G^{-1}F^{-1}T_y=T_y$ in the case where $E(y)$ is $T_1$-torsionfree and $T_2$-torsion, the other two cases being similar.  

By definition, we have $FG(E(y))\cong Fk_{1*}G_1j^*_1(E(y))$.  Since $E(y)$ is $T_1$-torsionfree and $T_2$-torsion, we have that $k_{1*}G_1j^*_1(E(y))$ is $S_1$-torsionfree and $S_2$-torsion.  Hence we have $FG(E(y))\cong j_{1*}F_1k^*_1k_{1*}G_1j^*_1(E(y))\cong j_{1*}F_1G_1j^*_1(E(y)).$  Since $j^*_1E(y)$ is an indecomposable injective $U_1$-module and $\M_1=(F_1,G_1)$ is right localizing, we have $F_1G_1j^*_1(E(y))\cong j^*_1E(y)^{(n_y)}$ for some positive integer $n_y$, by Proposition \ref{FG prop}.  Thus $FG(E(y))\cong j_{1*}j^*_1E(y)^{(n_y)}\cong E(y)^{(n_y)}$.  Since $G^{-1}F^{-1}T_y=\{M\in\Mod Y:\Hom_Y(M,FG(E(y)))=0\}$, we see that $G^{-1}F^{-1}T_y=T_y$ as claimed.

Finally, we show that $\M|_{U_i}=\M_i$, and again we only prove it for $i=1$, the proof for $i=2$ being similar.  From Lemma \ref{restriction lemma} and the equivalence $k_1^*k_{1*}\simeq\Id_{V_1}$ we have that $F|_{U_1}\simeq j_1^*Fk_{1*}$.  Now, if $E\in\I(V_1)$, then there is an injective $X$-module $E'$ with $k^*_1E'\cong E$, and in the decomposition of Lemma \ref{injective decomposition lemma}, we have $E'\cong E_1\oplus Q$; that is, the $S_2$-torsionfree and $S_1$-torsion term is zero.  Hence
\[F|_{U_1}(E)\cong j_1^*F(E_1\oplus Q)\cong j_1^*j_{1*}F_1k^*_1(E_1)\oplus j_1^*j_{12*}F_{12}k^*_{12}(Q).\]
Using \eqref{eq:torsionfree injective} and $j^*_1j_{1*}\simeq \Id_{U_1}$, we can write this as 
\[F|_{U_1}(E)\cong j_1^*j_{1*}F_1k^*_1(E_1)\oplus j_1^*j_{1*}F_{1}k^*_{1}(Q)\cong F_1k^*_1(E_1\oplus Q).\]
Since $k^*_1(E_1\oplus Q)\cong E$, we see that $F|_{U_1}(E)\cong F_1(E)$ for all $E\in\I(V_1)$.  A lengthy calculation using the definition of $F|_{U_1}$ and $F_1$ on morphisms shows that in fact $F_1\simeq F|_{U_1}$ as functors from $\I(V_1)$ to $\I(U_1)$.  Since both functors are exact, it follows that $F_1\simeq F|_{U_1}$ as functors from $\Mod V_1$ and $\Mod U_1$, and uniqueness of adjoints up to isomorphism show that $G_1\simeq G|_{U_1}$ as well.
\end{proof}

\begin{remark} By induction we can extend Theorem \ref{gluing thm} to finite weakly open covers satisfying the appropriate generalization of Notation \ref{notation}. The fact that the cover must be finite, and that $X$ and $Y$ must be noetherian, are hypotheses which are made necessary by the technical limitations of the proof employed.  We do not know whether it is possible to glue infinitely many bimodules, or whether the noetherian hypotheses can be relaxed.
\end{remark}

\section{Frobenius $X,X$-bimodules}
In this section we consider Frobenius $X,X$-bimodules for a space $X$, focusing in particular on two special classes of bimodules.  

\subsection{Localizing bimodules} We begin by imposing a local condition on a Frobenius $X,X$-bimodule $\M$ that enables us to restrict $\M$ to a Frobenius $U,U$-bimodule for each weakly open subspace $U$ of $X$.  The precise definition is as follows.

\begin{defn} A Frobenius bimodule $_X\M_X=(F,G)$ is called \emph{localizing} if $F^{-1}T\subseteq T$ and $G^{-1}T\subseteq T$ for all localizing subcategories $T$ of $\Mod X$.
\end{defn}

If $_X\M_X=(F,G)$ is a localizing bimodule, then we have that $G^{-1}F^{-1}T\subseteq F^{-1}T\subseteq T$ for all localizing subcategories $T$, and similarly $F^{-1}G^{-1}T\subseteq T$.  Thus $\M$ is both left and right localizing.  If $X$ is noetherian, then Proposition \ref{homeo prop} shows that there exists a homeomorphism $f:\Supp(F)\rightarrow\Supp(G)$ such that $F(E(x))\cong E(f(x))^{\rho_{\M}(x)}$.  We shall show that in fact $f$ is the identity map.  We begin with a lemma.

\begin{lemma}  Let $X$ be a noetherian space.  If $_X\M_X=(F,G)$ is localizing, then $\M$ preserves every dimension function on $\Mod X$.  
\end{lemma}

\begin{proof} Fix a dimension function $d$, and let $T_\a=\{M\in\Mod X:d(M)\leq\a\}$. If $d(M)=\a$, then $F(M)$ and $G(M)$ are each in $T_\a$, showing that $d(F(M))\leq d(M)$ and $d(G(M))\leq d(M)$ for all $X$-modules $M$.  Now suppose that $M$ is $\a$-critical and $F(M)\neq 0$, and let $N$ be a nonzero submodule of $F(M)$.  Since $\Hom_X(N,F(M))\neq 0$, we have $\Hom_X(G(N),M)\neq 0$.  Since $N\leq F(M)$, we have $d(G(N))\leq d(N)\leq d(F(M))\leq d(M)$.  On the other hand, the image of a nonzero $f:G(N)\rightarrow M$ has dimension $\a$, so that $d(G(N))\geq d(M)$.   Since every noetherian $X$-module has a critical composition series and $F$ is exact, we see that $d(F(M))=d(M)$ whenever $F(M)\neq 0$ and $M\in\mod X$.  Since $X$ is noetherian we see that $F$ is dimension preserving, and the proof for $G$ is similar.  
\end{proof}

\begin{prop} Let $_X\M_X=(F,G)$ be a localizing Frobenius bimodule over a noetherian space $X$.  Then $F(E(x))\cong E(x)^{(\rho_{\M}(x))}$ for all $x\in\Supp(F)$.  In particular, $\Supp(F)=\Supp(G)$.
\end{prop}

\begin{proof}  Fix $x\in\Supp(F)$.  Given $M\in T_x$, we have that $G(M)\in T_x$ because $\M$ is localizing.  Hence 
\[0=\Hom_X(G(M),E(x))\cong\Hom_X(M,F(E(x)))\cong\Hom_X(M,E(f(x)))^{(\rho_\M(x))}.\]  
Thus $M\in T_x$ implies that $M\in T_{f(x)}$.  In particular, if $N$ is a critical submodule of $E(f(x))$ (with respect to Krull dimension), it follows that $\Hom_X(N,E(x))\neq 0$.  Since $E(x)$ and $E(f(x))$ have the same critical dimension by Theorem \ref{decomp thm}, it follows that $N$ is isomorphic to a submodule of $E(x)$; i.e. $E(x)\cong E(f(x))$.  Thus $f$ is the identity map on $\Supp(F)$.  From this it follows immediately that $\Supp(F)=\Supp(G)$.
\end{proof}

The above result shows that localizing is \emph{not} the same as left and right localizing:  Any autoequivalence of $\Mod X$ which induces a nontrivial map $f$ on $\Inj(X)$ will be left and right localizing, but not localizing.

If $X$ is a noetherian space and $_X\M_X=(F,G)$ is a localizing Frobenius bimodule, then we have $F(E(x))\cong E(x)^{(\rho_\M(x))}$ and $G(E(x))\cong E(x)^{(\lambda_\M(x))}$ for all $x\in\Supp(F)$.  In fact, the left and right ranks of $\M$ are equal:

\begin{prop}  Let $_X\M_Y=(F,G)$ be a localizing Frobenius bimodule over a noetherian space $X$.  Then $\rho_\M(x)=\lambda_\M(x)$ for all $x\in\Supp(F)$.
\end{prop}

\begin{proof} Given $M,N\in\Mod X$, the adjunction isomorphism $\Hom_X(F(M),N)\cong \Hom_X(M,G(N))$ is actually an isomorphism of right $\End_X(M)$-bimodules, where $\End_X(M)$ acts on $\Hom_X(F(M),N)$ via the ring homomorphism $\End_X(M)\rightarrow\End_X(F(M))$.  

Suppose now that $x\in\Supp(F)$.  Then there are isomorphisms as right $\End_X(E(x))$-modules:
\begin{multline} \End_X(E(x))^{(\rho_\M(x))}\cong \Hom_X(E(x),E(x))^{(\rho_\M(x))}\cong\Hom_X(F(E(x)),E(x))\\
\cong \Hom_X(E(x),G(E(x)))\cong \Hom_X(E(x),E(x))^({\lambda_\M(x)})\cong \End_X(E(x))^{(\lambda_\M(x))}.\end{multline}
Since $\End_X(E(x))$ is local, it has the invariant basis property; hence $\rho_\M(x)=\lambda_\M(x)$.
\end{proof}

Let $U$ be a weakly open subspace of $\Mod X$ with  $\Mod U\simeq\Mod X/T$. Then since $T\subseteq F^{-1}T$, we see that $F^{-1}U\subseteq U$.  In particular Lemma \ref{restriction lemma} gives a $U,U$-bimodule $\bar \M=(\bar F,\bar G)$ such that $\bar Fj^*_U=j^*_UF$.  A slight modification of the proof of Proposition \ref{local Frobenius prop} shows that in fact $\bar \M$ is a Frobenius $U,U$-bimodule, and that $\bar Gj^*_U=j^*_UG$.  We shall denote $\bar \M$  by $\M|_U$, even though this notation is at odds with our previous usage. (Above, $\M|_U$ denoted an $F^{-1}U,U$-bimodule, while here it denotes a $U,U$-bimodule.) One reason that we make this notational change is that the ``new" $\M|_U$ is a more natural object of study in this context. For instance, we have the following.

\begin{lemma} If $\M$ is a localizing Frobenius $X,X$-bimodule, then $\M|_U$ is a localizing Frobenius $U,U$-bimodule for every weakly open subspace $U$ of $X$.\label{localizing restriction lemma}
\end{lemma}

\begin{proof}  Writing $\Mod U\simeq\Mod X/T$ for a localizing subcategory $T$ of $\Mod X$, the localizing subcategories of $\Mod U$ are identified with $S/T$, where $S$ is a localizing subcategory of $\Mod X$ containing $T$.  The result now follows by the fact that $\M$ is localizing, and the formulas $F|_Uj^*_U=j^*_UF$ and $G|_Uj^*_U=j^*_UG$, where $\M|_U=(F|_U,G|_U)$.
\end{proof}

 The following lemma shows that the notational ambiguity in $\M|_U$ vanishes when $F$ is faithful. 

\begin{lemma} Let $_X\M_X=(F,G)$ be a localizing Frobenius bimodule over a noetherian space $X$, and assume that $F$ is faithful.  Then $G$ is also faithful, and  $T=F^{-1}T=G^{-1}T$ for all localizing subcategories $T$ of $\Mod X$.  Consequently $F^{-1}U=G^{-1}U=U$ for all weakly open subspaces $U$ of $X$.
\end{lemma}

\begin{proof} Since $F$ is faithful, $\Supp(F)=\Inj(X)$. Since $\Supp(G)=\Supp(F)$, we also have that $\Supp(G)=\Inj(X)$, and hence $G$ is also faithful.  It suffices to show that $T_x=F^{-1}T_x=G^{-1}T_x$ for all $x\in\Inj(X)$.  Since $F(E(x))\cong E(x)^{(\rho_\M(x))}$ for all $x\in\Inj(X)$, we have that $\Hom_X(M,E(x))=0$ if and only if $\Hom_X(M,F(E(x))=0$, if and only if $\Hom_X(G(M),E(x))=0$.  Thus $T_x=G^{-1}T_x$, and similarly $T_x=F^{-1}T_x$. 
\end{proof}

It is possible to modify the proof of Theorem \ref{gluing thm} to prove a gluing theorem for localizing bimodules over a noetherian space $X$, using the bimodules $\M|_U$ introduced in this section:  

\begin{thm} Let $X$ be a noetherian space with weakly open cover $\{U_1,\dots, U_n\}$.  Write $\Mod U_i\simeq\Mod X/T_i$, and assume that each $T_i$ is closed under injective envelopes.  If $\M_i$ is a localizing Frobenius $U_i,U_i$-bimodule for each $i$, such that $\M_i|_{U_i\cap U_j}=\M_j|_{U_i\cap U_j}$ for all $i$ and $j$, then there exits a localizing Frobenius $X,X$-bimodule $\M$ with $\M|_{U_i}\cong \M_i$ for all $i$.
\end{thm}

The details of the proof are left to the reader. 

\subsection{Centralizing bimodules} Localizing Frobenius bimodules are close to being ``commutative," in the sense that the induced map $f$ on the support of $F$ is the identity.  There is another type of commutativity condition that can be imposed on Frobenius $X,X$-bimodules; namely, if $_X\M_X=(F,G)$, one can require that the action of $F$ and/or $G$ commutes with the center of $\Mod X$. Recall that the \emph{center} of an abelian category is the ring of natural transformations of the identity functor.  If $X$ is a space, we write $\Cent(X)$ for the center of $\Mod X$. If $X\simeq \Mod R$ is affine, then there is an isomorphism $\Cent(\Mod R)\cong\Cent(R)$, given by sending $z\in\Cent(R)$ to the natural transformation $\tau(z)$ given by multiplication by $z$.

\begin{defn} An Frobenius bimodule $_X\M_X=(F,G)$ is \emph{centralizing} if $F$ commutes with the center of $X$; that is, if $F(\tau_M)=\tau_{F(M)}$ for all $\tau\in\Cent(X)$ and all $M\in \Mod X$.
\end{defn}

In contrast to most of our previous definitions, the definition of centralizing is one-sided.  However the following result shows that in fact, being centralizing is a two-sided condition.

\begin{lemma} If $\M$ is centralizing then so is $\M^*$.\label{central dual lemma}
\end{lemma}

\begin{proof}  Writing $\M=(F,G)$, we need to show that $G(\tau_M)=\tau_{G(M)}$ for all $\tau\in\Cent(X)$.  Let $\varepsilon:FG\rightarrow \Id_X$ denote the counit of the adjoint pair $(F,G)$, and let $\nu:\Hom_X(-,G(-))\rightarrow\Hom_X(F(-),-)$ denote the corresponding natural transformation of bifunctors.  If $\tau\in\Cent(X)$, then $\tau_M\varepsilon_M=\varepsilon_M\tau_{FG(M)}=\varepsilon_MF(\tau_{G(M)})$.  Basic properties of the natural transformation $\nu$ imply that $\tau_M\varepsilon_M=\nu(G(\tau_M))$ and $\varepsilon_MF(\tau_{G(M)})=\nu(\tau_{G(M)})$.  Since $\nu$ is a bijection we obtain $G(\tau_M)=\tau_{G(M)}$ as claimed.
\end{proof}

\begin{prop}  Let $R$ be a ring, and let $_RM_R$ be a Frobenius bimodule.  
\begin{enumerate}
\item $M$ is centralizing if and only if $rm=mr$ for all $m\in M$ and $r\in\Cent(R)$.
\item If $R$ is commutative and $M$ is centralizing, then $M$ is localizing.
\end{enumerate}\label{central=localizing prop}
\end{prop}

\begin{proof} (1) Given $N\in\Mod R$, we have $F(N)=N\otimes_RM$.  If $n\otimes m$ is a basic tensor in $F(N)$, then $\tau(z)_{F(N)}(n\otimes m)=n\otimes mz$.  On the other hand, $F(\tau(z)_N)(n\otimes m)=\tau(z)_N\otimes_RM(n\otimes m)=nz\otimes m=n\otimes zm$.  Thus $M$ is centralizing if and only if $n\otimes zm=n\otimes mz$ for all $N\in\Mod R$, $m\in M$, and $z\in \Cent(R)$.  Clearly, if $zm=mz$ for all $z\in\Cent(R)$ and $m\in M$, then $M$ is centralizing.  Applying the above characterization to the case $N=R$ shows that if $M$ is centralizing, then $zm=mz$ for all $z\in \Cent(R)$ and all $m\in M$.

(2) Since $M$ is centralizing, the left and right actions of $R$ on $M$ coincide by part (1).  In particular $M$ is isomorphic to a bimodule direct summand of $R^{(n)}$ for some $n$.  This implies that $F(N)=N\otimes_RM$ is isomorphic to a summand of $N^{(n)}$ for all $N\in\Mod R$, and similarly $G(N)=\Hom_R(M,N)$ is isomorphic to a summand of $N^{(n)}$ for all $N\in\Mod R$.   From this it follows readily that $M$ is localizing.
\end{proof}

Part (2) of the above proposition is special to commutative rings.  Indeed, for many spaces $X$ it is possible for an $X,X$-bimodule to be centralizing without being localizing.  The following is a noncommutative affine example.

\begin{example} Let $R$ be the free algebra $k\langle x,y\rangle$ over a field $k$. Then $k=\Cent(R)$, and a Frobenius bimodule $M$ is centralizing if and only if the left and right $k$ actions are the same. Let $\vphi$ be the automorphism of $R$ interchanging $x$ and $y$, and let $M={_1R}_\vphi$. (That is, the right of action of $R$ on $M$ is given by $m\cdot r=m\vphi(r)$.) Then $M$ is clearly centralizing, but is not localizing.  To see this, note that $R/xR\otimes_RM\cong R/yR$, so that $E(R/xR)\otimes_RM\cong E(R/yR)\not\cong E(R/xR)$.\qed
\end{example}

It is also interesting to note that the analogue of part (2) of Proposition \ref{central=localizing prop} need not hold for a noetherian scheme $X$.  To see this, first note that the center of $\QCoh(\O_X)$ is isomorphic to the ring of global sections $\Gamma(X,\O_X)$:  Given a global section $s$, we define a natural transformation $\tau(s)$, where $\tau(s)_\sh{F}:\sh{F}\rightarrow\sh{F}$ is induced by multiplication by $s$.  Conversely, if $\{U_i:i\in I\}$ is an open affine cover of $X$ and $\tau\in\Cent(\QCoh(\O_X))$, then $\tau$ restricts to an element of $\Cent(\Mod\O_X(U_i))$ for all $i$.  Thus $\tau_\sh{F}$ is given locally by multiplication by a section $s_i\in\O_X(U_i)$; these sections satisfy the necessary compatibilities to be lifted to a global section of $\O_X$.  

Now, let $k$ be a field and let $X=\mathbb{P}^1_k$.  If $\vphi$ is the automorphism of $k[x,y]$ interchanging $x$ and $y$, then $\vphi$ induces an autoequivalence of $\QCoh(\O_X)$, which commutes with the center, because $\Gamma(X,\O_X)=k$.  However, this autoequivalence is not a localizing bimodule, since it induces a nontrivial automorphism on the underlying point set of $X$.

If $_X\M_X$ is a localizing Frobenius bimodule that is also centralizing, one may ask whether or not $\M|_U$ is centralizing for a weakly open subspace $U$ of $X$.  In general the answer to this question is ``no," because $\Cent(U)$ can be much larger than $\Cent(X)$.  The following example illustrates this phenomenon.

\begin{example}  Let $K/k$ be an extension of fields, and let $\vphi$ be a nontrivial $k$-automorphism of $K$.  Consider the ring $R=\begin{pmatrix} k&0\\K&K\end{pmatrix}$.  Then $\vphi$ induces an automorphism of $R$, by acting componentwise.  Let $M$ be the Frobenius bimodule ${_1R_\vphi}$, so that $M^*={_1R_{\vphi^{-1}}}$.   Note first that the center of $R$ is $k$, embedded as the diagonal matrices, and since the action of $\vphi$ on $k$ is trivial we see by Proposition \ref{central=localizing prop}(1) that $M$ is centralizing.  

Note also that there are two simple $R$-modules up to isomorphism; namely $S_1=(k\ 0)$ and $S_2=(0\ K)=(K\ K)/(K\ 0)$. Here $R$ acts on $S_1$ via its $(1,1)$-component $k$ and $R$ acts on $S_2$ via its $(2,2)$-component $K$, and so we write $S_1\cong k$ and $S_2\cong K$. 
There are exactly three nonzero localizing subcategories of $\Mod R$: $T_1$, consisting of direct sums of copies of $k$, $T_2$, consisting of direct sums of copies of $K$, and $\Mod R$.

Now, $k\otimes_RM\cong k_{\vphi}\cong k$ and $K\otimes_RM\cong K_\vphi\cong K$, and similarly $k\otimes_RM^*\cong k$ and $K\otimes_RM^*\cong K$, so that $M$ is seen to be a localizing bimodule.  If $U$ is the weakly open subspace with $\Mod U\simeq \Mod R/T_1$, then $\Mod U\simeq\Mod K$, and under this equivalence $M|_U\cong {_1K_\vphi}$.  Since $\Cent(U)\cong K$ and $\vphi$ is a nontrivial automorphism of $K$, we see that $M|_U$ is not centralizing.\qed
\end{example}

In light of the above example it is natural to identify those localizing bimodules $_X\M_X$ such that $\M|_U$ is centralizing for all $U$.

\begin{defn} A localizing Frobenius bimodule $_X\M_X$ is called \emph{locally centralizing} if $\M|_U$ is centralizing for all weakly open subspaces $U$ of $X$. (Here $\M|_U$ denotes the $U,U$-bimodule defined earlier in this section.)
\end{defn}

\begin{prop} Let $\M$ be a localizing bimodule over a noetherian space $X$.  Then $\M$ is locally centralizing if and only if $\M_x$ is centralizing for all $x\in\Inj(X)$, where $\M_x=\M|_{X_x}$.
\end{prop}

\begin{proof} One direction is clear, so assume that $\M_x$ is centralizing for all $x\in\Inj(X)$.  Given a weakly open subspace $U$ of $\Mod X$, write $\Mod U\simeq \Mod X/T$ for a localizing subcategory $T$.  By Lemma \ref{localizing subcategory lemma}, we can write $T=\bigcap_{x\in\Sigma}T_x$ for some $\Sigma\subseteq\Inj(X)$.  If we let $j:U\rightarrow X$ denote the inclusion, then there are weakly open immersions $\a_x:X_x\rightarrow U$ such that $j\a_x=j_x$, for all $x\in\Sigma$ \cite[Lemma 6.12]{Smith subspaces}. Since $\{X_x:x\in\Sigma\}$ is a weakly open cover for $U$, the functor $\Phi:\Mod U\rightarrow\prod_{x\in\Sigma}\Mod X_x$ given by $\Phi(M)=(\a_x^*M)_{x\in\Sigma}$ is faithful.

Let $\M_x=(F_x,G_x)$ and $\M|_U=(F|_U,G|_U)$, so that $F_xj^*_x=j^*_xF$ and $F|_Uj^*=j^*F$.  Then the fact that $j\a_x=j_x$ implies that $F_x\a_x^*j^*=\a_xj^*F=\a^*_xF|_Uj^*$, so that $F_x\a_x^*=\a_x^*F|_U$.  Given $\tau\in\Cent(U)$, we define $\tau(x)\in\Cent(X_x)$ by the rule $\tau(x)_N=\a^*_x(\tau_M)$, where $M\in\Mod U$ satisfies $\a^*_xM=N$.  (The fact that $\tau(x)\in\Cent(X_x)$ is left to the reader.)  Then we have 
\begin{multline}\Phi(F|_U(\tau_M))=(\a^*_xF|_U(\tau_M))_{x\in\Sigma}=(F_x\a_x^*(\tau_M))_{x\in\Sigma}=\\(\tau(x)_{F_x\a_x^*(M)})_{x\in \Sigma}=(\tau(x)_{\a^*_xF|_U(M)})_{x\in\Sigma}=\Phi(\tau_{F|_U(M)}).\end{multline}
Since $\Phi$ is faithful we conclude that $\tau_{F|_U(M)}=F|_U(\tau_M)$.  Thus $\M|_U$ is centralizing.
\end{proof}

We can use the techniques of section 6 to give a characterization of those locally centralizing bimodules over a noetherian scheme $X$.  It turns out that the locally centralizing bimodules are precisely those that come from locally free sheaves of finite rank on $X$, lending credence to the idea that locally centralizing Frobenius bimodules are ``commutative."

\begin{prop} Let $X$ be a noetherian scheme.  Then there is a one-to-one correspondence between isomorphism classes of locally centralizing Frobenius bimodules on $X$ and isomorphism classes of locally free $\O_X$-modules of finite rank, given by $(F,G)\mapsto F(\O_X)$ and $\E\mapsto(-\otimes_{\O_X}\E,\HOM_{\O_X}(\E,-))$.  
\end{prop}

\begin{proof}  The proof follows very closely the proof of Theorem \ref{right localizing scheme thm}.  Let $(F,G)$ be a locally centralizing Frobenius bimodule on $X$, and let $\{U_i:i\in I\}$ be an open affine cover for $X$.  Then $(F|_{U_i},G|_{U_i})$ is a centralizing Frobenius bimodule on $\QCoh(\O_{U_i})\simeq\Mod \O_X(U_i)$.  If we set $R_i=\O_X(U_i)$, then there is a central Frobenius $R_i$-bimodule $M_i$ such that $F|_{U_i}\simeq-\otimes_{\O_{U_i}}\tilde M_i$, where $\tilde M_i$ denotes the sheaf of $\O_{U_i}$-modules associated to $M_i$.

As in the proof of Theorem \ref{right localizing scheme thm}, the $\tilde M_i$ can be glued to give a coherent sheaf $\E$ on $X$.  The key observation is that, since the left and right actions of each $R_i$ agree on each $M_i$, we can view $\E$ as an ordinary sheaf on $X$ and not as a sheaf $X,X$-bimodule.  The fact that $F\simeq -\otimes_{\O_X}\E$ again follows as in Theorem \ref{right localizing scheme thm}, where now we can use the ordinary tensor product of $\O_X$-modules, instead of the bimodule tensor product of section 6.  Since $F(\O_X)\cong\O_X\otimes_{\O_X}\E\cong \E$, we see that $(F,G)\mapsto F(\O_X)$ assigns a locally free $\O_X$-module of finite rank to each locally centralizing Frobenius bimodule on $X$, and isomorphic bimodules are sent to isomorphic locally free sheaves.  

On the other hand, the functorial isomorphisms $\sh{F}\otimes_{\O_X}\E\cong\HOM_{\O_X}(\E^*,\sh{F})$ and $\sh{F}\otimes_{\O_X}\E^*\cong\HOM_{\O_X}(\E,\sh{F})$ for a quasicoherent $\O_X$-module $\sh{F}$ show that $(-\otimes_{\O_X}\E,\HOM_{\O_X}(\E,-))$ is a Frobenius bimodule on $X$, which is clearly locally centralizing.  The fact that $\E$ determines $-\otimes_{\O_X}\E$ up to isomorphism establishes the bijection.  
\end{proof}

\section{The category $\Frob(X,Y)$}
In this section we study $\Frob(X,Y)$, the category of all Frobenius $X,Y$-bimodules, and several of its full subcategories.  We view $\Frob(X,Y)$ as a full subcategory of $\BIMOD(X,Y)$, so that a morphism $\tau:\M\rightarrow\sh{N}$ of Frobenius $X,Y$-bimodules is a natural transformation $\tau:\HOM_Y(\sh{N},-)\rightarrow \HOM_Y(\M,-)$.  Note that if $\M$ and $\sh{N}$ are in $\BIMOD(X,Y)$, then $\Hom(\M,\sh{N})$ may be a proper class.  However, if $\M$ and $\sh{N}$ are Frobenius bimodules, then \cite[Lemma 5.1]{CGN} shows that $\Hom(\M,\sh{N})$ is actually a set.  

One difficulty that we need to address is the fact that the adjoint to a Frobenius functor $G:\Mod Y\rightarrow \Mod X$ is only determined up to natural equivalence.  So we let $FP(X,Y)$ denote the category of Frobenius pairs $(F,G)$, where $F:\Mod X\rightarrow \Mod Y$ and $G:\Mod Y\rightarrow \Mod X$.  A morphism of pairs $\tau:(F_1,G_1)\rightarrow(F_2,G_2)$ is a natural transformation $\tau:G_2\rightarrow G_1$.  Then the forgetful functor $FP(X,Y)\rightarrow \Frob(X,Y)$ sending $(F,G)$ to $G$ is an equivalence of categories.  (Since $\Frob(X,Y)$ has a set of isomorphism classes and each $G\in\Frob(X,Y)$ has a set of left adjoints, we choose one left adjoint for each isomorphism class and use this to define a functor $\Frob(X,Y)\rightarrow FP(X,Y)$.)

Given Frobenius pairs $(F_1,G_1)$, $(F_2,G_2)$, and $(F_3,G_3)$ in $FP(X,Y)$, let $\eta^i:\Id_X\rightarrow G_iF_i$ and $\varepsilon^i:F_iG_i\rightarrow \Id_Y$ denote the unit and counit for the adjoint pair $(F_i,G_i)$, and let $\theta^i:\Id_Y\rightarrow F_iG_i$ and $\xi^i:G_iF_i\rightarrow \Id_X$ denote the unit and counit for the adjoint pair $(G_i,F_i)$, for each of $i=1,2,3$.  We define contravariant functors $(-)^*:FP(X,Y)\rightarrow FP(Y,X)$ and $(-)^\dagger:FP(X,Y)\rightarrow FP(Y,X)$, as follows.

On objects, we set $(F,G)^*=(F,G)^\dagger=(G,F)$.  If $\tau:(F_1,G_1)\rightarrow (F_2,G_2)$ is a morphism, then we define $\tau^*$ and $\tau^\dagger$ by the formulas:
\begin{equation}
\begin{split}
\tau_M^*&=\varepsilon^1_{F_2(M)}F_1(\tau_{F_2(M)}\eta^2_M)\\
\tau^\dagger_M&=F_2(\xi^1_M\tau_{F_1(M)})\theta^2_{F_1(M)}.
\end{split}
\end{equation}

\begin{prop} Each of $(-)^*$ and $(-)^\dagger$ is a contravariant functor.
\end{prop}

\begin{proof} These  follow from standard properties of the various units and counits involved.   We first show that $\tau^*:F_1\rightarrow F_2$ is a natural transformation.  Given $f:M\rightarrow N$, we have
\begin{equation}
\begin{split}
\tau^*_NF_1(f)&=\varepsilon^1_{F_2(N)}F_1(\tau_{F_2(N)}\eta^2_Nf)\\
&=\varepsilon^1_{F_2(N)}F_1(\tau_{F_2(N)}G_2F_2(f)\eta^2_M)\\
&=\varepsilon^1_{F_2(N)}F_1(G_1F_2(f)\tau_{F_2(M)}\eta^2_M)\\
&=\varepsilon^1_{F_2(N)}F_1G_1F_2(f)F_1(\tau_{F_2(M)}\eta^2_M)\\
&=F_2(f)\varepsilon^1_{F_2(M)}F_1(\tau_{F_2(M)}\eta^2_M)=F_2(f)\tau^*_M.
\end{split}
\end{equation}
Next, suppose that $\tau:(F_1,G_1)\rightarrow (F_2,G_2)$ and $\sigma:(F_2,G_2)\rightarrow (F_2,G_3)$ are morphisms in $FP(X,Y)$, so that $\tau^*_M=\varepsilon^1_{F_2(M)}F_1(\tau_{F_2(M)}\eta^2_M)$, $\sigma^*_M=\varepsilon^2_{F_3(M)}F_2(\sigma_{F_3(M)}\eta^3_M)$, and $(\tau\sigma)^*_M=\varepsilon^1_{F_3(M)}F_1((\tau\sigma)_{F_3(M)}\eta^3_M)$.
We compute:
\begin{equation}
\begin{split}
\sigma^*_M\tau^*_M&=\sigma^*_M\varepsilon^1_{F_2(M)}F_1(\tau_{F_2(M)}\eta^2_M)\\
&=\varepsilon^1_{F_3(M)}F_1G_1(\sigma^*_M)F_1(\tau_{F_2(M)}\eta^2_M)\\
&=\varepsilon^1_{F_3(M)}F_1(G_1(\sigma^*_M)\tau_{F_2(M)}\eta^2_M)\\
&=\varepsilon^1_{F_3(M)}F_1(\tau_{F_3(M)}G_2(\sigma^*_M)\eta^2_M)\\
&=\varepsilon^1_{F_3(M)}F_1(\tau_{F_3(M)}G_2(\varepsilon^2_{F_3(M)})G_2F_2(\sigma_{F_3(M)}\eta^3_M)\eta^2_M)\\
&=\varepsilon^1_{F_3(M)}F_1(\tau_{F_3(M)}G_2(\varepsilon^2_{F_3(M)})\eta^2_{G_2F_3(M)}\sigma_{F_3(M)}\eta^3_M)\\
&=\varepsilon^1_{F_3(M)}F_1(\tau_{F_3(M)}\sigma_{F_3(M)}\eta^3_M)=(\tau\sigma)^*_M.
\end{split}
\end{equation}
The proof for $(-)^\dagger$ is similar and left to the reader.
\end{proof}

\begin{thm} $(-)^*$ and $(-)^\dagger$ induce a duality between $\Frob(X,Y)$ and $\Frob(Y,X)$.\label{duality thm}\end{thm}

\begin{proof} We show that $(-)^*:FP(X,Y)\rightarrow FP(Y,X)$ and $(-)^\dagger:FP(Y,X)\rightarrow FP(X,Y)$ give a duality between $FP(X,Y)$ and $FP(Y,X)$. The result then follows by the equivalence between $\Frob(X,Y)$ and $FP(X,Y)$ for spaces $X$ and $Y$.

We know from the above calculations that each of $(-)^*$ and $(-)^\dagger$ are contravariant functors. We shall show that $(-)^{*\dagger}=\Id_{FP(X,Y)}$.  

On objects, we have $(F_1,G_1)^{*\dagger}=(G_1,F_1)^\dagger=(F_1,G_1)$.  If $\tau:(F_1,G_1)$ is a morphism, then we have 
\begin{equation}
\begin{split}
(\tau^*)^\dagger_M&=G_1(\varepsilon^2_M\tau^*_{G_2(M)})\eta^1_{G_2(M)}\\
&=G_1(\varepsilon^2_M\varepsilon^1_{F_2G_2(M)}F_1(\tau_{F_2G_2(M)}\eta^2_{G_2(M)}))\eta^1_{G_2(M)}\\
&=G_1(\varepsilon^2_M)G_1(\varepsilon^1_{F_2G_2(M)})G_1F_1(\tau_{F_2G_2(M)}\eta^2_{G_2(M)})\eta^1_{G_2(M)}\\
&=G_1(\varepsilon^2_M)G_1(\varepsilon^1_{F_2G_2(M)})\eta^1_{G_1F_2G_2(M)}\tau_{F_2G_2(M)}\eta^2_{G_2(M)}\\
&=G_1(\varepsilon^2_M)\tau_{F_2G_2(M)}\eta^2_{G_2(M)}\\
&=\tau_MG_2(\varepsilon^2_M)\eta^2_{G_2(M)}=\tau_M.
\end{split}
\end{equation}
Thus $\tau^{*\dagger}=\tau$ and $(-)^{*\dagger}=\Id_{FP(X,Y)}$.  In a similar way one can show that $(-)^{\dagger *}=\Id_{FP(Y,X)}$; the details are left to the reader.
\end{proof}

Having introduced a number of different types of Frobenius bimodules above, it is natural to ask what categorical properties these various bimodules have.  Thus we introduce the following full subcategories of $\Frob(X,Y)$:

\begin{enumerate}
\item[(a)] $\DimPres(X,Y)=\{\M\in\Frob(X,Y):\mbox{$\M$ is dimension preserving}\}$ (assuming $X$ and $Y$ are equipped with dimension functions $d$ and $\d$ respectively).
\item[(b)] $\RLoc(X,Y)=\{\M\in\Frob(X,Y):\mbox{$\M$ is right localizing}\}$.
\item[(c)] $\LLoc(X,Y)=\{\M\in\Frob(X,Y):\mbox{$\M$ is left localizing}\}$.
\end{enumerate}

The following proposition summarizes some of the relevant properties of the above subcategories.

\begin{prop}  Let $X$ and $Y$ be spaces, and consider the functors $(-)^*$ and $(-)^\dagger$ on $\Frob(X,Y)$ defined above.
\begin{enumerate}
\item $(-)^*$ and $(-)^\dagger$ induce a duality between $\RLoc(X,Y)$ and $\LLoc(Y,X)$.
\item If $X$ and $Y$ are equipped with dimension functions $d$ and $\d$ respectively, then $(-)^*$ and $(-)^\dagger$ induce a duality between $\DimPres(X,Y)$ and $\DimPres(Y,X)$.
\item $\DimPres(X,Y)$ is an additive subcategory of $\Frob(X,Y)$. 
\end{enumerate}
\end{prop}

\begin{proof} Part (1) is immediate from the fact that $\M=(F,G)$ is right localizing if and only if $\M^*=(G,F)$ is left localizing.  Similarly, part (2) follows from the observation that $\M$ is dimension preserving if and only if $\M^*$ is. Finally, it is clear that a direct sum of two dimension preserving bimodules is again dimension preserving, proving (3).
\end{proof}

In general, $\RLoc(X,Y)$ will not be an additive subcategory of $\Frob(X,Y)$.  Indeed, let $X$ and $Y$ be noetherian, and let $\M_1=(F_1,G_1)$ and $\M_2=(F_2,G_2)$ be right localizing Frobenius bimodules, with each of $F_1$ and $F_2$ faithful, such that $f_1\neq f_2$, where $f_i:\Inj(X)\rightarrow \Inj(Y)$ are the continuous functions of Theorem \ref{continuous thm}.  Then it is clear that $\M_1\oplus\M_2$ is not right localizing.  

However, suppose that $\M_1$ and $\M_2$ are in addition dimension preserving, and that $\ker G_1$ and $\ker G_2$ are closed under injective envelopes.  Then it follows from Proposition \ref{automatically continuous} that, if $f_1=f_2$, then $\M_1\oplus\M_2$ is a right localizing Frobenius bimodule.   

If $X$ and $Y$ are noetherian spaces and $f:\Inj(X)\rightarrow \Inj(Y)$ is surjective and continuous, then we define $\RLoc_f(X,Y)$ to be the full subcategory consisting of those right localizing bimodules which are faithful and for which $F(E(x))\cong E(f(x))^{(n_x)}$.  Then the above remarks show that $\RLoc_f(X,Y)\cap\DimPres(X,Y)$ is an additive subcategory of $\Frob(X,Y)$.

We next consider properties of Frobenius $X,X$-bimodules.  Note that $\Frob(X)$ is a monoidal category, where $\otimes_X$ is composition of functors.  Since the composition of dimension preserving bimodules is easily seen to be dimension preserving, we see that $\DimPres(X)$ is a monoidal subcategory of $\Frob(X)$.  We introduce two more full subcategories of $\Frob(X)$.

\begin{enumerate}
\item[(d)] $\Loc(X)=\{\M\in\Frob(X):\mbox{$\M$ is localizing}\}$.
\item[(e)] $\LCent(X)=\{\M\in\Frob(X):\mbox{$\M$ is locally centralizing}\}$.
\end{enumerate}

\begin{prop} If $X$ is a space, then each of $\Loc(X)$ and $\LCent(X)$ is an additive monoidal subcategory of $\Frob(X)$,  self-dual under $(-)^*$ and $(-)^\dagger$.
\end{prop}

\begin{proof}  Let $\M_1=(F_1,G_1)$ and $\M_2=(F_2,G_2)$ be Frobenius bimodules on $X$. It is immediate from the definitions that, if $\M_1$ and $\M_2$ are both localizing or both locally centralizing, then $\M_1\otimes_X\M_2=(F_2F_1,G_1G_2)$ is again localizing or locally centralizing. We  show that $\M_1\oplus\M_2$ is localizing if and only if $\M_1$ and $\M_2$ are localizing.  If $T$ is a localizing subcategory of $\Mod X$ and $M\in T$, then $F_1(M)\oplus F_2(M)\in T$ if and only if each of $F_1(M)$ and $F_2(M)$ are in $T$.  A similar claim holds for  $G_1(M)\oplus G_2(M)$.  If $U$ is a weakly open subspace of $\Mod X$, then $F_1|_U\oplus F_2|_U$ commutes with $\Cent(U)$ if and only if each of $F_1|_U$ and $F_2|_U$ do.  Thus $\M_1\oplus\M_2$ is locally centralizing if and only if each of $\M_1$ and $\M_2$ is locally centralizing.

Finally, the duality claim is immediate for localizing bimodules, since $\M$ is localizing if and only if $\M^*$ is localizing, and follows from Lemma \ref{central dual lemma} when $\M$ is locally centralizing.
\end{proof}

\begin{prop}  Let $U$ be a weakly open subspace of a space $X$.  Then restriction to $U$ as defined in section 9 defines a functor $(-)|_U:\Loc(X)\rightarrow \Loc(U)$.
\end{prop}

\begin{proof}  Write $j:U\rightarrow X$ for the inclusion.  We have already shown in Lemma \ref{localizing restriction lemma} that if $\M\in\Loc(X)$, then $\M|_U\in\Loc(U)$.  We define $(-)|_U$ on morphisms as follows.  Recall that a morphism $\tau:\M_1\rightarrow \M_2$ is a natural transformation $\tau:G_2\rightarrow G_1$, where $\M_1=(F_1,G_1)$ and $\M_2=(F_2,G_2)$.  Then we define a natural transformation $\tau|_U:G_2|_U\rightarrow G_1|_U$ by $(\tau|_U)_M=j^*(\tau_{j_*M})$ for all $M\in\Mod U$.

We show that $\tau|_U$ is in fact a natural transformation.  Given $M,N\in\Mod U$ and $f\in\Hom_U(M,N)$, there exists $\tilde f\in\Hom_X(j_*M,j_*N)$ with $j^*(\tilde f)=f$.  Thus we have
\begin{multline}(\tau|_U)_N G_2|_U(f)=(\tau|_U)_N G_2|_Uj^*(\tilde f)=(\tau|_U)_N j^*G_2(\tilde f)=j^*(\tau_{j_*M}G_2(\tilde f))\\
=j^*(G_1(\tilde f)\tau_{j_*N})=j^*G_1(\tilde f)j^*(\tau_{j_*N})=G_1|_U(f)(\tau|_U)_N.
\end{multline}
To see that $(-)|_U$ is a functor, we need to show that $(\tau\sigma)|_U=\tau|_U\sigma|_U$ for all $\sigma:\M_1\rightarrow\M_2$ and $\tau:\M_2\rightarrow\M_3$.  Given $M\in\Mod U$, we have by definition \[(\tau\sigma|_U)_M=j^*(\tau\sigma_{j_*M})=j^*(\tau_{j_*M}\sigma_{j_*M})=j^*(\tau_{j_*M})j^*(\sigma_{j_*M})=(\tau|_U)_M(\sigma|_U)_M.\]
Thus $(-)|_U$ is a functor.
\end{proof}

\section{Noncommutative vector bundles?}
As we mentioned in the introduction, one of the main motivations for studying Frobenius bimodules was the desire to formulate a general definition of ``vector bundle" on a noncommutative space $X$.  In this final section we discuss some possible definitions, weighing the various pros and cons that they offer.  Before doing so, we give a list of desirable properties that noncommutative vector bundles ``should" have.  To fix notation, let $X$ be a space and let $\E$ be a (as yet undefined) vector bundle on $X$.  The category of vector bundles on $X$ will be denoted by $\Vect(X)$.

\begin{enumerate}
\item $\E$ should be two-sided; that is, $\E$ should be an $X,X$-bimodule.  

\item It should be possible to define the dual bundle $\E^*$ to a vector bundle $\E$, and $\E$ should be reflexive: $\E^{**}\cong\E$.

\item It should be possible to study $\E$ locally; that is, given a weakly open (or perhaps open) subspace $U$ of $X$, there should be a way to define the restriction $\E|_U$ of $\E$.  Also, if $\{U_i:i\in I\}$ is a (weakly) open cover of $X$, it should be possible to patch together information from the various $\E|_{U_i}$ to obtain information about $\E$.
\item Vector bundles should exhibit \emph{noncommutative} phenomena, even over commutative spaces.  
 \item Given a map of spaces $f:Z\rightarrow X$, it should be possible to pull $\E$ back to a vector bundle $f^*\E$ on $Z$, so that the projection formula holds: $f_*(M\otimes_Zf^*\E)\cong f_*M\otimes_X\E$ for $M\in\Mod Z$.
\item It should be possible to carry out many, if not all, of the standard constructions in algebraic geometry such as tensor algebras, symmetric algebras, exterior algebras, etc. of a vector bundle $\E$.

\end{enumerate}

Our contention is that it is not unreasonable to look for $\Vect(X)$ as a full subcategory of $\Frob(X)$; that is, vector bundles on $X$ should be sought among the Frobenius $X,X$-bimodules.  Indeed, conditions (1) and (2) on the above list become automatic for Frobenius bimodules.  We have also seen that for certain (but not all) full subcategories of $\Frob(X)$, conditions (3) and (4) also hold.  We discuss here three potential definitions of vector bundle made in terms of  Frobenius bimodules, and discuss advantages and disadvantages to all three possibilities.

\subsection{Frobenius $X,X$-bimodules} The most obvious definition is to declare a vector bundle to be a Frobenius bimodule; that is, $\Vect(X)=\Frob(X)$.  As already noted, conditions (1) and (2) are immediate, and  since $\Frob(X)$ is a monoidal category, it is possible to construct the tensor algebra of a Frobenius bimodule, addressing one of the elements of condition (6).  Finally, Frobenius bimodules exhibit noncommutative phenomena, addressing condition (4).  For example, if $\E$ is an ordinary vector bundle on a scheme $X$ and $\vphi$ is an automorphism of $X$, then the twisted sheaf ${_1\E_\vphi}$ is a Frobenius bimodule.    

The main drawbacks to adopting this definition of vector bundle are conditions (3) and (5).  If $\M=(F,G)$ is a Frobenius bimodule over $X$ and $U$ is a weakly open subspace of $X$, then $\M|_U$ is an $F^{-1}U,U$-bimodule, but need not be Frobenius in general.  As we have seen above, the theory of Frobenius bimodules becomes much richer in situations where we know that $\M|_U$ is again Frobenius. It is unclear how suitable general Frobenius bimodules are to local analysis.

This leaves condition (5) and, as the following example shows, the desire to be able to pull back is in some sense incompatible with the desire for noncommutative phenomena articulated in (4).

\begin{example} Let $R$ be a commutative ring, let $\vphi$ be an automorphism of $R$, and let $I$ be an ideal of $R$ with $\vphi^{-1}(I)\neq I$.  We let $X=\Mod R$ and $Z=\Mod R/I$, and $\E={_1R_\vphi}$ (i.e. we twist the action of $R$ on the right by $\vphi$).  Note that $\E$ is a Frobenius bimodule; indeed, given $M\in\Mod R$, then $M\otimes_R\E\cong M_\vphi$; that is, tensoring with $\E$ twists the action of $R$ on $M$ by $\vphi$.  If $i:Z\rightarrow X$ denotes the inclusion, then we claim that there is no Frobenius bimodule $i^*\E$ on $Z$  such that the projection formula holds.  

To see this, consider $M=R/I$.  Since the functor $i_*$ just views an $R/I$-module as an $R$-module, we see that $i_*(R/I)\otimes_R\E\cong (R/I)_\vphi\cong R/\vphi^{-1}(I)$.  On the other hand, $i_*(R/I\otimes_{R/I}i^*\E)$ is necessarily annihilated by $I$. Since $\vphi^{-1}(I)\neq I$ by construction, we see that $i^*\E$ cannot exist.\qed\label{can't pull back ex}
\end{example}

The above example is in some sense a major disappointment, because if we want to define vector bundles as bimodules over $X$, anyone would agree that the autoequivalences of $\Mod X$ should be included in this definition.  We shall see below that a more restrictive definition of vector bundle would preclude this behavior.  

It seems worthwhile to point out that it is too much to expect the projection formula in (5) to hold for all maps of spaces, because the notion of map of space is too general.   For example, let $Z$ be a weakly closed subspace of a space $X$.  Then there is a map of spaces $f:X\rightarrow Z$, where $f_*=i^!$ and $f^*=i_*$. (Here $i^!$ is a right adjoint to the inclusion $i_*$ and is called the \emph{support functor}.)  The projection formula for this map of spaces would then read $i^!(M\otimes_Xi_*\E)\cong i^!M\otimes_Z\E$, where $\E$ is a vector bundle on $Z$.  We then present the following explicit example, which shows that one cannot hope to pull back $\E$ in this setting.

\begin{example} Let $X=\Mod R$, where $R$ is a commutative polynomial ring in $\geq 2$ variables, and let $Z=\Mod R/I$, where $I$ is chosen so that $R/I\cong R_1\oplus R_2$ for rings $R_1,R_2$. Then $\E=R_1$ is a projective $R/I$-module and hence a (commutative) vector bundle on $\Mod R/I$.  We claim that $i_*\E$ cannot exist; indeed, if it did, then it would be a projective, hence free, module over $R$.  If we consider $R_2$ as an $R$-module, then $i^!(R_2\otimes_R i_*\E)\cong i^!R_2^{(n)}$, where $n$ is the rank of $i_*\E$.  But $i^!R_2\otimes_{R/I}\E=0$, since $\E$ annihilates $R_2$.  \qed
\end{example}

The difficulty in this example is that, while $f:\Mod R\rightarrow \Mod R/I$ is a map of spaces, it is not a ``geometric" map; i.e. it is not induced by a ring homomorphism.  Thus it seems that the definition of map of noncommutative spaces may need to be modified, or that any pullback formula can only be expected to hold for a restricted class of maps between spaces.

\subsection{Right or left localizing Frobenius $X,X$-bimodules}  In order to have a definition of noncommutative vector bundle that enables local study as in condition (3), we might look to left and/or right localizing Frobenius bimodules for our definition.  In this case it seems advantageous to distinguish between the two, and speak of \emph{left} and \emph{right} vector bundles over $X$, and define a vector bundle to be both a left and right vector bundle: $\Vect(X)=\RLoc(X)\cap\LLoc(X)$.

While this definition makes the earlier results of the paper available for studying vector bundles locally, it has some drawbacks with regard to the other conditions.  For instance $\RLoc(X)$ and $\LLoc(X)$ are not additive categories.  One possible remedy is to restrict attention to the case of bimodules $\M=(F,G)$ with $F$ and $G$ faithful, and study the categories $\RLoc_f(X)$ and/or $\LLoc_f(X)$ as defined above.  While these categories are additive, they restrict the study of vector bundles to one map $f$ at a time.  Also, $\RLoc_f(X)$ and $\LLoc_f(X)$ are not monoidal categories: if $\M\in\RLoc_f(X)$ and $\sh{N}\in\RLoc_g(X)$, then $\M\otimes_X\sh{N}\in\RLoc_{gf}(X)$.  Another possible remedy is to define a left vector bundle to be a finite direct sum of left localizing bimodules, and similar for a right vector bundle.  This definition gives additive monoidal categories, but would presumably make local analysis more difficult.  

\subsection{Localizing Frobenius $X,X$-bimodules} Finally, we could define a vector bundle to be a localizing Frobenius bimodule; that is $\Vect(X)=\Loc(X)$.  Since $\Loc(X)$ is an additive monoidal category, we have the minimum requirement needed to carry out the constructions in condition (6).  Also, it is possible to study such vector bundles locally as in condition (3), and in fact one can always restrict a localizing bimodule over $X$ to a localizing bimodule over $U$ for any weakly open subspace $U$ of $X$, using the restriction defined in section 10. 

The main drawback to this definition is that it may not be sufficiently ``noncommutative" as desired in condition (4).  For instance, many nontrivial autoequivalences of $\Mod X$ would not be vector bundles under this definition.  In particular, if $\vphi$ is an automorphism of a commutative ring $R$ which induces a nontrivial permutation of the prime ideals of $R$, then $_1R_\vphi$ would not be a vector bundle under this definition.

While this admittedly is a strike against adopting the definition $\Vect(X)=\Loc(X)$, we do point out one positive observation:  The fact that $_1R_\vphi$ is not a vector bundle under this definition means that Example \ref{can't pull back ex} is \emph{not} a counterexample to pulling back.  We do not know if it is possible to pull back a localizing bimodule for a suitable class of maps of spaces as in condition (5).

\appendix
\section{Basic results and definitions in abelian categories}
We collect here in an appendix several basic definitions and results concerning abelian categories.  All categories considered will be complete and cocomplete.  We begin by recalling some finiteness conditions. (See \cite{Mitchell} or \cite{Popescu} for further details.)

\begin{defn} Let $\A$ be an abelian category.  An object $M$ of $\A$ is called \emph{finitely presented} if $\Hom_\A(M,-)$ commutes with arbitrary direct limits.  Dually, we say that $M$ is \emph{finitely copresented} if it is finitely presented in the opposite category $\A^{\rm op}$. Writing this condition out explicitly, $M$ is finitely copresented if $\Hom_\A(\varprojlim N_i,M)\cong \varinjlim\Hom_\A(N_i,M)$.  

Finally, we call $M$ \emph{finitely generated} if, whenever we have $\sum_{i\in I}M_i=M$ for submodules $M_i$ of $M$, there exists a finite $J\subseteq I$ with $\sum_{i\in J}M_i=M$. Dually, $M$ is \emph{finitely cogenerated} if $M$ is finitely generated in $\A^{\rm op}$; thus $M$ is finitely cogenerated if and only if whenever $\bigcap_{i\in I}M_i=0$ for submodules $M_i$ of $M$, we have $\bigcap_{i\in J}M_i=0$ for some finite $J\subseteq I$.
\end{defn}

\begin{lemma}  Any noetherian $\A$-module is finitely presented. Dually, any artinian $\A$-module is finitely copresented.
\end{lemma}

\begin{proof} This is Exercise 1 on p. 370 of \cite{Popescu}.\end{proof}

\begin{lemma} If $M$ is finitely copresented then $M$ is finitely cogenerated.  Dually if $M$ is finitely presented then it is finitely generated.
\end{lemma}

\begin{proof}  Let $\{M_i:i\in I\}$ be a collection of submodules of $M$ with $\bigcap_{i\in I}M_i=0$.  If $J$ is a finite subset of $I$, then we define $M_J$ to be $\bigcap_{i\in J} M_i$.  Now, for every finite $J\subseteq I$ we have an exact sequence
\begin{equation}
0\rightarrow M_J\rightarrow M\rightarrow \prod_{i\in J}M/M_i.\end{equation}
Each of the modules fits into a natural inverse system over the collection of finite subsets of $I$, the one for $M$ simply being the constant system.  Since $\varprojlim$ is left exact and $\varprojlim \prod_{i\in J}M/M_i=\prod_{i\in I}M/M_i$, we obtain the exact sequence
\begin{equation}
0\rightarrow \varprojlim M_J\rightarrow M\rightarrow \prod_{i\in I}M/M_i.\end{equation}
It follows that $\varprojlim M_J\cong \bigcap_{i\in I}M_i=0$.  Since $M$ is finitely copresented, we have $\varinjlim\Hom_\A(M_J,M)\cong \Hom_\A(\varprojlim M_J,M)=0$.  Thus, given any morphism $\vphi:M_J\rightarrow M$, there exists a $J_0\supseteq J$ such that the restriction of $\vphi$ to $M_{J_0}$ is $0$.  This applies in particular to the natural inclusion functor $M_J\rightarrow M$; that is, the inclusion $M_{J_0}\rightarrow M$ is $0$ for some $J_0\supseteq J$.  This says precisely that $\bigcap_{i\in J_0}M_i=0$, and $M$ is finitely cogenerated.

The dual statement follows by dualizing the above proof.
\end{proof}

Finally, suppose that $\A$ is a Grothendieck category; that is, $\A$ is an Ab5 abelian category with a generator.  Then we have the following.

\begin{prop} If $\A$ is a Grothendieck category, then there exist simple $\A$-modules.
\end{prop}

\begin{proof} Let $G$ be a generator for $\A$.  The Gabriel-Popescu Theorem \cite{Gab-Pop} states that the functor $\Hom_\A(G,-)$ induces an equivalence between $\A$ and a quotient category of $\Mod R$, where $R$ is the ring $\End_\A(G)$.  Now, let $\Sigma=\{M_i:i\in I\}$ be the lattice of all proper subobjects of $G$.  Then since $\Hom_\A(G,-)$ is left exact, the image of $\Sigma$ in $\Mod R$ is a sublattice of the lattice of proper right ideals of $R$.  This shows that if $\{G_j:j\in J\}$ is a chain of proper subobjects of $G$, then $\bigcup_{j\in J}G_j$ is a proper subobject of $G$.  Zorn's Lemma can then be applied to show that $\Sigma$ has maximal elements.  If $M$ is a maximal subobject of $G$, then $G/M$ is a simple $\A$-module.
\end{proof}

\end{document}